\documentclass[a4paper, 12pt]{article}
\usepackage[T2A]{fontenc}
\usepackage[english,russian]{babel}
\usepackage{amsmath}
\usepackage{amsfonts,amssymb,mathrsfs,amscd}
\usepackage{bm}

\newcommand{\Kappa}{\operatorname{\mathcal K}}
\newcommand{\sgn}{\operatorname{sign}}
\newcommand{\loc}{\operatorname{loc}}
\newcommand{\inter}{\operatorname{int}}
\newcommand{\e}{\mathfrak e}
\renewcommand{\l}{\mathfrak l}

\newcommand{\mn}{\operatorname{\bm m}}
\newcommand{\m}{\mathfrak m}
\newcommand{\n}{\mathfrak n}
\newcommand{\M}{\mathfrak M}
\newcommand{\s}{\bm s}

\newcommand{\card}{\operatorname{card}}

\newcommand{\supp}{\operatorname{supp}}
\newcommand{\supvrai}{\operatornamewithlimits{sup\,vrai}}
\newcommand{\N}{\mathbb N}
\newcommand{\Z}{\mathbb Z}
\newcommand{\R}{\mathbb R}
\newcommand{\Nu}{\mathcal N}

\newcommand{\D}{\mathcal D}
\newcommand{\J}{\mathcal J}

\renewcommand{\L}{\mathcal L}
\newcommand{\mes}{\operatorname{mes}}
\allowdisplaybreaks

\begin{document}

\author{ С. Н. Кудрявцев }
\title{Продолжение функций из пространств
Никольского -- Бесова смешанной гладкости \\
за пределы областей определённого вида}
\date{}
\maketitle
\begin{abstract}
В статье рассмотрены пространства Никольского и Бесова
с нормами, в определении которых вместо смешанных модулей непрерывности
известных порядков определённых смешанных производных функций
используются "$L_p$-усреднённые" смешанные модули непрерывности функций
соответствующих порядков. Для таких пространств функций, заданных в
областях определённого вида. построены непрерывные линейные
отображения их в обычные пространства Никольского и Бесова смешанной гладкости
в $ \R^d, $ являющиеся операторами продолжения функций, что влечёт совпадение тех и
других пространств, заданных в таких областях. При этом существенно расширен
класс пространств Никольского -- Бесова функций смешанной гладкости, для которых
получены теоремы о продолжении такого рода. При определённых условиях
установлена непрерывность операторов частного дифференцирования, действующих
из рассматриваемых пространств функций смешанной гладкости в пространства
Лебега.
\end{abstract}

Ключевые слова: пространства Никольского -- Бесова смешанной гладкости,
продолжение функций, эквивалентные нормы
\bigskip

\centerline{Введение}
\bigskip

В работе рассмотрены некоторые задачи теории функциональных
пространств для пространств Никольского и Бесова функций смешанной
гладкости, заданных в областях определённого вида.
Остановимся подробнее на содержании работы.

При $ d \in \N, \alpha \in \R_+^d, 1 \le p,\theta < \infty $ для области $ D
\subset \R^d $ вводятся в рассмотрение
пространства функций $ (S_{p,\theta}^\alpha B)^\prime(D) ((S_p^\alpha H)^\prime(D)) $
с нормами
\begin{multline*}
\| f \|_{(S_{p,\theta}^\alpha B)^\prime(D)} = \\
\max \biggl(\| f \|_{L_p(D)},
\max_{J \subset \{1,\ldots,d\}: J \ne \emptyset}
\left(\int_{(0,\infty) \times \ldots \times (0,\infty)}
(\prod_{j \in J} t_j^{-1 -\theta \alpha_j})
(\Omega^{\prime l \chi_J}(f, t^J)_{L_p(D)})^{\theta} \,dt^J
\right)^{1/\theta}\biggr),
\end{multline*}
$$
\| f \|_{(S_p^\alpha H)^\prime(D)} = \max(\| f \|_{L_p(D)},
\max_{J \subset \{1, \ldots, d\}: J \ne \emptyset}
\sup_{t^J: t_j > 0, j \in J} (\prod_{j \in J} t_j^{-\alpha_j})
\Omega^{\prime l \chi_J}(f, t^J)_{L_p(D)}),
$$
где
\begin{multline*}
\Omega^{\prime l \chi_J}(f, t^J)_{L_p(D)} =
\biggl((\prod_{j \in J} (2 t_j)^{-1}) \int_{ \{\xi^J: |\xi_j| \le t_j, j \in J\}}
\| (\prod_{j \in J} \Delta_{\xi_j e_j}^{l_j}) f\|_{L_p(D_{\xi}^{l \chi_J})}^p d\xi^J\biggr)^{1 /p}, \\
l \in \N^d - \text{ вектор с компонентами } \\
l_j = \min \{m \in \N: \alpha_j < m \}, j =1,\ldots,d, \\
\text{ а для } J = \{j_1, \ldots, j_k\}: 1 \le j_1 < \ldots < j_k \le d, \text{ и }
x \in \R^d \text{ вектор } x^J = (x_{j_1},\ldots,x_{j_k}).
\end{multline*}

Для областей $ D \subset \R^d, $ имеющих определённое строение,
построены непрерывные линейные отображения пространств $
(S_{p,\theta}^\alpha B)^\prime(D) ((S_p^\alpha H)^\prime(D)) $ в
обычные пространства Бесова $ S_{p,\theta}^\alpha B(\R^d) $
(Никольского $ S_p^\alpha H(\R^d)), $ являющиеся операторами
продолжения функций, что влечет совпадение тех и других
пространств в упомянутых областях (см. п. 2.3.). Из публикаций, в
которых рассматривается подобная задача для таких пространств,
автору известны лишь его собственные работы [1] и [2], в которых в
качестве $ D $ рассматриваются куб $ I^d $ и области, так
называемого, $m$-рода. Среди близких работ по этой тематике
отметим [3] -- [12] (см. также имеющуюся там литературу), в
которых изучается вопрос о продолжении за пределы области
определения с сохранением класса гладких функций из изотропных и
неизотропных пространств, а также функций из пространств смешанной
гладкости. Отметим, что средства и схема построения операторов
продолжения функций и вывод метрических соотношений, применяемых
ниже для доказательства непрерывности таких операторов, отличаются
от тех, что использовались в упомянутых работах. Добавим ещё, что
нормы, определяющие пространства функций смешанной гладкости,
рассматриваемые в [12], отличаются от $ \| \cdot
\|_{(S_{p,\theta}^\alpha B)^\prime(D)}, \| \cdot \|_{(S_p^\alpha
H)^\prime(D)}. $ Кроме того, теорема 2 из [12] о продолжении
функций смешанной гладкости из рассматриваемых там пространств
установлена лишь при $ 1 < p \le \theta < \infty, $ а класс
областей определения функций, удовлетворяющих условиям теоремы 2
из [12], содержится в классе областей, удовлетворяющих условиям
теоремы 2.3.1 (см. п. 2.3.).

Статья состоит из введения и двух параграфов, в первом из которых
рассматриваются некоторые вспомогательные средства для решения
объявленной задачи, а во втором -- средства построения и конструкция
подходящих операторов продолжения, а также доказательство их свойств.
\bigskip

\centerline{\S 1. Предварительные сведения и вспомогательные утверждения}
\bigskip

1.1. В этом пункте вводятся обозначения, относящиеся к
функциональным пространствам, рассматриваемым в
настоящей работе, а также приводятся некоторые факты, необходимые
в дальнейшем.

Для $ d \in \N $ через $ \Z_+^d $ обозначим множество
$$
\Z_+^d = \{\lambda = (\lambda_1, \ldots, \lambda_d) \in \Z^d:
\lambda_j \ge0, j=1, \ldots, d\}.
$$
Обозначим также при  $ d \in \N $ для $ l \in \Z_+^d $ через $ \Z_+^d(l) $
множество
\begin{equation*} \tag{1.1.1}
\Z_+^d(l) = \{ \lambda  \in \Z_+^d: \lambda_j \le l_j, j=1, \ldots, d\}.
\end{equation*}

Для $  d \in  \N, l \in \Z_+^d $ через $ \mathcal P^{d,l} $ будем
обозначать пространство вещественных  полиномов, состоящее из всех
функций $ f: \R^d \mapsto \R $ вида
$$
f(x) = \sum_{\lambda \in \Z_+^d(l)} a_{\lambda} \cdot x^{\lambda}, x \in \R^d,
$$
где $ a_{\lambda} \in \R, x^{\lambda} = x_1^{\lambda_1} \ldots x_d^{\lambda_d},
\lambda \in \Z_+^d(l). $

Для множества $ A $ из топологического пространства $ T $ через $ \overline A $
обозначается замыкание множества $ A, $ а через $ \inter A $  -- его внутренность.
В $ \R^d $ рассматривается топология, порождённая нормой
$$
\|x\| = \max_{j =1,\ldots,d} |x_j|, x = (x_1,\ldots,x_d) \in \R^d.
$$

Для измеримого по Лебегу множества $ D \subset \R^d $ при $ 1 \le p \le \infty $
через $ L_p(D), $ как обычно, обозначается
пространство  всех  вещественных измеримых на $ D $ функций $f,$
для которых определена норма
$$
\|f\|_{L_p(D)} = \begin{cases} (\int_D |f(x)|^p dx)^{1/p}, 1 \le p < \infty; \\
\supvrai_{x \in D}|f(x)|, p = \infty. \end{cases}
$$

Отметим здесь неоднократно используемое в дальнейшем неравенство
\begin{equation*} \tag{1.1.2}
| \sum_{j =1}^n x_j|^a \le \sum_{j =1}^n |x_j|^a, x_j \in \R, j =1,\ldots,n,
n \in \N, 0 \le a \le 1.
\end{equation*}

Для $ x,y \in \R^d $ положим $ xy = x \cdot y = (x_1 y_1, \ldots, x_d y_d), $
а для $ x \in \R^d $ и $ A \subset \R^d $ определим
$$
x A = x \cdot A = \{xy: y \in A\}.
$$
Для $ A, B \subset \R^d $  положим ещё
$$
A +B = \{ x +y: x \in A, y \in B\}.
$$

Для $ x \in \R^d: x_j \ne 0, $ при $ j=1,\ldots,d,$ положим
$ x^{-1} = (x_1^{-1},\ldots,x_d^{-1}). $

При $ d \in \N $ для $ x,y \in \R^d $ будем писать $ x \le y (x < y), $ если
для каждого $ j=1,\ldots,d $ выполняется неравенство $ x_j \le y_j (x_j < y_j). $

При $ d \in \N $ для $ x \in \R^d $ положим
$$
x_+ = ((x_1)_+, \ldots, (x_d)_+),
$$
где $ t_+ = \frac{1} {2} (t +|t|), t \in \R. $

Обозначим через $ \R_+^d $ множество $ x \in \R^d, $ для которых
$ x_j >0 $ при $ j=1,\ldots,d,$ и для $ a \in \R_+^d, x \in \R^d $
положим $ a^x = a_1^{x_1} \ldots a_d^{x_d}.$

При $ d \in \N $ определим множества
$$
I^d = \{x \in \R^d: 0 < x_j < 1,j=1,\ldots,d\},
$$
$$
\overline I^d = \{x \in \R^d: 0 \le x_j \le 1,j=1,\ldots,d\},
$$
$$
B^d = \{x \in \R^d: -1 \le x_j \le 1,j=1,\ldots,d\}.
$$

Через $ \e $ будем обозначать вектор в $ \R^d, $ задаваемый
равенством $ \e = (1,\ldots,1). $

При $ d \in \N $ для $ \lambda \in \Z_+^d $ через $ \D^\lambda $
будем обозначать оператор дифференцирования $ \D^\lambda =
\frac{\D^{|\lambda|}} {\D x_1^{\lambda_1} \ldots \D x_d^{\lambda_d}}, $
где $ |\lambda| = \sum_{j=1}^d \lambda_j. $

Теперь приведём некоторые факты, относящиеся к полиномам, которыми
мы будем пользоваться ниже.

В [13] содержится такое утверждение.

Лемма 1.1.1

Пусть $ d \in \N, l \in \Z_+^d, \lambda \in \Z_+^d, 1 \le p,q \le \infty,
\rho, \sigma \in \R_+^d. $ Тогда существует константа
$ c_1(d,l,\lambda,\rho, \sigma) >0 $ такая, что для любых
измеримых по Лебегу множеств $ D,Q \subset \R^d, $  для которых
можно найти $ \delta \in \R_+^d $ и $ x^0 \in \R^d $ такие, что
$ D \subset (x^0 +\rho \delta B^d) $ и $ (x^0 +\sigma \delta I^d) \subset Q, $
для любого полинома $ f \in \mathcal P^{d,l} $
выполняется неравенство
    \begin{equation*} \tag{1.1.3}
\| \D^\lambda f\|_{L_q(D)} \le c_1 \delta^{-\lambda -p^{-1} \e +q^{-1} \e}
\|f\|_{L_p(q)}.
   \end{equation*}

Далее, напомним, что для открытого множеста $ D \subset \R^d $ и вектора
$ h \in \R^d $ через $ D_h $ обозначается множество
$$
D_h = \{x \in D: x +th \in D \ \forall t \in \overline I\}.
$$

Для функции $ f, $ заданной на открытом множестве $ D \subset \R^d, $ и
вектора $ h \in \R^d $ определим на $ D_h $ её разность $ \Delta_h f $
с шагом $ h, $ полагая
$$
(\Delta_h f)(x) = f(x+h) -f(x), x \in D_h,
$$
а для $ l \in \N: l \ge 2, $ на $ D_{lh} $ определим $l$-ую
разность $ \Delta_h^l f $ функции $ f $ с шагом $ h $ равенством
$$
(\Delta_h^l f)(x) = (\Delta_h (\Delta_h^{l-1} f))(x), x \in
D_{lh},
$$
положим также $ \Delta_h^0 f = f. $

Как известно, справедливо равенство
$$
(\Delta_h^l f)(\cdot) = \sum_{k=0}^l C_l^k (-1)^{l-k} f(\cdot +kh),
C_l^k =\frac{l!} {k! (l-k)!}.
$$

При $ d \in \N $ для $ j=1,\ldots,d$ через $ e_j $ будем
обозначать вектор $ e_j = (0,\ldots,0,1_j,0,\ldots,0).$

Как показано в [14], справедлива следующая лемма.

Лемма 1.1.2

Пусть $ d \in \N, l \in \Z_+^d. $ Тогда
для любых $ \delta \in \R_+^d $ и $ x^0 \in \R^d $ для $ Q = x^0 +\delta I^d $
существует единственный линейный оператор
$ P_{\delta, x^0}^{d,l}: L_1(Q) \mapsto \mathcal P^{d,l}, $
обладающий следующими свойствами:

1) для $ f \in \mathcal P^{d,l} $ имеет место равенство
\begin{equation*} \tag{1.1.4}
P_{\delta, x^0}^{d,l}(f \mid_Q) = f,
  \end{equation*}

2)
\begin{equation*}
\ker P_{\delta,x^0}^{d,l} = \biggl\{\,f \in L_1(Q):
\int \limits_{Q} f(x) g(x) \,dx =0\ \forall g \in \mathcal P^{d,l}\,\biggr\},
\end{equation*}

причём существуют константы $ c_2(d,l) >0 $ и $ c_3(d,l) >0 $
такие, что

   3) при $ 1 \le p \le \infty $ для $ f \in L_p(Q) $ справедливо неравенство
  \begin{equation*} \tag{1.1.5}
\|P_{\delta, x^0}^{d,l} f \|_{L_p(Q)} \le c_2 \|f\|_{L_p(Q)},
  \end{equation*}

4) при $ 1 \le p < \infty $ для $ f \in L_p(Q) $ выполняется неравенство
   \begin{equation*} \tag{1.1.6}
  \| f -P_{\delta, x^0}^{d,l} f \|_{L_p(Q)} \le c_3 \sum_{j=1}^d
\delta_j^{-1/p} \biggl(\int_{\delta_j B^1} \int_{Q_{(l_j +1) \xi e_j}}
|\Delta_{\xi e_j}^{l_j +1} f(x)|^p dx d\xi\biggr)^{1/p}.
\end{equation*}

Теперь определим пространства функций, изучаемые в настоящей работе (ср. с
[15], [16]). Но прежде введём некоторые обозначения.

При $ d \in \N $ для $ x \in \R^d $ через $\s(x) $ обозначим
множество $ \s(x) = \{j =1,\ldots,d: x_j \ne 0\}, $ а для
множества $ J \subset \{1,\ldots,d\} $ через $ \chi_J $ обозначим
вектор из $ \R^d $ с компонентами
$$
(\chi_J)_j = \begin{cases} 1, & \text{ для } j \in J; \\
0, & \text{ для } j \in (\{1,\ldots,d\} \setminus J). \end{cases}
$$

При $ d \in \N $ для $ x \in \R^d $ и $ J = \{j_1,\ldots,j_k\}
\subset \N: 1 \le j_1 < j_2 < \ldots < j_k \le d, $ через $ x^J $
обозначим вектор $ x^J = (x_{j_1},\ldots,x_{j_k}) \in \R^k, $ а
для множества $ A \subset \R^d $ положим $ A^J = \{x^J: x \in A\}. $

Для открытого множества $ D \subset \R^d $ и векторов $ h \in \R^d $ и $ l \in
\Z_+^d $ через $ D_h^l $ обозначим множество
\begin{multline*}
D_h^l = (\ldots (D_{l_d h_d e_d})_{l_{d-1} h_{d-1} e_{d-1}}
\ldots)_{l_1 h_1 e_1} = \{ x \in D: x +tlh \in D \ \forall t \in
\overline I^d\} = \\ \{ x \in D: (x +\sum_{j \in \s(l)} t_j l_j h_j e_j) \in D \
\forall t^{\s(l)} \in (\overline I^d)^{\s(l)} \}.
\end{multline*}

Пусть $ d \in \N, D $ -- открытое множество в $ \R^d $ и $ 1 \le p \le \infty. $
Тогда для $ f \in L_p(D), h \in \R^d $ и $ l \in \Z_+^d $ определим в $ D_h^l $
смешанную разность функции $ f $ порядка $ l, $ соответствующую вектору $ h, $
равенством
\begin{multline*}
(\Delta_h^l f)(x) = \biggl(\biggl(\prod_{j=1}^d \Delta_{h_j e_j}^{l_j}\biggr) f\biggr)(x)
= \biggl(\biggl(\prod_{j \in \s(l)} \Delta_{h_j e_j}^{l_j}\biggr) f\biggr)(x) = \\
\sum_{k \in \Z_+^d(l)} (-\e)^{l-k} C_l^k f(x+kh), x \in D_h^l,
\end{multline*}
где $ C_l^k = \prod_{j=1}^d C_{l_j}^{k_j}. $

Имея в виду, что для $ f \in L_p(D), l \in \Z_+^d $ и векторов
$ h,h^\prime \in \R^d: h^{\s(l)} = (h^\prime)^{\s(l)}, $ соблюдается
соотношение
$$
\| \Delta_h^l f\|_{L_p(D_h^l)} = \| \Delta_{h^\prime}^l
f\|_{L_p(D_{h^\prime}^l)}, 1 \le p \le \infty,
$$
определим при $ 1 \le p \le \infty $ для функции $ f \in L_p(D) $ смешанный
модуль непрерывности в $ L_p(D) $ порядка $ l \in \Z_+^d $ равенством
$$
\Omega^l (f,t^{\s(l)})_{L_p(D)} = \supvrai_{ \{ h \in \R^d:
h^{\s(l)} \in t^{\s(l)} (B^d)^{\s(l)} \}} \| \Delta_h^l f\|_{L_p(D_h^l)},
t^{\s(l)} \in (\R_+^d)^{\s(l)}.
$$
Кроме того, при тех же условиях введём в рассмотрение для функции $ f $
"усреднённый" смешанный модуль непрерывности в $ L_p(D) $ порядка $ l, $
полагая
\begin{multline*}
\Omega^{\prime l} (f, t^{\s(l)})_{L_p(D)} = \begin{cases}
\biggl((2 t^{\s(l)})^{-\e^{\s(l)}}
\int_{ t^{\s(l)} (B^d)^{\s(l)}} \| \Delta_\xi^l f\|_{L_p(D_\xi^l)}^p
d \xi^{\s(l)}\biggr)^{1 /p} = \\
\biggl((2 t^{\s(l)})^{-\e^{\s(l)}}
\int_{ (t B^d)^{\s(l)}} \int_{D_\xi^{l \chi_{\s(l)}}}
| \Delta_\xi^{l \chi_{\s(l)}} f(x)|^p dx
d \xi^{\s(l)}\biggr)^{1 /p}, p \ne \infty; \\
\Omega^l (f,t^{\s(l)})_{L_p(D)}, p = \infty,
\end{cases} \\
t^{\s(l)} \in (\R_+^d)^{\s(l)}.
\end{multline*}

Из приведенных определений видно, что
\begin{multline*} \tag{1.1.7}
\Omega^{\prime l} (f, t^{\s(l)})_{L_p(D)} \le
\Omega^l (f, t^{\s(l)})_{L_p(D)}, t^{\s(l)} \in (\R_+^d)^{\s(l)}, \\
f \in L_p(D), 1 \le p \le \infty, l \in \Z_+^d, \\
D \text{ -- произвольное открытое множество в } \R^d.
\end{multline*}

Пусть теперь $ d \in \N, \alpha \in \R_+^d, 1 \le p \le \infty, D $ --
область в $ \R^d $ и вектор $ \l \in \Z_+^d $ такой, что $ \l < \alpha. $
Тогда зададим вектор $ l = l(\alpha) \in \N^d, $ полагая
$ l_j = (l(\alpha))_j = \min \{m \in \N: \alpha_j < m \}, j =1,\ldots,d, $
и через $ (S_p^\alpha H)^{\l}(D) $ обозначим пространство всех функций
$ f \in L_p(D), $ обладающих тем свойством, что для любого непустого
множества $ J \subset \{1,\ldots,d\} $ обобщённая производная
$ \D^{\l \chi_J} f \in L_p(D) $ и выполняется условие
\begin{multline*}
\sup_{t^J \in (\R_+^d)^J} (t^J)^{-(\alpha^J -\l^J)}
\Omega^{(l -\l) \chi_J}(\D^{\l \chi_J} f,t^J)_{L_p(D)} = \\
\sup_{t^J \in (\R_+^d)^J}
\biggl(\prod_{j \in J} t_j^{-(\alpha_j -\l_j)}\biggr) \Omega^{(l -\l) \chi_J}(\D^{\l \chi_J} f,
t^{\s((l -\l) \chi_J)})_{L_p(D)} < \infty.
\end{multline*}

В пространстве $ (S_p^\alpha H)^{\l}(D) $ вводится норма
\begin{multline*}
\| f \|_{(S_p^\alpha H)^{\l}(D)} = \\
\max \biggl(\| f \|_{L_p(D)}, \max_{J \subset \{1, \ldots, d\}:
J \ne \emptyset} \sup_{t^J \in (\R_+^d)^J} (t^J)^{-(\alpha^J -\l^J)}
\Omega^{(l -\l) \chi_J}(\D^{\l \chi_J} f, t^J)_{L_p(D)}\biggr), \\
f \in (S_p^\alpha H)^{\l}(D).
\end{multline*}
При тех же условиях на $ \alpha, p, D $ обозначим через
$ (S_p^\alpha H)^\prime(D) $ пространство всех функций
$ f \in L_p(D), $ обладающих тем свойством, что для любого непустого
множества $ J \subset \{1,\ldots,d\} $ соблюдается условие
$$
\sup_{t^J \in (\R_+^d)^J} (t^J)^{-\alpha^J}
\Omega^{\prime l \chi_J}(f, t^J)_{L_p(D)} = \sup_{t^J \in (\R_+^d)^J}
\biggl(\prod_{j \in J} t_j^{-\alpha_j}\biggr) \Omega^{\prime l \chi_J}(f,
t^{\s(l \chi_J)})_{L_p(D)} < \infty.
$$

В пространстве $ (S_p^\alpha H)^\prime(D) $ задаётся норма
\begin{multline*}
\| f \|_{(S_p^\alpha H)^\prime(D)} = \\
\max \biggl(\| f \|_{L_p(D)}, \max_{J \subset \{1, \ldots, d\}:
J \ne \emptyset} \sup_{t^J \in (\R_+^d)^J} (t^J)^{-\alpha^J}
\Omega^{\prime l \chi_J}(f, t^J)_{L_p(D)}\biggr), \\
f \in (S_p^\alpha H)^\prime(D).
\end{multline*}

Для области $ D \subset \R^d $ при $ \alpha \in \R_+^d,\ 1 \le p \le \infty,
\theta \in \R: 1 \le \theta < \infty, \l \in \Z_+^d: \l < \alpha, $
полагая, как и выше, $ l = l(\alpha) \in \N^d, $
через $ (S_{p,\theta}^\alpha B)^{\l}(D) $ обозначим пространство всех функций
$ f \in L_p(D), $ у которых для
любого непустого множества $ J \subset \{1,\ldots,d\} $ обобщённая производная
$ \D^{\l \chi_J} f \in L_p(D) $ и соблюдается условие
\begin{multline*}
\int_{(\R_+^d)^J} (t^J)^{-\e^J -\theta (\alpha^J -\l^J)}
(\Omega^{(l -\l) \chi_J}(\D^{\l \chi_J} f, t^J)_{L_p(D)})^\theta dt^J = \\
\int_{(\R_+^d)^J} \biggl(\prod_{j \in J} t_j^{-1 -\theta (\alpha_j -\l_j)}\biggr)
(\Omega^{(l -\l) \chi_J}(\D^{\l \chi_J} f, t^{\s((l -\l) \chi_J)})_{L_p(D)})^\theta
\prod_{j \in J} dt_j < \infty.
\end{multline*}

В пространстве $ (S_{p,\theta}^\alpha B)^{\l}(D) $ фиксируется норма
\begin{multline*}
\| f \|_{(S_{p,\theta}^\alpha B)^{\l}(D)} = \\
\max \biggl(\| f\|_{L_p(D)}, \max_{J \subset \{1, \ldots, d\}: J \ne \emptyset}
\biggl(\int_{(\R_+^d)^J} (t^J)^{-\e^J -\theta (\alpha^J -\l^J)}
(\Omega^{(l -\l) \chi_J}(\D^{\l \chi_J} f, t^J)_{L_p(D)})^\theta dt^J\biggr)^{1/\theta}\biggr), \\
f \in (S_{p,\theta}^\alpha B)^{\l}(D).
\end{multline*}

При $ \theta = \infty $ положим $ (S_{p,\infty}^\alpha B)^{\l}(D) =
(S_p^\alpha H)^{\l}(D). $

При тех же условиях на параметры обозначим через
$ (S_{p,\theta}^\alpha B)^\prime(D) $ пространство всех функций
$ f \in L_p(D), $ которые для любого непустого множества $ J \subset \{1,\ldots,d\} $
подчинены условию
\begin{multline*}
\int_{(\R_+^d)^J} (t^J)^{-\e^J -\theta \alpha^J}
(\Omega^{\prime l \chi_J}(f, t^J)_{L_p(D)})^\theta dt^J = \\
\int_{(\R_+^d)^J} \biggl(\prod_{j \in J} t_j^{-1 -\theta \alpha_j}\biggr)
(\Omega^{\prime l \chi_J}(f, t^{\s(l \chi_J)})_{L_p(D)})^\theta
\prod_{j \in J} dt_j < \infty.
\end{multline*}

В пространстве $ (S_{p,\theta}^\alpha B)^\prime(D) $ определяется норма
\begin{multline*}
\| f \|_{(S_{p,\theta}^\alpha B)^\prime(D)} = \\
\max \biggl(\| f\|_{L_p(D)},
\max_{J \subset \{1, \ldots, d\}: J \ne \emptyset}
\biggl(\int_{(\R_+^d)^J} (t^J)^{-\e^J -\theta \alpha^J}
(\Omega^{\prime l \chi_J}(f, t^J)_{L_p(D)})^\theta dt^J\biggr)^{1/\theta}\biggr), \\
 \ f \in (S_{p,\theta}^\alpha B)^\prime(D).
\end{multline*}
При $ \theta = \infty $ положим $ (S_{p,\infty}^\alpha B)^\prime(D) =
(S_p^\alpha H)^\prime(D). $

В случае, когда вектор $ \l = \l(\alpha) \in \Z_+^d $ имеет компоненты
$ (\l(\alpha))_j = \max\{m \in \Z_+: m < \alpha_j\}, j =1,\ldots,d,$
пространство $ (S_{p,\theta}^\alpha B)^{\l}(D) ((S_p^\alpha H)^{\l}(D))$ обычно
обозначается $ S_{p,\theta}^\alpha B(D) (S_p^\alpha H(D)).$

Будет полезна следующая лемма (см., например, [1]).

Лемма 1.1.3

Пусть $ d \in \N, l \in \Z_+^d, 1 \le p < \infty, D $ -- открытое множество в $ \R^d. $
Тогда для функции $ f \in L_p(D), $ у которой обобщённая производная
$ \D^l f \in L_p(D), $ при $ h \in \R^d $ верно неравенство
\begin{equation*} \tag{1.1.8}
\| \Delta_h^l f \|_{L_p(D_h^l)} \le \biggl(\prod_{j \in \s(l)} |h_j|^{l_j}\biggr)
\| \D^l f \|_{L_p(D)}.
\end{equation*}

В тех же обозначениях, что и выше, принимая во внимание то обстоятельство,
что для $ f \in (S_{p,\theta}^\alpha B)^{\l}(D) $ для $ J \subset
\{1,\ldots,d\}: J \ne \emptyset, $ и $ t^J \in (\R_+^d)^J $
справедлива оценка (см., например, [17])
\begin{multline*}
(t^J)^{-(\alpha^J -\l^J)} \Omega^{(l -\l) \chi_J}(\D^{\l \chi_J} f, t^J)_{L_p(D)} = \\
\biggl(\int_{\{\tau^J \in (\R_+^d)^J: t_j < \tau_j < 2t_j, j \in
J\}} (t^J)^{-\e^J -\theta (\alpha^J -\l^J)} (\Omega^{(l -\l) \chi_J}(\D^{\l \chi_J} f,
t^J)_{L_p(D)})^\theta d\tau^J\biggr)^{1/\theta}  \le \\
\biggl(\int_{\substack{\{\tau^J \in (\R_+^d)^J:\\  t_j < \tau_j < 2t_j, j \in
J\}}} (\prod_{j \in J} 2^{1+\theta (\alpha_j -\l_j)}) (\tau^J)^{-\e^J
-\theta (\alpha^J -\l^J)} (\Omega^{(l -\l) \chi_J}(\D^{\l \chi_J} f, \tau^J)_{L_p(D)})^\theta
d\tau^J\biggr)^{1/\theta} \le \\
(\prod_{j \in J} 2^{\alpha_j -\l_j +1/\theta}) \biggl(\int_{(\R_+^d)^J}
(\tau^J)^{-\e^J -\theta (\alpha^J -\l^J)} (\Omega^{(l -\l) \chi_J}(\D^{\l \chi_J} f,
\tau^J)_{L_p(D)})^\theta d\tau^J\biggr)^{1/\theta},
\end{multline*}
заключаем, что
\begin{multline*}
(S_{p,\theta}^\alpha B)^{\l}(D) \subset (S_p^\alpha H)^{\l}(D), \\
D \text{ -- область в } \R^d, \alpha \in \R_+^d, 1 \le p \le \infty,
1 \le \theta < \infty, \l \in \Z_+^d: \l < \alpha.
\end{multline*}

Учитывая, что для $ f \in (S_{p,\theta}^\alpha B)^{\prime}(D),
t^J \in (\R_+^d)^J (J \subset \{1, \ldots, d\}: J \ne \emptyset) $
выполняется неравенство
\begin{multline*}
(t^J)^{-\alpha^J} \Omega^{\prime l \chi_J}(f, t^J)_{L_p(D)} = \\
\biggl(\int_{t^J +t^J (I^d)^J} (t^J)^{-\e^J -\theta \alpha^J} (\Omega^{\prime l \chi_J}(f,
t^J)_{L_p(D)})^\theta d\tau^J\biggr)^{1/\theta} = \\
\biggl(\int_{t^J +t^J (I^d)^J} (t^J)^{-\e^J -\theta \alpha^J} ((2 t^J)^{-\e^J}
\int_{ (t B^d)^J} \| \Delta_\xi^{l \chi_J} f\|_{L_p(D_\xi^{l \chi_J})}^p
d \xi^J)^{\theta /p} d\tau^J\biggr)^{1/\theta} \le \\
\biggl(\int_{t^J +t^J (I^d)^J} (\prod_{j \in J} 2^{1 +\theta \alpha_j})
(\tau^J)^{-\e^J -\theta \alpha^J} ((\tau^J)^{-\e^J}
\int_{ (\tau B^d)^J} \| \Delta_\xi^{l \chi_J} f\|_{L_p(D_\xi^{l \chi_J})}^p
d \xi^J)^{\theta /p} d\tau^J\biggr)^{1/\theta} = \\
\biggl(\int_{t^J +t^J (I^d)^J} (\prod_{j \in J} 2^{1 +\theta \alpha_j})
(\tau^J)^{-\e^J -\theta \alpha^J} (\prod_{j \in J} 2^{\theta /p})
(\Omega^{\prime l \chi_J}(f, \tau^J)_{L_p(D)})^\theta
d\tau^J\biggr)^{1/\theta} \le \\
(\prod_{j \in J} 2^{\alpha_j +1/\theta +1 /p}) \biggl(\int_{(\R_+^d)^J}
(\tau^J)^{-\e^J -\theta \alpha^J} (\Omega^{\prime l \chi_J}(f,
\tau^J)_{L_p(D)})^\theta d\tau^J\biggr)^{1/\theta},
\end{multline*}
заключаем, что
\begin{equation*} \tag{1.1.9}
(S_{p, \theta}^\alpha B)^\prime(D) \subset
(S_p^\alpha H)^\prime(D),
\end{equation*}
и
\begin{equation*} \tag{1.1.10}
\| f\|_{(S_p^\alpha H)^\prime(D)} \le c_4(\alpha)
\| f\|_{(S_{p, \theta}^\alpha B)^\prime(D)},
\end{equation*}
где $ c_4(\alpha) = \prod_{j=1}^d 2^{2+\alpha_j}. $

Из (1.1.7) и (1.1.8) следует, что для $ f \in (S_{p, \theta}^\alpha B)^{\l}(D), $
при $ \l \in \Z_+^d: \l < \alpha, l = l(\alpha), $ для
$ J \subset \{1,\ldots,d\}: J \ne \emptyset, $ при $ t^J \in (\R_+^d)^J, $
справедливо неравенство

\begin{multline*}
\Omega^{\prime l \chi_J}(f, t^J)_{L_p(D)} \le
\Omega^{l \chi_J}(f, t^J)_{L_p(D)} = \\
\supvrai_{\xi \in \R^d: \xi^J \in
(t B^d)^J} \| \Delta_\xi^{\l \chi_J +(l -\l) \chi_J} f
\|_{L_p(D_\xi^{\l \chi_J +(l -\l) \chi_J})} \le \\
\supvrai_{\xi \in \R^d: \xi^J \in
(t B^d)^J} \| \Delta_\xi^{\l \chi_J}(\Delta_\xi^{(l -\l) \chi_J} f)
\|_{L_p((D_\xi^{(l -\l) \chi_J})_\xi^{\l \chi_J})} \le \\
\supvrai_{\xi \in \R^d: \xi^J \in (t B^d)^J} \biggl(\prod_{j \in \s(\l \chi_J)} |\xi_j|^{\l_j}\biggr)
\| \D^{\l \chi_J} \Delta_\xi^{(l -\l) \chi_J} f
\|_{L_p(D_\xi^{(l -\l) \chi_J})} = \\
\supvrai_{\xi \in \R^d: \xi^J \in (t B^d)^J} \biggl(\prod_{j \in \s(\l \chi_J)} |\xi_j|^{\l_j}\biggr)
\| \Delta_\xi^{(l -\l) \chi_J} \D^{\l \chi_J} f
\|_{L_p(D_\xi^{(l -\l) \chi_J})} \le \\
\biggl(\prod_{j \in J} t_j^{\l_j}\biggr) \supvrai_{\xi \in \R^d: \xi^J \in (t B^d)^J}
\| \Delta_\xi^{(l -\l) \chi_J} \D^{\l \chi_J} f
\|_{L_p(D_\xi^{(l -\l) \chi_J})} = \\
\biggl(\prod_{j \in J} t_j^{\l_j}\biggr) \Omega^{(l -\l) \chi_J}(\D^{\l \chi_J} f, t^J)_{L_p(D)} =
(t^J)^{\l^J} \Omega^{(l -\l) \chi_J}(\D^{\l \chi_J} f, t^J)_{L_p(D)},
\end{multline*}
и, значит,
\begin{equation*} \tag{1.1.11}
(S_{p, \theta}^\alpha B)^{\l}(D) \subset (S_{p, \theta}^\alpha B)^\prime(D)
\end{equation*}
и
\begin{multline*} \tag{1.1.12}
\| f\|_{(S_{p, \theta}^\alpha B)^\prime(D)} \le
\| f\|_{(S_{p, \theta}^\alpha B)^{\l}(D)}, \\
\ f \in (S_{p, \theta}^\alpha B)^{\l}(D), \
\alpha \in \R_+^d, 1 \le p < \infty, 1 \le \theta \le \infty, \\
D \text{ -- произвольная область в } \R^d. \l \in \Z_+^d: \l < \alpha.
\end{multline*}

Обозначим через $ C^\infty(D) $ пространство бесконечно
дифференцируемых функций на открытом множестве $ D \subset \R^d, $ а через
$ C_0^\infty(D) $ -- пространство функций $ f \in C^\infty(\R^d), $
каждая из которых имеет компактный носитель $ \supp f \subset D. $
И ещё через $ L_1^{\loc}(D) $ обозначим пространство вещественных локально
суммируемых функций на открытом множестве $ D \subset \R^d, $ т.е.
пространство вещественных функций на $ D, $ суммируемых на любом компакте,
лежащем в $ D. $

В заключение этого пункта введём ещё несколько обозначений.

Для банахова пространства $ X $ (над $ \R$) обозначим
$ B(X) = \{x \in X: \|x\|_X \le 1\}. $

Для банаховых пространств $ X,Y $ через $ \mathcal B(X,Y) $
обозначим банахово пространство, состоящее из непрерывных линейных
операторов $ T: X \mapsto Y, $ с нормой
$$
\|T\|_{\mathcal B(X,Y)} = \sup_{x \in B(X)} \|Tx\|_Y.
$$
Отметим, что если $ X=Y,$ то $ \mathcal B(X,Y) $ является банаховой алгеброй.
\bigskip

1.2. В этом пункте приведём некоторые вспомогательные утверждения,
которые используются в следующем параграфе.

Следующее утверждение установлено в [18].

Лемма 1.2.1

Пусть $ d \in \N, 1 \le p < \infty. $ Тогда

1) при  $ j=1,\ldots,d$ для любого непрерывного линейного оператора
$ T: L_p(\R) \mapsto L_p(\R) $ существует единственный непрерывный
линейный оператор $ \mathcal T^j: L_p(\R^d) \mapsto L_p(\R^d), $ для которого
для любой функции $ f \in L_p(\R^d) $ почти для всех $ (x_1,\ldots,x_{j-1},
x_{j+1},\ldots,x_d) \in \R^{d-1} $ в $ L_p(\R) $ имеет место равенство
\begin{equation*} \tag{1.2.1}
(\mathcal T^j f)(x_1,\ldots, x_{j-1},\cdot,x_{j+1},\ldots,x_d) =
(T(f(x_1,\ldots,x_{j-1},\cdot,x_{j+1},
\ldots,x_d)))(\cdot),
\end{equation*}

2) при этом, для каждого $ j=1,\ldots,d $ отображение $ V_j^{L_p}, $ которое
каждому оператору $ T \in \mathcal B(L_p(\R), L_p(\R)) $ ставит в соответствие
оператор $ V_j^{L_p}(T) = \mathcal T^j \in \mathcal B(L_p(\R^d), L_p(\R^d)), $
удовлетворяющий (1.2.1), является непрерывным гомоморфизмом банаховой алгебры
$ \mathcal B(L_p(\R), L_p(\R)) $ в банахову алгебру $ \mathcal B(L_p(\R^d),
L_p(\R^d)), $

3) причём для любых операторов $ S,T \in \mathcal B(L_p(\R), L_p(\R)) $
при любых $ i,j =1,\ldots,d: i \ne j, $ выполняется равенство
\begin{equation*} \tag{1.2.2}
(V_i^{L_p}(S) V_j^{L_p}(T)) f = (V_j^{L_p}(T) V_i^{L_p}(S)) f, f \in L_p(\R^d).
\end{equation*}

Отметим, что в дальнейшем символы $ L_p $ в обозначениях операторов $ V_j^{L_p} $
будем опускать.

Для $ \phi \in L_\infty(\R) $ через $ M_\phi $ обозначим линейный оператор
в пространстве локально суммируемых функций на $ \R, $ определяемый
равенством $ M_\phi g = \phi g, $ где $ g \in L_1^{\loc}(\R). $

Через $ \chi_A $ будем обозначать характеристическую функцию множества
$ A \subset \R^d. $

Так же, как соответствующие усверждения в [17], устанавливаются
леммы 1.2.2, 1.2.3.

Лемма 1.2.2

Пусть $ d \in \N, l \in \N^d, 1 \le p < \infty. $ Тогда существует константа
$ c_1(d,l) >0 $ такая, что для любых $ x^0, X^0 \in \R^d, \delta, \Delta \in
\R_+^d $ таких, что $ Q = (x^0 +\delta I^d) \subset (X^0 +\Delta I^d), $ для
любой функции $ f \in L_p(\R^d), $ для любого $ \xi \in \R^d, $ для любых
множеств $ J, J^\prime \subset \{1,\ldots,d\} $ имеет место неравенство
\begin{multline*} \tag{1.2.3}
\biggl\| \Delta_\xi^{l \chi_J} ((\prod_{j \in J^\prime}
V_j(E -M_{\chi_{X_j^0 +\Delta_j I}}
P_{\delta_j, x_j^0}^{1,l_j -1})) f)\biggr\|_{L_p(Q_\xi^{l \chi_J})} \le \\
c_1 \biggl(\prod_{j \in J^\prime \setminus J} \delta_j^{-1/p}\biggr)
\biggl(\int_{(\delta B^d)^{J^\prime \setminus J}} \int_{ Q_\xi^{l
\chi_{J \cup J^\prime}}} |(\Delta_\xi^{l \chi_{J \cup J^\prime}}
f)(x)|^p dx d\xi^{J^\prime \setminus J}\biggr)^{1/p},
\end{multline*}
где $ E $ -- тождественный оператор.

Доказательство.

Сначала убедимся в справедливости (1.2.3) в случае, когда
множества $ J $ и $ J^\prime $ не пересекаются. А именно, покажем,
что если выполнены условия леммы 1.2.2 и $ J \cap J^\prime =
\emptyset, $ то
\begin{multline*} \tag{1.2.4}
\|\Delta_\xi^{l \chi_J} ((\prod_{j \in J^\prime} V_j(E
-M_{\chi_{X_j^0 +\Delta_j I}}
P_{\delta_j, x_j^0}^{1, l_j -1})) f)\|_{L_p(Q_\xi^{l \chi_J})}\\
\le c_2^k \biggl(\prod_{j \in J^\prime} \delta_j^{-1/p}\biggr)
\biggl(\int_{(\delta B^d)^{J^\prime}} \int_{ Q_\xi^{l \chi_{J \cup
J^\prime}}} |(\Delta_\xi^{l \chi_{J \cup J^\prime}} f)(x)|^p dx
d\xi^{J^\prime}\biggr)^{1/p},
\end{multline*}
где $ k = \card J^\prime, $ а $ c_2(d,l) = \max_{j=1,\ldots,d}
c_3(1, l_j-1) (см. (1.1.6)). $

Доказательство (1.2.4) проведём по индукции относительно $ k. $ При
$ k=1, $ а, следовательно, $ J^\prime = \{j\} (j \in \{1,\ldots,d\}
\setminus J), $ используя теорему Фубини, (1.2.1), (1.1.6), в условиях леммы имеем
\begin{multline*}
\|\Delta_\xi^{l \chi_J} (V_j(E -M_{\chi_{X_j^0 +\Delta_j I}}
P_{\delta_j, x_j^0}^{1, l_j -1}) f)\|_{L_p (Q_\xi^{l \chi_J})}^p  =\\
\int_{\substack{ \prod_{i=1,\ldots,d: i \ne j} \\ (x_i^0
+\delta_i I)_{l_i (\chi_J)_i \xi_i}}} \int_{(x_j^0 +\delta_j I)}
|((E -M_{\chi_{X_j^0 +\Delta_j I}}
P_{\delta_j, x_j^0}^{1, l_j -1}) (\Delta_\xi^{l \chi_J} f)
(x_1, \ldots, x_{j-1}, \cdot, x_{j+1}, \ldots,x_d))(x_j)|^p \\
\times dx_j dx_1 \ldots dx_{j-1} dx_{j+1} \ldots dx_d =\\
\int_{\substack{ \prod_{i=1,\ldots,d: i \ne j} \\ (x_i^0
+\delta_i I)_{l_i (\chi_J)_i \xi_i}}} \int_{(x_j^0 +\delta_j I)}
|((E -P_{\delta_j, x_j^0}^{1, l_j -1}) (\Delta_\xi^{l \chi_J} f)
(x_1, \ldots, x_{j-1}, \cdot, x_{j+1}, \ldots,x_d))(x_j)|^p \\
\times dx_j dx_1 \ldots dx_{j-1} dx_{j+1} \ldots dx_d \le \\
c_2^p \delta_j^{-1} \int_{\delta_j B^1} \int_{ Q_\xi^{l \chi_{J \cup
J^\prime}}} |(\Delta_\xi^{l \chi_{J \cup J^\prime}} f)(x)|^p dx
d\xi_j.
\end{multline*}
Отсюда следует (1.2.4) в случае, когда $ J^\prime = \{j\}, $ а $ j \in
\{1,\ldots,d\} \setminus J. $

Предположим, что при $ k \in \N $ неравенство (1.2.4) имеет место
при соблюдении условий леммы 1.2.2 для $ J,J^\prime \subset \{1.\ldots,d\}:
J \cap J^\prime = \emptyset $ и $ \card J^\prime \le k. $ Покажем, что тогда оно
справедливо в условиях леммы 1.2.2 в ситуации, когда $ J,J^\prime \subset
\{1,\ldots,d\}: J \cap J^\prime = \emptyset, $ и $ \card J^\prime = k+1. $

Для этого, фиксируя для таких множеств $ J,J^\prime $ элемент $ j \in J^\prime $
и применяя (1.2.4) сначала с множествами $ J $ и $ \{j\} $ (см. (1.2.2)), а
затем с множествами $ J \cup \{j\} $ и $ J^\prime \setminus \{j\}, $ получаем
(1.2.4) в рассматриваемой ситуации, что и завершает доказательство (1.2.4) в общей ситуации.

Для вывода (1.2.3) заметим, что для $ d,l,p,f,\delta,x^0, \Delta, X^0,$
удовлетворяющих условиям леммы 1.2.2, при $ j=1,\ldots,d $ и $ \xi_j \in \R $
в виду п.2 леммы 1.2.1 в $ (X^0 +\Delta I^d)_{l_j \xi_j e_j} $
справедливо равенство
\begin{equation*} \tag{1.2.5}
\Delta_{\xi_j e_j}^{l_j} (V_j(E -M_{\chi_{X_j^0 +\Delta_j I}}
P_{\delta_j, x_j^0}^{1,l_j-1}) f) = \Delta_{\xi_j e_j}^{l_j} f -
\Delta_{\xi_j e_j}^{l_j} (V_j(M_{\chi_{X_j^0 +\Delta_j I}}
P_{\delta_j, x_j^0}^{1,l_j-1})f).
\end{equation*}
При этом, в силу (1.2.1) почти для всех
$ (x_1,\ldots,x_{j-1},x_{j+1},\ldots,x_d) \in \R^{d-1} $ для почти
всех $ x_j \in (X_j^0 +\Delta_j I)_{l_j \xi_j} $ выполняется равенство
\begin{multline*}
(\Delta_{\xi_j e_j}^{l_j} (V_j (M_{\chi_{X_j^0 +\Delta_j I}}
P_{\delta_j, x_j^0}^{1,l_j-1}) f))(x_1,\ldots,x_{j-1},x_j,x_{j+1},\ldots,x_d)
= \\
(\Delta_{\xi_j}^{l_j} ((M_{\chi_{X_j^0 +\Delta_j I}}
P_{\delta_j,x_j^0}^{1,l_j-1}) f(x_1,\ldots,x_{j-1},
\cdot,x_{j+1},\ldots,x_d)))(x_j) = \\
(\Delta_{\xi_j}^{l_j} (P_{\delta_j,x_j^0}^{1,l_j-1}
f(x_1,\ldots,x_{j-1},\cdot,x_{j+1},\ldots,x_d)))(x_j) =0.
\end{multline*}

Сопоставляя сказанное с (1.2.5), приходим к выводу, что почти для
всех $ x \in (X^0 +\Delta I^d)_{l_j \xi_j e_j} $ имеет место равенство
$ (\Delta_{\xi_j e_j}^{l_j} (V_j(E -M_{\chi_{X_j^0 +\Delta_j I}}
P_{\delta_j, x_j^0}^{1,l_j-1}) f))(x) = (\Delta_{\xi_j e_j}^{l_j} f)(x). $

Для доказательства (1.2.3) покажем, что при $ d \in \N, l \in \N^d,
1 \le p < \infty, x^0, X^0 \in \R^d, \delta, \Delta \in \R_+^d $ для
$ f \in L_p(\R^d), \xi \in \R^d $ и $ \mathcal J \subset \{1,\ldots,d\} $ почти
для всех $ x \in (X^0 +\Delta I^d)_\xi^{l \chi_{\mathcal J}} $ выполняется
равенство
\begin{equation*} \tag{1.2.6}
\biggl((\prod_{j \in \mathcal J} \Delta_{\xi_j e_j}^{l_j}) ((\prod_{j \in
\mathcal J} V_j(E -M_{\chi_{X_j^0 +\Delta_j I}}
P_{\delta_j, x_j^0}^{1,l_j-1})) f)\biggr)(x) =
\biggl((\prod_{j \in \mathcal J} \Delta_{\xi_j e_j}^{l_j}) f\biggr)(x).
\end{equation*}

Равенство (1.2.6) установим по индукции относительно $ \card \mathcal J. $

Как показано выше, в случае, когда $ \card \mathcal J =1, $
равенство (1.2.6) справедливо. Предположим, что при $ k \in \N $
оно верно в случае, когда $ \card \mathcal J \le k. $ Проверим,
что тогда оно соблюдается и в ситуации, когда $ \card \mathcal J = k+1. $

Фиксируя $ j \in \mathcal J $ и применяя (1.2.6)  сначала с
множеством $ \{j\}, $ а затем с множеством $ \mathcal J \setminus
\{j\}, $ с учётом (1.2.2) получаем (1.2.6) в случае, когда $ \card
\mathcal J = k+1, $ что и завершает вывод (1.2.6) в общем случае.

Теперь при $ d \in \N, l \in \N^d, 1 \le p < \infty, x^0, X^0 \in \R^d,
\delta, \Delta \in \R_+^d $ для $ f \in L_p(\R^d), \xi \in \R^d $ и любых
множеств $ J,J^\prime \subset \{1,\ldots,d\}, $ используя (1.2.2) и
(1.2.6), имеем
\begin{multline*} \tag{1.2.7}
\biggl(\Delta_\xi^{l \chi_J} ((\prod_{j \in J^\prime} V_j(E
-M_{\chi_{X_j^0 +\Delta_j I}} P_{\delta_j, x_j^0}^{1,l_j -1})) f)\biggr)(x) =\\
\biggl(\Delta_\xi^{l \chi_J} ((\prod_{j \in J^\prime \setminus J}
V_j(E -M_{\chi_{X_j^0 +\Delta_j I}} P_{\delta_j, x_j^0}^{1,l_j-1})) f)\biggr)(x), x \in
(X^0 +\Delta I^d)_\xi^{l \chi_J}.
\end{multline*}

Соединяя (1.2.7) с (1.2.4), получаем (1.2.3). $ \square $

Лемма 1.2.3

Пусть $ d \in \N, l \in \Z_+^d, 1 \le p < \infty $ и $ \rho, \sigma \in \R_+^d. $
Тогда существует константа $ c_3(d,l,\rho,\sigma) >0 $ такая, что для любых
$ x^0, \bm x^0, X^0 \in \R^d $ и $ \delta, \bm \delta, \Delta \in \R_+^d: $
$ (x^0 +\sigma \delta I^d) \subset (\bm x^0 +\bm \delta I^d) \subset
(\bm x^0 +\rho \delta I^d) \cap (X^0 +\Delta I^d), $
для любого множества $ J \subset \{1,\ldots,d\}, $ для $ f \in L_p(\R^d) $ соблюдается неравенство
\begin{equation*} \tag{1.2.8}
\biggl\| (\prod_{j \in J} V_j(E -M_{\chi_{X_j^0 +\Delta_j I}}
P_{\sigma_j \delta_j,x_j^0}^{1,l_j})) f \biggr\|_{L_p(D)} \le
c_3 \|f\|_{L_p(D)},
\end{equation*}
где $ D = (\bm x^0 +\bm \delta I^d). $

Доказательство.

Понятно, что (1.2.8) достаточно доказать в случае, когда $ J =
\{j\}, j \in \{1,\ldots,d\}. $ Переходя к выводу (1.2.8) в этом
случае, заметим, что в условиях леммы 1.2.3 имеет место включение
\begin{equation*} \tag{1.2.9}
(\bm x^0_j +\bm \delta_j I) \subset (\bm x^0_j +\rho_j \delta_j I) \subset
(x_j^0 +\rho_j \delta_j B^1).
\end{equation*}

Далее, в виду п. 2 леммы 1.2.1 имеем
\begin{equation*} \tag{1.2.10}
\| V_j(E -M_{\chi_{X_j^0 +\Delta_j I}} P_{\sigma_j \delta_j, x_j^0}^{1,l_j}) f\|_{L_p(D)} \le
\|f\|_{L_p(D)} +\| V_j(M_{\chi_{X_j^0 +\Delta_j I}}
P_{\sigma_j \delta_j,x_j^0}^{1,l_j}) f\|_{L_p(D)}.
\end{equation*}

Оценивая второе слагаемое в правой части (1.2.10), с помощью теоремы Фубини и
соотношений (1.2.1), (1.1.3) (с учётом (1.2.9)), (1.1.5), находим, что
\begin{multline*} \tag{1.2.11}
\| V_j(M_{\chi_{X_j^0 +\Delta_j I}} P_{\sigma_j \delta_j, x_j^0}^{1,l_j}) f\|_{L_p(D)}^p =\\
\int_{ \prod_{i=1,\ldots,d: i \ne j} (\bm x_i^0 +\bm \delta_i I)}
\int_{(\bm x^0_j +\bm \delta_j I)} | \chi_{X_j^0 +\Delta_j I}(x_j)
(P_{\sigma_j \delta_j,x_j^0}^{1,l_j}
f(x_1,\ldots,x_{j-1},\cdot,x_{j+1},\ldots,x_d))(x_j)|^p \\
\times dx_j dx_1 \ldots dx_{j-1} dx_{j+1} \ldots dx_d =\\
\int_{ \prod_{i=1,\ldots,d: i \ne j} (\bm x_i^0 +\bm \delta_i I)}
\int_{(\bm x^0_j +\bm \delta_j I)} | (P_{\sigma_j \delta_j,x_j^0}^{1,l_j}
f(x_1,\ldots,x_{j-1},\cdot,x_{j+1},\ldots,x_d))(x_j)|^p \\
\times dx_j dx_1 \ldots dx_{j-1} dx_{j+1} \ldots dx_d \le\\
\int_{ \prod_{i=1,\ldots,d: i \ne j} (\bm x_i^0 +\bm \delta_i I)}
c_1(1,l_j,0,\rho_j,\sigma_j)^p \\
\times \int_{(x_j^0 +\sigma_j \delta_j I)}
|(P_{\sigma_j \delta_j,x_j^0}^{1,l_j}
f(x_1,\ldots,x_{j-1},\cdot,x_{j+1},\ldots,x_d))(x_j)|^p \\
\times dx_j dx_1 \ldots dx_{j-1} dx_{j+1} \ldots dx_d \le \\
c_4^p \int_{ \prod_{i=1,\ldots,d: i \ne j} (\bm x_i^0 +\bm \delta_i I)}
\int_{(x_j^0 +\sigma_j \delta_j I)}
| f(x_1,\ldots,x_{j-1},x_j,x_{j+1},\ldots,x_d)|^p \\
\times dx_j dx_1 \ldots dx_{j-1} dx_{j+1} \ldots dx_d \le 
c_4^p \|f\|_{L_p(D)}^p.
\end{multline*}
Соединяя (1.2.10) с (1.2.11), приходим к (1.2.8) в
ситуации, когда $ J = \{j\}, j \in \{1,\ldots,d\}. \square $

В [1] доказана следующая лемма.

Лемма 1.2.4

Пусть $ d \in \N, l \in \Z_+^d, \delta, \Delta \in \R_+^d, x^0, X^0 \in \R^d. $
Тогда при $ 1 \le p < \infty $ для любой функции $ f \in L_p(\R^d) $ почти для
всех $ x \in \R^d $ справедливо равенство
\begin{equation*} \tag{1.2.12}
\chi_{X^0 +\Delta I^d}(x) (P_{\delta, x^0}^{d,l}(f \mid_{x^0 +\delta I^d}))(x)
= \biggl((\prod_{j=1}^d V_j(M_{\chi_{X_j^0 +\Delta_j I}}
P_{\delta_j, x_j^0}^{1,l_j})) f\biggr)(x).
\end{equation*}
\bigskip

1.3. В этом пункте вводятся в рассмотрение пространства кусочно-
полиномиальных функций и операторы в них, которые используются для построения
средств приближения функций из изучаемых нами пространств.
Но сначала приведём некоторые вспомогательные сведения.

Для $ d \in \N, y \in \R^d $ положим
$$
\mn(y) = \min_{j=1,\ldots,d} y_j
$$
и для банахова пространства $ X, $ вектора $ x \in X $ и семейства
$ \{x_\kappa \in X, \kappa \in \Z_+^d\} $ будем писать $ x =
\lim_{ \mn(\kappa) \to \infty} x_\kappa, $ если для любого $ \epsilon >0 $
существует $ n_0 \in \N $ такое, что для любого $ \kappa \in \Z_+^d, $
для которого $ \mn(\kappa) > n_0, $
справедливо неравенство $ \|x -x_\kappa\|_X < \epsilon. $

Пусть $ X $ -- банахово пространство (над $ \R $), $ d \in \N $ и
$ \{ x_\kappa \in X: \kappa \in \Z_+^d\} $ -- семейство векторов.
Тогда под суммой ряда $ \sum_{\kappa \in \Z_+^d} x_\kappa $ будем
понимать вектор $ x \in X, $ для которого выполняется равенство
$ x = \lim_{\mn(k) \to \infty} \sum_{\kappa \in \Z_+^d(k)} x_\kappa $ (см. (1.1.1)).

При $ d \in \N $ через $ \Upsilon^d $ обозначим множество
$$
\Upsilon^d = \{ \epsilon \in \Z^d: \epsilon_j \in \{0,1\},
j=1,\ldots,d\}.
$$

Как показано в [14], имеет место следующая лемма.

   Лемма 1.3.1

Пусть $ X $ -- банахово пространство, а вектор $ x \in X $ и
семейство $ \{x_\kappa \in X: \kappa \in \Z_+^d\} $ таковы, что
$ x = \lim_{ \mn(\kappa) \to \infty} x_\kappa, $ Тогда для семейства
$ \{ \mathcal X_\kappa \in X, \kappa \in \Z_+^d \}, $ элементы которого определяются
равенством
$$
\mathcal X_\kappa = \sum_{\epsilon \in \Upsilon^d: \s(\epsilon)
\subset \s(\kappa)} (-\e)^\epsilon x_{\kappa -\epsilon}, \kappa
\in \Z_+^d,
$$
справедливо равенство
$$
x = \sum_{\kappa \in \Z_+^d} \mathcal X_\kappa.
$$

Замечание

Как отмечалось в [14], если для семейства векторов $ \{x_\kappa \in X,
\kappa \in \Z_+^d\} $ банахова пространства $ X $ ряд
$ \sum_{\kappa \in \Z_+^d} \| x_\kappa \|_X $ сходится, то для любой
последовательности подмножеств $ \{Z_n \subset \Z_+^d, n \in \Z_+\}, $
таких, что $ \card Z_n < \infty, Z_n \subset Z_{n+1}, n \in \Z_+, $
и $ \cup_{ n \in \Z_+} Z_n = \Z_+^d, $
в $ X $ соблюдается равенство
$ \sum_{\kappa \in \Z_+^d} x_\kappa =
\lim_{ n \to \infty} \sum_{\kappa \in Z_n} x_\kappa. $

Введём в рассмотрение систему разбиений единицы на открытых множествах,
используемую для построения средств приближения функций из изучаемых
пространств. Для этого обозначим через $ \psi^{1,0} $ характеристическую функцию
интервала $ I, $ т.е. функцию, определяемую равенством
$$
\psi^{1,0}(x) = \begin{cases} 1, & \text{ для } x \in I; \\
0, & \text{ для } x \in \R \setminus I.
\end{cases}
$$
При $ m \in \N $ положим
$$
\psi^{1,m}(x) = \int_I \psi^{1, m-1}(x -y) dy \ (\text{см. } [19]),
$$
а для $ d \in \N, m \in \Z_+^d $ определим
$$
\psi^{d,m}(x) = \prod_{j=1}^d \psi^{1,m_j}(x_j), x =
(x_1,\ldots,x_d) \in \R^d.
$$

Для $ d \in \N, m,n \in \Z^d: m \le n, $ обозначим
\begin{equation*} \tag{1.3.1}
\Nu_{m,n}^d = \{ \nu \in \Z^d: m \le \nu \le n \} = \prod_{j=1}^d
\Nu_{m_j,n_j}^1.
\end{equation*}

Опираясь на определения, используя индукцию, нетрудно проверить
следующие свойства функций $ \psi^{d,m}, d \in \N, m \in \Z_+^d. $

1) При $ d \in \N, m \in \Z_+^d $
$$
\sgn \psi^{d,m}(x) = \begin{cases} 1, \text{ для } x \in ((m +\e) I^d); \\
0, \text{ для } x \in \R^d \setminus ((m +\e) I^d),
\end{cases}
$$

2) при $ d \in \N, m \in \Z_+^d $ для каждого $ \lambda \in
\Z_+^d(m) $ (обобщённая) производная $ \D^\lambda \psi^{d,m} \in
L_\infty(\R^d), $

3) при $ d \in \N, m \in \Z_+^d $ почти для всех $ x \in \R^d $
справедливо равенство
$$
\sum_{\nu \in \Z^d} \psi^{d,m}(x -\nu) =1,
$$

4) при $ m \in \N $ для всех $ x \in \R $ (при $ m =0 $ почти для всех $ x \in \R $)
имеет место равенство
\begin{equation*} \tag{1.3.2}
\psi^{1,m}(x) = \sum_{\mu \in \Nu_{0, m+1}^1} a_{\mu}^m
\psi^{1,m}(2x -\mu),
\end{equation*}
где $ a_\mu^m = 2^{-m} C_{m+1}^\mu. $
Используя разложение Ньютона для $ (1+1)^{m+1} $  и $ (-1+1)^{m+1}, $ легко
проверить, что при $ m \in \Z_+ $ выполняются равенства
\begin{equation*} \tag{1.3.3}
\sum_{\mu \in \Nu_{0,m +1}^1 \cap (2 \Z)} a_\mu^m =1,
\sum_{\mu \in \Nu_{0,m +1}^1 \cap (2 \Z +1)} a_\mu^m =1.
\end{equation*}

При $ d \in \N $ для $ t \in \R^d $ через $ 2^t $ будем обозначать
вектор $ 2^t = (2^{t_1}, \ldots, 2^{t_d}). $

Для $ d \in \N, m,\kappa \in \Z_+^d, \nu \in \Z^d $ обозначим
$$
g_{\kappa, \nu}^{d,m}(x) = \psi^{d,m}(2^\kappa x -\nu) =
\prod_{j =1}^d \psi^{1,m_j}( 2^{\kappa_j} x_j -\nu_j), x \in \R^d.
$$
Из первого среди приведенных выше свойств функций $ \psi^{d,m} $
следует, что при $ d \in \N, m,\kappa \in \Z_+^d, \nu \in \Z^d $ носитель
\begin{equation*} \tag{1.3.4}
\supp g_{\kappa,\nu}^{d,m} =
2^{-\kappa} \nu +2^{-\kappa} (m +\e) \overline I^d.
\end{equation*}
При $ d \in \N, \kappa \in \Z_+^d, \nu \in \Z^d $ обозначим
\begin{equation*} \tag{1.3.5}
Q_{\kappa, \nu}^d = 2^{-\kappa} \nu +2^{-\kappa} I^d,
\overline Q_{\kappa, \nu}^d = 2^{-\kappa} \nu +2^{-\kappa} \overline I^d.
\end{equation*}

Отметим некоторые полезные для нас свойства носителей функций
$ g_{\kappa,\nu}^{d,m}. $

При $ d \in \N, m,\kappa \in \Z_+^d $ для каждого $ \nu^\prime \in \Z^d $
имеет место равенство
\begin{equation*} \tag{1.3.6}
\{ \nu \in \Z^d: Q_{\kappa, \nu^\prime}^d \cap
\supp g_{\kappa, \nu}^{d,m} \ne \emptyset\} = \nu^\prime +\Nu_{-m,0}^d.
\end{equation*}

Из свойства 3) функций $ \psi^{d,m} $ вытекает, что при $ d \in \N,
m, \kappa \in \Z_+^d $ для любого  открытого множества $ U \subset \R^d $
почти для всех $ x \in U $ соблюдается равенство
\begin{equation*} \tag{1.3.7}
\sum_{ \nu \in \Z^d: \supp g_{\kappa, \nu}^{d,m} \cap U \ne \emptyset}
g_{\kappa, \nu}^{d,m}(x) =1.
\end{equation*}

При $ d \in \N $ для $ x,y \in \R^d $ будем обозначать
$$
(x,y) = \sum_{j =1}^d x_j y_j.
$$

Имея в виду свойство 2) функций $ \psi^{d,m}, $ отметим, что при
$ d \in \N, m,\kappa \in \Z_+^d, \nu \in \Z^d, \lambda \in \Z_+^d(m)$
выполняется равенство
\begin{multline*} \tag{1.3.8}
\| \D^\lambda g_{\kappa, \nu}^{d,m} \|_{L_\infty (\R^d)} =
2^{(\kappa, \lambda)} \| \D^\lambda \psi^{d,m} \|_{L_\infty(\R^d)} =
c_1(d,m,\lambda) 2^{(\kappa, \lambda)}.
\end{multline*}
Введём в рассмотрение следующие пространства кусочно-полиномиальных
функций.
При $ d \in \N, l \in \Z_+^d, m \in \N^d, \kappa \in \Z_+^d $ и открытого
множества $ U \subset \R^d, $ полагая
\begin{equation*} \tag{1.3.9}
N_\kappa^{d,m,U} = \{\nu \in \Z^d: \supp g_{\kappa, \nu}^{d,m}
\cap U \ne \emptyset\},
\end{equation*}
через $ \mathcal P_\kappa^{d,l,m,U} $ обозначим линейное пространство,
состоящее из функций $ f: \R^d \mapsto \R, $ для каждой из которых существует
набор полиномов
$ \{f_\nu \in \mathcal P^{d,l}, \nu \in N_\kappa^{d,m,U}\} $ такой, что
для $ x \in \R^d $ выполняется равенство
\begin{equation*} \tag{1.3.10}
f(x) = \sum_{\nu \in N_\kappa^{d,m,U}} f_\nu(x) g_{\kappa,\nu}^{d,m}(x).
\end{equation*}

Замечание

В случае, когда $ \card N_\kappa^{d,m,U} < \infty, $ не требуется пояснять, что
понимается под суммой $ \sum_{\nu \in N_\kappa^{d,m,U}} f_\nu(x) g_{\kappa,\nu}^{d,m}(x). $
При этом нетрудно проверить, что при $ d \in \N, l \in \Z_+^d, m \in \N^d,
\kappa \in \Z_+^d $ и ограниченного открытого множества $ U \subset \R^d $
отображение, которое каждому набору полиномов $ \{f_\nu \in \mathcal P^{d, l},
\nu \in N_\kappa^{d,m,U} \} $ ставит в соответствие функцию $ f, $ задаваемую
равенством (1.3.10), является изоморфизмом прямого произведения
$ \card N_\kappa^{d,m,U} $ экземпляров пространства $ \mathcal P^{d,l} $
на пространство $ \mathcal P_\kappa^{d,l,m,U}. $

Если же $ \card N_\kappa^{d,m,U} = \infty, $ то фиксируя некоторое
биективное отображение $ \N \ni s \mapsto \nu^s \in N_\kappa^{d,m,U}, $ и
учитывая, что для $ x \in \R^d $ в силу (1.3.4) $ \card \{s \in \N:
g_{\kappa,\nu^s}^{d,m}(x) \ne 0\} < \infty, $ положим
\begin{equation*}
\sum_{\nu \in N_\kappa^{d,m,U}} f_\nu(x) g_{\kappa,\nu}^{d,m}(x) =
\sum_{s =1}^\infty f_{\nu^s}(x) g_{\kappa,\nu^s}^{d,m}(x) =
\sum_{s \in \N: g_{\kappa,\nu^s}^{d,m}(x) \ne 0} f_{\nu^s}(x) g_{\kappa,\nu^s}^{d,m}(x),
\end{equation*}
причём для любого компакта $ K \subset \R^d $ для $ x \in K $ сумма
$$
\sum_{\nu \in N_\kappa^{d,m,U}} f_\nu(x) g_{\kappa,\nu}^{d,m}(x) =
\sum_{\nu \in N_\kappa^{d,m,U}: K \cap \supp g_{\kappa,\nu}^{d,m} \ne \emptyset}
f_\nu(x) g_{\kappa,\nu}^{d,m}(x),
$$
а вторая сумма ввиду (1.3.4) содержит конечное число слагаемых (суммируемых на
$ K $ функций), тем самым, имеет место включение $ \mathcal P_\kappa^{d,l,m,U}
\subset L_1^{\loc}(\R^d). $

Опираясь на (1.3.2), так же, как лемма 1.2.1 из [10], устанавливается следующая
лемма.

Лемма 1.3.2

Пусть $ d \in \N, l \in \Z_+^d, m \in \N^d, \kappa \in \Z_+^d, U $ --
открытое множество в $ \R^d. $ Тогда при $ j =1,\ldots,d $
линейный оператор $ H_\kappa^{j,d,l,m,U}:
\mathcal P_\kappa^{d,l,m,U} \mapsto \mathcal P_{\kappa +e_j}^{d,l,m,U}, $
значение которого на функции $ f \in \mathcal P_\kappa^{d,l,m,U}, $ задаваемой
равенством (1.3.10), определяется соотношением
\begin{multline*} \tag{1.3.11}
(H_\kappa^{j,d,l,m,U} f)(x) = \\
\sum_{\nu \in N_{\kappa +e_j}^{d,m,U}}
\biggl(\sum_{\substack{\nu^\prime \in N_\kappa^{d,m,U}, \mu_j \in \Nu_{0, m_j +1}^1: \\
 2 \nu^\prime_j +\mu_j = \nu_j, \nu^\prime_i = \nu_i, i = 1,\ldots,d, i \ne j }} a_{\mu_j}^{m_j}
f_{\nu^\prime}(x)\biggr) g_{\kappa +e_j,\nu}^{d,m}(x), x \in \R^d,
\end{multline*}
обладает тем свойством, что для $ f \in \mathcal P_\kappa^{d,l,m,U} $ выполняется равенство
\begin{equation*}
(H_\kappa^{j,d,l,m,U} f) \mid_{U} = f \mid_{U}.
\end{equation*}

Отметим, что если в формулировке леммы 1.3.2 множество $ U $ не ограниченное, то
оператор $ H_\kappa^{j,d,l,m,U}, $ вообще говоря, многозначный.

Ниже понадобятся следующие объекты.

При $ m \in \N^d, \epsilon \in \Upsilon^d, \nu \in \Z^d $ обозначим
через $ \M_{\epsilon}^m(\nu) $ множество наборов чисел
\begin{multline*} \tag{1.3.12}
\M_{\epsilon}^m(\nu) = \{ \m^{\epsilon} = \{ \m_j \in \Nu_{0, m_j +1}^1,
j \in \s(\epsilon)\}: \\
(\nu_j -\m_j) /2 \in \Z \ \forall j \in \s(\epsilon)\} = \\
\prod_{j \in \s(\epsilon)} \{ \m_j \in \Nu_{0, m_j +1}^1: (\nu_j -\m_j) /2 \in \Z\} =
\prod_{j \in \s(\epsilon)} \M_1^{m_j}(\nu_j),
\end{multline*}
и каждой паре $ \nu \in \Z^d, \m^{\epsilon} \in \M_{\epsilon}^m(\nu) $ сопоставим
$ \n_{\epsilon}(\nu,\m^{\epsilon}) \in \Z^d, $ полагая
\begin{equation*} \tag{1.3.13}
(\n_{\epsilon}(\nu,\m^{\epsilon}))_j = \begin{cases} (\nu_j -\m_j) /2,
j \in \s(\epsilon); \\
\nu_j, j \in \{1,\ldots,d\} \setminus \s(\epsilon).
\end{cases}
\end{equation*}

В дальнейшем будет полезен следующий факт. При $ d \in \N, m \in
\N^d $ для $ \nu \in \Z^d, \epsilon, \epsilon^\prime \in
\Upsilon^d: \s(\epsilon) \cap \s(\epsilon^\prime) = \emptyset, $
(а, значит, $ \epsilon +\epsilon^\prime \in \Upsilon^d,
\s(\epsilon +\epsilon^\prime) = \s(\epsilon) \cup
\s(\epsilon^\prime), $) и любых $ \m^{\epsilon +\epsilon^\prime}
\in \M_{\epsilon +\epsilon^\prime}^m(\nu), \m^{\epsilon} \in
\M_{\epsilon}^m(\nu), \m^{\epsilon^\prime} \in
\M_{\epsilon^\prime}^m(\nu): (\m^{\epsilon +\epsilon^\prime})_j =
(\m^{\epsilon})_j, j \in \s(\epsilon), (\m^{\epsilon
+\epsilon^\prime})_j = (\m^{\epsilon^\prime})_j, j \in
\s(\epsilon^\prime), $ выполняется равенство
\begin{equation*} \tag{1.3.14}
\n_{\epsilon +\epsilon^\prime}(\nu, \m^{\epsilon +\epsilon^\prime}) =
\n_{\epsilon^\prime}(\n_{\epsilon}(\nu, \m^{\epsilon}), \m^{\epsilon^\prime}).
\end{equation*}

Замечание.

При $ d \in \N, m \in \N^d $ для $ \kappa \in \Z_+^d,
\epsilon \in \Upsilon^d: \s(\epsilon) \subset \s(\kappa), $ для открытого
множества $ U \subset \R^d $ и $ \nu \in N_{\kappa}^{d,m,U}, \m^{\epsilon} \in
\M_{\epsilon}^m(\nu) $ имеет место включение
\begin{equation*} \tag{1.3.15}
\n_{\epsilon}(\nu,\m^{\epsilon}) \in N_{\kappa -\epsilon}^{d,m,U}.
\end{equation*}

В самом деле, при соблюдении условий замечания, ввиду (1.3.4) и открытости $ U $
выбирая $ x \in U \cap (2^{-\kappa} \nu +2^{-\kappa} (m +\e) I^d), $
получаем, что
\begin{equation*}
2^{-\kappa_i} \nu_i < x_i < 2^{-\kappa_i} \nu_i +2^{-\kappa_i}
(m_i +1) \text{ при } i = 1,\ldots,d,
\end{equation*}
откуда с учётом (1.3.13) имеем
\begin{multline*}
2^{-(\kappa -\epsilon)_i} (\n_{\epsilon}(\nu,\m^{\epsilon}))_i =
2^{-\kappa_i} \nu_i < x_i < 2^{-\kappa_i} \nu_i +
2^{-\kappa_i} (m_i +1) = \\
2^{-(\kappa -\epsilon)_i} (\n_{\epsilon}(\nu,\m^{\epsilon}))_i +
2^{-(\kappa -\epsilon)_i} (m_i +1) \text{ при } i \in \{1,\ldots,d\} \setminus \s(\epsilon); \\
2^{-(\kappa -\epsilon)_j} (\n_{\epsilon}(\nu,\m^{\epsilon}))_j =
2^{-\kappa_j +1} (\nu_j -\m_j) /2 =
2^{-\kappa_j} \nu_j -2^{-\kappa_j} \m_j \le \\
2^{-\kappa_j} \nu_j < x_j < 2^{-\kappa_j} \nu_j +2^{-\kappa_j} (m_j +1) = \\
2^{-\kappa_j +1} (\nu_j -\m_j) /2 +
2^{-\kappa_j +1} (\m_j /2) +2^{-\kappa_j +1} (m_j +1) /2 = \\
2^{-\kappa_j +1} (\nu_j -\m_j) /2 +2^{-\kappa_j +1} (\m_j +m_j +1) /2 \le
2^{-\kappa_j +1} (\nu_j -\m_j) /2 +2^{-\kappa_j +1} (m_j +1) = \\
2^{-(\kappa -\epsilon)_j} (\n_{\epsilon}(\nu,\m^{\epsilon}))_j +
2^{-(\kappa -\epsilon)_j} (m_j +1) \text{ при } j \in \s(\epsilon),
\end{multline*}
т.е. $ x \in U \cap \inter (\supp g_{\kappa -\epsilon, \n_{\epsilon}(\nu,\m^{\epsilon})}^{d,m}), $
а, следовательно, $ \n_{\epsilon}(\nu,\m^{\epsilon}) \in N_{\kappa -\epsilon}^{d,m,U}. $

Для формулировки леммы 1.3.3 при $ d \in \N $ и $ j \in \{1,\ldots,d\} $
обозначим через $ \eta^j: \R^d \times \R^d \mapsto \R^d $
отображение, определяемое соотношением
$$
(\eta^j(\xi,x))_i = \begin{cases} \xi_i, i =1, \ldots, j; \\
x_i, i =j +1, \ldots, d,
\end{cases} \xi,x \in \R^d.
$$
Отметим также, что лемма 1.3.3 устанавливается с помощью леммы 1.3.2. и
её доказательство упрощённо повторяет доказательство леммы 1.2.2 из [10].

Лемма 1.3.3

Пусть $ d \in \N, l \in \Z_+^d, m \in \N^d, \kappa \in \Z_+^d, \epsilon \in
\Upsilon^d: \s(\epsilon) \subset \s(\kappa), U $ -- открытое множество в
$ \R^d. $ Тогда линейный оператор $ H_{\kappa, \kappa -\epsilon}^{d,l,m,U}:
\mathcal P_{\kappa -\epsilon}^{d,l,m,U} \mapsto \mathcal P_\kappa^{d,l,m,U}, $
значение которого для $ f \in \mathcal P_{\kappa -\epsilon}^{d,l,m,U} $
определяется равенством
\begin{equation*} \tag{1.3.16}
H_{\kappa, \kappa -\epsilon}^{d,l,m,U} f = \begin{cases} f, \text{ при }
\epsilon =0; \\
(\prod_{j \in \s(\epsilon)} H_{\eta^j(\kappa -\epsilon, \kappa)}^{j,d,l,m,U}) f,
\text{ при } \epsilon \ne 0, (\text{ см. } (1.3.11)),
\end{cases}
\end{equation*}
обладает следующими свойствами:

1) для $ f \in \mathcal P_{\kappa -\epsilon}^{d,l,m,U} $ выполняется равенство
\begin{equation*} \tag{1.3.17}
(H_{\kappa, \kappa -\epsilon}^{d,l,m,U} f) \mid_{U} = f \mid_{U};
\end{equation*}

2) для $ f \in \mathcal P_{\kappa -\epsilon}^{d,l,m,U} $ вида
\begin{equation*} \tag{1.3.18}
f = \sum_{\nu^\prime \in N_{\kappa -\epsilon}^{d,m,U}}
f_{\kappa -\epsilon, \nu^\prime} g_{\kappa -\epsilon, \nu^\prime}^{d,m},
\{f_{\kappa -\epsilon, \nu^\prime} \in \mathcal P^{d,l},
\nu^\prime \in N_{\kappa -\epsilon}^{d,m,U}\},
\end{equation*}
имеет место представление
\begin{equation*} \tag{1.3.19}
H_{\kappa, \kappa -\epsilon}^{d,l,m,U} f = \sum_{\nu \in N_{\kappa}^{d,m,U}}
f_{\kappa,\nu} g_{\kappa,\nu}^{d,m},
\end{equation*}
где
\begin{equation*} \tag{1.3.20}
f_{\kappa,\nu} = \sum_{\m^{\epsilon} \in \M_{\epsilon}^m(\nu)} A_{\m^{\epsilon}}^m
f_{\kappa -\epsilon, \n_{\epsilon}(\nu,\m^{\epsilon})},
\end{equation*}
а
\begin{equation*} \tag{1.3.21}
A_{\m^{\epsilon}}^m = \prod_{i \in \s(\epsilon)}
a_{\m_i}^{m_i}, (\text{ см. } (1.3.2)), \m^{\epsilon} \in \M_{\epsilon}^m(\nu),
\nu \in N_{\kappa}^{d,m,U}.
\end{equation*}

Отметим ещё частный случай леммы 1.2.3 из [10] (см. (1.3.3), (1.3.21)).

Лемма 1.3.4

При $ d \in \N, \nu \in \Z^d, \epsilon \in \Upsilon^d, m \in \N^d $ имеет место
равенство
\begin{equation*} \tag{1.3.22}
\sum_{\m^{\epsilon} \in \M_{\epsilon}^m(\nu)} A_{\m^{\epsilon}}^m =1.
\end{equation*}
\bigskip

\centerline{\S 2. Основные результаты}
\bigskip

2.1. В этом пункте будут построены средства приближения на некоторых
открытых подмножествах области определения функций из рассматриваемых
нами пространств, которые используются в следующем пункте.

При $ d \in \N, l, \kappa \in \Z_+^d, \nu \in \Z^d $ определим линейный
оператор $ S_{\kappa,\nu}^{d,l}: L_1(Q_{\kappa,\nu}^d) \mapsto \mathcal P^{d,l}, $
полагая $ S_{\kappa, \nu}^{d,l} = P_{\delta, x^0}^{d,l} $ при
$ \delta = 2^{-\kappa}, x^0 = 2^{-\kappa} \nu $ (см. лемму 1.1.2 и (1.3.5)).

При $ d \in \N $ для открытого множества $ D \subset \R^d $ обозначим через
$ L_1^{\loc \square}(D) $ пространство всех функций $ f \in L_1^{\loc}(D), $ для
которых при любых $ \kappa \in \Z_+^d, \nu \in \Z^d $ таких, что
$ Q_{\kappa,\nu}^d \subset D, $ соблюдается включение
$ f \mid_{Q_{\kappa,\nu}^d} \in L_1(Q_{\kappa,\nu}^d). $
Отметим, что в ситуации, когда $ Q_{\kappa,\nu}^d \subset D, $ для $ f \in L_1^{\loc \square}(D) $ вместо
$ S_{\kappa, \nu}^{d,l}(f \mid_{Q_{\kappa, \nu}^d}) $ будем писать
$ S_{\kappa, \nu}^{d,l} f.$

При $ d \in \N $ для области $ D \subset \R^d, $ её открытого подмножества
$ U \subset D $ и $ \kappa \in \Z_+^d $ таких, что множество
\begin{equation*} \tag{2.1.1}
\{\nu^\prime \in \Z^d: Q_{\kappa,\nu^\prime}^d \subset D\} \ne \emptyset,
\end{equation*}
а $ m \in \N^d, $ фиксируем некоторое отображение
\begin{equation*} \tag{2.1.2}
\nu_\kappa = \nu_\kappa^{d,m,D,U}: N_\kappa^{d,m,U} \ni \nu \mapsto
\nu_\kappa^{d,m,D,U}(\nu) \in \{\nu^\prime \in \Z^d:
Q_{\kappa,\nu^\prime}^d \subset D\} \text{ (см. } (1.3.9)),
\end{equation*}
и при $ l \in \Z_+^d $ определим линейный оператор $
E_\kappa^{d,l,m,D,U,\nu_\kappa}: L_1^{\loc \square}(D) \mapsto
\mathcal P_\kappa^{d,l,m,U} $ (см. п. 1.3.) равенством
\begin{equation*} \tag{2.1.3}
E_\kappa^{d,l,m,D,U,\nu_\kappa} f = \sum_{\nu \in N_\kappa^{d,m,U}}
(S_{\kappa, \nu_\kappa^{d,m,D,U}(\nu)}^{d,l} f ) g_{\kappa, \nu}^{d,m},
f \in L_1^{\loc \square}(D).
\end{equation*}

Замечание

Учитывая, что при $ d \in \N, \kappa \in \Z_+^d, \nu \in \Nu_{0, 2^\kappa -\e}^d $
клетка $ (2^{-\kappa} \nu +2^{-\kappa} I^d) \subset I^d, $ а, значит,
\begin{multline*}
Q_{\kappa^0,\nu^0}^d = (2^{-\kappa^0}  \nu^0 +2^{-\kappa^0} I^d) \supset
(2^{-\kappa^0} \nu^0 +2^{-\kappa^0} (2^{-\kappa} \nu +2^{-\kappa} I^d)) = \\
(2^{-\kappa^0 -\kappa} (2^\kappa \nu^0 +\nu) +2^{-\kappa^0 -\kappa} I^d) =
Q_{\kappa^0 +\kappa, 2^\kappa \nu^0 +\nu}^d, \kappa^0 \in \Z_+^d,
\nu^0 \in \Z^d,
\end{multline*}
замечаем, что если при $ d \in \N $ для области $ D \subset \R^d $ и $ \kappa^0
\in \Z_+^d $ соблюдается (2.1.1) при $ \kappa^0 $ вместо $ \kappa, $
то для $ \kappa \in \Z_+^d $ имеет место (2.1.1) при $ \kappa^0 +\kappa $ вместо
$ \kappa. $

Будет полезно следующее утверждение.

Лемма 2.1.1

Пусть $ d \in \N, \lambda \in \Z_+^d, D $ --- область в $ \R^d $ и
функция $ f \in C^\infty(D), $ а $ g \in L_1^{\loc}(D), $ причём
для каждого $ \mu \in \Z_+^d(\lambda) $ (см. (1.1.1)) обобщённая производная
$ \D^\mu g \in L_1^{\loc}(D). $ Тогда в пространстве обобщённых функций в
области $ D $ имеет место соотношение
\begin{equation*} \tag{2.1.4}
\D^\lambda (fg) = \sum_{ \mu \in \Z_+^d(\lambda)} C_\lambda^\mu
\D^{\lambda -\mu} f \D^\mu g \in \L_1^{\loc}(D).
\end{equation*}

Доказательство.

В условиях леммы построим последовательность ограниченных областей
$ D_n \subset \R^d, n \in \N, $ обладающих следующими свойсвами:
$$
\overline D_n \subset D, D_n \subset D_{n +1}, n \in \N, \\
\cup_{n \in \N} D_n = D.
$$
Тогда для каждого $ n \in \N $ имеем:
\begin{equation*}
f \mid_{D_n} \in C^\infty(D_n), \\
\D^\mu (g \mid_{D_n}) = (\D^\mu g) \mid_{D_n} \in L_1(D_n), \mu
\in \Z_+^d(\lambda).
\end{equation*}
Теперь для $ \phi \in C_0^\infty(D) $ вследствие компактности $ \supp \phi \subset D $
фиксируем $ n \in \N, $ для которого $ \supp \phi \subset D_n, $ и, значит,
$ \phi \in C_0^\infty(D_n). $ Тогда применяя лемму 1.3.2 из [10], получаем
\begin{multline*}
\langle \D^\lambda (fg), \phi \rangle =
\langle (\D^\lambda (fg)) \mid_{D_n}, \phi \rangle =
\langle \D^\lambda ((fg) \mid_{D_n}), \phi \rangle =
\langle \D^\lambda (f \mid_{D_n} g \mid_{D_n}), \phi \rangle \\ 
=\int_{D_n} \D^\lambda (f \mid_{D_n} g \mid_{D_n}) \phi dx =
\int_{D_n} (\sum_{ \mu \in \Z_+^d(\lambda)} C_\lambda^\mu
\D^{\lambda -\mu} (f \mid_{D_n}) \D^\mu (g \mid_{D_n})) \phi dx \\
=\sum_{ \mu \in \Z_+^d(\lambda)} C_\lambda^\mu
\int_{D_n} \D^{\lambda -\mu} (f \mid_{D_n}) \D^\mu (g \mid_{D_n}) \phi dx \\
=\sum_{ \mu \in \Z_+^d(\lambda)} C_\lambda^\mu
\int_{D_n} (\D^{\lambda -\mu} f) \mid_{D_n} (\D^\mu g) \mid_{D_n} \phi dx \\
=\sum_{ \mu \in \Z_+^d(\lambda)} C_\lambda^\mu
\int_D \D^{\lambda -\mu} f \D^\mu g \phi dx =
\int_D (\sum_{ \mu \in \Z_+^d(\lambda)} C_\lambda^\mu
\D^{\lambda -\mu} f \D^\mu g) \phi dx,
\end{multline*}
откуда в силу произвольности $ \phi \in C_0^\infty(D) $ следует (2.1.4). $ \square $

Замечание

При $ d \in \N, l \in \Z_+^d, m \in \N^d, \kappa \in \Z_+^d $ для открытого
множества $ U \subset \R^d, $ функции $ f \in \mathcal P_\kappa^{d,l,m,U} $
вида (1.3.10), $ \lambda \in \Z_+^d(m) $ в пространстве обобщённых функций на
$ \R^d $ имеет место соотношение
\begin{equation*} \tag{2.1.5}
\D^\lambda f = \sum_{\nu \in N_\kappa^{d,m,U}} \D^\lambda (f_\nu
g_{\kappa,\nu}^{d,m}) \in L_1^{\loc}(\R^d),
\end{equation*}
где сумма, вообще говоря ряда, в правой части (2.1.5) понимается поточечно.

В самом деле, в описанных условиях, учитывая, что вследствие (1.3.8) при
$ \nu \in N_\kappa^{d,m,U}, \mu \in \Z_+^d(\lambda) $ справедливо включение
$ \D^\mu g_{\kappa, \nu}^{d,m} \in L_\infty(\R^d) \subset \L_1^{\loc}(\R^d), $
согласно (2.1.4) получаем, что
\begin{equation*}
\D^\lambda (f_\nu g_{\kappa,\nu}^{d,m}) \in L_1^{\loc}(\R^d), \nu \in N_\kappa^{d,m,U}.
\end{equation*}
Принимая во внимание это обстоятельство, а также тот факт, что для любого
компакта $ K \subset \R^d $ число
$$
\card \{\nu \in N_\kappa^{d,m,U}: \supp g_{\kappa,\nu}^{d,m} \cap K \ne
\emptyset\} < \infty,
$$
и $ \supp \D^\lambda (f_\nu g_{\kappa,\nu}^{d,m}) \subset \supp g_{\kappa,\nu}^{d,m},
\nu \in N_\kappa^{d,m,U}, $ для $ \phi \in C_0^\infty(\R^d) $ имеем
\begin{multline*}
\langle \D^\lambda f, \phi \rangle = (-1)^{|\lambda|} \langle f, \D^\lambda \phi \rangle =
(-1)^{|\lambda|} \int_{\R^d} (\sum_{\nu \in N_\kappa^{d,m,U}} f_\nu
g_{\kappa,\nu}^{d,m}) \D^\lambda \phi dx \\
= (-1)^{|\lambda|} \int_{\supp \phi}
\sum_{\nu \in N_\kappa^{d,m,U}} f_\nu g_{\kappa,\nu}^{d,m} \D^\lambda \phi dx = \\
(-1)^{|\lambda|} \int_{\supp \phi} \sum_{\nu \in N_\kappa^{d,m,U}: \supp \phi
\cap \supp g_{\kappa,\nu}^{d,m} \ne \emptyset} f_\nu g_{\kappa,\nu}^{d,m} \D^\lambda \phi dx = \\
(-1)^{|\lambda|} \sum_{\nu \in N_\kappa^{d,m,U}: \supp \phi \cap \supp g_{\kappa,\nu}^{d,m}
\ne \emptyset} \int_{\supp \phi} f_\nu g_{\kappa,\nu}^{d,m} \D^\lambda \phi dx = \\
\sum_{\nu \in N_\kappa^{d,m,U}: \supp \phi \cap \supp g_{\kappa,\nu}^{d,m}
\ne \emptyset} (-1)^{|\lambda|} \int_{\supp \phi} f_\nu g_{\kappa,\nu}^{d,m}
\D^\lambda \phi dx = \\
\sum_{\nu \in N_\kappa^{d,m,U}: \supp \phi \cap \supp g_{\kappa,\nu}^{d,m}
\ne \emptyset} (-1)^{|\lambda|} \int_{\R^d} f_\nu g_{\kappa,\nu}^{d,m}
\D^\lambda \phi dx = \\
\sum_{\nu \in N_\kappa^{d,m,U}: \supp \phi \cap \supp g_{\kappa,\nu}^{d,m}
\ne \emptyset} \int_{\R^d} \D^\lambda( f_\nu g_{\kappa,\nu}^{d,m}) \phi dx = \\
\int_{\R^d} \sum_{\nu \in N_\kappa^{d,m,U}: \supp \phi \cap \supp g_{\kappa,\nu}^{d,m}
\ne \emptyset} \D^\lambda( f_\nu g_{\kappa,\nu}^{d,m}) \phi dx = \\
\int_{\R^d} (\sum_{\nu \in N_\kappa^{d,m,U}: \supp \phi \cap \supp g_{\kappa,\nu}^{d,m}
\ne \emptyset} \D^\lambda( f_\nu g_{\kappa,\nu}^{d,m})) \phi dx = \\
\int_{\supp \phi} (\sum_{\nu \in N_\kappa^{d,m,U}: \supp \phi \cap \supp g_{\kappa,\nu}^{d,m}
\ne \emptyset} \D^\lambda( f_\nu g_{\kappa,\nu}^{d,m})) \phi dx = \\
\int_{\supp \phi} (\sum_{\nu \in N_\kappa^{d,m,U}} \D^\lambda( f_\nu
g_{\kappa,\nu}^{d,m})) \phi dx = 
\int_{\R^d} (\sum_{\nu \in N_\kappa^{d,m,U}} \D^\lambda( f_\nu
g_{\kappa,\nu}^{d,m})) \phi dx,
\end{multline*}
откуда следует соотношение (2.1.5).

Предложение 2.1.2

Пусть $ d \in \N, l \in \Z_+^d, m \in \N^d, \lambda \in \Z_+^d(m), $ а область
$ D \subset \R^d $ и её открытое подмножество $ U \subset D $ таковы, что
существуют константы $ \kappa^0 = \kappa^0(d,m,D,U) \in \Z_+^d, \gamma^0 =
\gamma^0(d,m,D,U) \in \R_+^d, $ для которых при любом $ \kappa \in \Z_+^d $
существует отображение $ \nu_{\kappa^0 +\kappa} =
\nu_{\kappa^0 +\kappa}^{d,m,D,U}: N_{\kappa^0 +\kappa}^{d,m,U} \mapsto \Z^d, $
обладающее тем свойством, что для каждого
$ \nu \in N_{\kappa^0 +\kappa}^{d,m,U} $ справедливо включение
\begin{equation*} \tag{2.1.6}
Q_{\kappa^0 +\kappa,\nu_{\kappa^0 +\kappa}^{d,m,D,U}(\nu)}^d \subset D \cap
(2^{-\kappa^0 -\kappa} \nu +\gamma^0 2^{-\kappa^0 -\kappa} B^d).
\end{equation*}
И пусть $ 1 \le p \le q \le \infty, $ а также, если множество $ U $ --
ограниченно, пусть ещё $ 1 \le q < p \le \infty. $ Тогда существует константа
$ c_1(d,l,m,D,U,\lambda,p,q) > 0 $ такая, что для любой функции $ f \in
L_p(D) $ при $ \kappa \in \Z_+^d $ справедливо неравенство
\begin{equation*} \tag{2.1.7}
\| \D^\lambda E_{\kappa^0 +\kappa}^{d,l,m,D,U,\nu_{\kappa^0 +\kappa}} f \|_{L_q(\R^d)} \le
c_1 2^{(\kappa, \lambda +(p^{-1} -q^{-1})_+ \e)} \| f\|_{L_p(D)}.
\end{equation*}

Доказательство.

Сначала отметим свойства некоторых вспомогательных множеств и других
объектов, которые понадобятся для доказательства предложения.

Для $ d \in \N, \kappa \in \Z_+^d, m \in \N^d $ и открытого множества
$ U \subset \R^d $ обозначим множество
$$
G_\kappa^{d,m,U} = \cup_{\nu \in N_\kappa^{d,m,U}} \supp g_{\kappa,\nu}^{d,m},
$$
для которого справедливо представление
\begin{multline*} \tag{2.1.8}
G_\kappa^{d,m,U} = (\cup_{n \in \Z^d: Q_{\kappa,n}^d \cap G_\kappa^{d,m,U} \ne
\emptyset} Q_{\kappa,n}^d) \cup (G_\kappa^{d,m,U} \cap A_\kappa^d), \\
\text{ где } \mes A_\kappa^d =0, A_\kappa^d \cap Q_{\kappa,n}^d = \emptyset,
Q_{\kappa,n}^d \cap Q_{\kappa,n^\prime}^d = \emptyset, n,n^\prime \in
\Z^d: n \ne n^\prime.
\end{multline*}
Понятно, что для $ n \in \Z^d $ такого, что пересечение $ Q_{\kappa,n}^d \cap
G_\kappa^{d,m,U} \ne \emptyset, $ множество
$ \{\nu \in N_\kappa^{d,m,U}: Q_{\kappa,n}^d \cap
\supp g_{\kappa,\nu}^{d,m} \ne \emptyset\} \ne \emptyset, $ при этом, используя
(1.3.4), (1.3.5), легко проверить, что для $ \nu \in N_\kappa^{d,m,U}:
Q_{\kappa,n}^d \cap \supp g_{\kappa,\nu}^{d,m} \ne \emptyset, $ верно
соотношение
\begin{equation*} \tag{2.1.9}
Q_{\kappa,n}^d \subset \overline Q_{\kappa,n}^d \subset
\supp g_{\kappa,\nu}^{d,m} \subset G_\kappa^{d,m,U}.
\end{equation*}

Отметим ещё, что в условиях предложения при $ \kappa \in \Z_+^d $
для $ n \in \Z^d: Q_{\kappa^0 +\kappa,n}^d \cap G_{\kappa^0 +\kappa}^{d,m,U} \ne \emptyset,
\nu \in N_{\kappa^0 +\kappa}^{d,m,U}: \supp g_{\kappa^0 +\kappa, \nu}^{d,m} \cap Q_{\kappa^0 +\kappa,n}^d
\ne \emptyset, $ вследствие (2.1.6) и замкнутости $ B^d, $ а также благодаря
(2.1.9) и (1.3.4) (при $ \kappa^0 +\kappa $ вместо $ \kappa $ ) справедливо
включение
\begin{equation*} \tag{2.1.10}
Q_{\kappa^0 +\kappa,\nu_{\kappa^0 +\kappa}^{d,m,D,U}(\nu)}^d \subset
\overline Q_{\kappa^0 +\kappa,\nu_{\kappa^0 +\kappa}^{d,m,D,U}(\nu)}^d
\subset (2^{-\kappa^0 -\kappa} n +\gamma^1 2^{-\kappa^0 -\kappa} B^d)
\end{equation*}
с $ \gamma^1 = \gamma^1(d,m,D,U) > \e. $

В условиях предложения при $ \kappa \in \Z_+^d, n \in \Z^d:
Q_{\kappa^0 +\kappa,n}^d \cap G_{\kappa^0 +\kappa}^{d,m,U} \ne \emptyset, $
зададим
\begin{equation*} \tag{2.1.11}
x_{\kappa^0 +\kappa,n}^{d,m,D,U} = 2^{-\kappa^0 -\kappa} n -\gamma^1 2^{-\kappa^0 -\kappa}; \\
\delta_{\kappa^0 +\kappa,n}^{d,m,D,U} = 2 \gamma^1 2^{-\kappa^0 -\kappa}
\end{equation*}
и определим клетку $ D_{\kappa^0 +\kappa,n}^{d,m,D,U} $ равенством
\begin{equation*} \tag{2.1.12}
D_{\kappa^0 +\kappa,n}^{d,m,D,U} = x_{\kappa^0 +\kappa,n}^{d,m,D,U} +
\delta_{\kappa^0 +\kappa,n}^{d,m,D,U} I^d = \inter (2^{-\kappa^0 -\kappa} n +\gamma^1 2^{-\kappa^0 -\kappa} B^d).
\end{equation*}

Из приведенных определений с учётом того, что $ \gamma^1 > \e, $ видно, что
соблюдается включение
\begin{equation*} \tag{2.1.13}
Q_{\kappa^0 +\kappa, n}^d \subset D_{\kappa^0 +\kappa,n}^{d,m,D,U},
n \in \Z^d: Q_{\kappa^0 +\kappa,n}^d \cap G_{\kappa^0 +\kappa}^{d,m,U} \ne
\emptyset, \kappa \in \Z_+^d.
\end{equation*}
Учитывая (2.1.13), (2.1.12), нетрудно видеть, что в условиях предложения
существует константа $ c_2(d,m,D,U) >0 $
такая, что при $ \kappa \in \Z_+^d $ для каждого $ x \in \R^d $ число
\begin{equation*} \tag{2.1.14}
\card \{ n \in \Z^d: Q_{\kappa^0 +\kappa,n}^d \cap G_{\kappa^0 +\kappa}^{d,m,U} \ne \emptyset,
x \in D_{\kappa^0 +\kappa,n}^{d,m,D,U} \} \le c_2.
\end{equation*}

Из (2.1.10) и (2.1.12) следует, что при $ \kappa \in \Z_+^d$ для $ n \in \Z^d:
Q_{\kappa^0 +\kappa,n}^d \cap G_{\kappa^0 +\kappa}^{d,m,U} \ne \emptyset, \nu \in N_{\kappa^0 +\kappa}^{d,m,U}:
\supp g_{\kappa^0 +\kappa, \nu}^{d,m} \cap Q_{\kappa^0 +\kappa,n}^d \ne
\emptyset, $ имеет место включение
\begin{equation*} \tag{2.1.15}
Q_{\kappa^0 +\kappa,\nu_{\kappa^0 +\kappa}^{d,m,D,U}(\nu)}^d \subset
D_{\kappa^0 +\kappa,n}^{d,m,D,U}.
\end{equation*}

Из (1.3.6) с учётом (1.3.1) вытекает, что при $ \kappa \in \Z_+^d $ для
$ n \in \Z^d: Q_{\kappa^0 +\kappa,n}^d \cap G_{\kappa^0 +\kappa}^{d,m,U} \ne
\emptyset, $ верно неравенство
\begin{equation*} \tag{2.1.16}
\card \{ \nu \in N_{\kappa^0 +\kappa}^{d,m,U}: \supp g_{\kappa^0 +\kappa, \nu}^{d,m} \cap
Q_{\kappa^0 +\kappa,n}^d \ne \emptyset \} \le c_3(d,m).
\end{equation*}

Теперь в условиях предложения при $ 1 \le p \le q < \infty, \kappa \in \Z_+^d $
для $ f \in L_p(D), $ принимая во внимание (2.1.3), (2.1.5), (2.1.4),
имеем
\begin{multline*} \tag{2.1.17}
\| \D^\lambda E_{\kappa^0 +\kappa}^{d,l,m,D,U,\nu_{\kappa^0 +\kappa}} f \|_{L_q(\R^d)} = \\
\biggl\| \D^\lambda (\sum_{\nu \in N_{\kappa^0 +\kappa}^{d,m,U}}
(S_{\kappa^0 +\kappa, \nu_{\kappa^0 +\kappa}^{d,m,D,U}(\nu)}^{d,l} f)
g_{\kappa^0 +\kappa, \nu}^{d,m}) \biggr\|_{L_q(\R^d)} = \\
\biggl\| \sum_{\nu \in N_{\kappa^0 +\kappa}^{d,m,U}} \D^\lambda
((S_{\kappa^0 +\kappa, \nu_{\kappa^0 +\kappa}^{d,m,D,U}(\nu)}^{d,l} f)
g_{\kappa^0 +\kappa, \nu}^{d,m}) \biggr\|_{L_q(\R^d)} = \\
\biggl\| \sum_{\nu \in N_{\kappa^0 +\kappa}^{d,m,U}} \sum_{ \mu \in \Z_+^d(\lambda)}
C_\lambda^\mu (\D^\mu S_{\kappa^0 +\kappa, \nu_{\kappa^0 +\kappa}^{d,m,D,U}(\nu)}^{d,l} f)
\D^{\lambda -\mu} g_{\kappa^0 +\kappa, \nu}^{d,m} \biggr\|_{L_q(\R^d)} = \\
\biggl\| \sum_{ \mu \in \Z_+^d(\lambda)} C_\lambda^\mu
\sum_{\nu \in N_{\kappa^0 +\kappa}^{d,m,U}}
(\D^\mu S_{\kappa^0 +\kappa, \nu_{\kappa^0 +\kappa}^{d,m,D,U}(\nu)}^{d,l} f)
\D^{\lambda -\mu} g_{\kappa^0 +\kappa, \nu}^{d,m} \biggr\|_{L_q(\R^d)} \le \\
\sum_{ \mu \in \Z_+^d(\lambda)} C_\lambda^\mu
\biggl\| \sum_{\nu \in N_{\kappa^0 +\kappa}^{d,m,U}}
(\D^\mu S_{\kappa^0 +\kappa, \nu_{\kappa^0 +\kappa}^{d,m,D,U}(\nu)}^{d,l} f)
\D^{\lambda -\mu} g_{\kappa^0 +\kappa, \nu}^{d,m} \biggr\|_{L_q(\R^d)}.
\end{multline*}

Оценивая правую часть (2.1.17), при $ \mu \in \Z_+^d(\lambda) $ с учётом
(2.1.8) получаем
\begin{multline*} \tag{2.1.18}
\biggl\| \sum_{\nu \in N_{\kappa^0 +\kappa}^{d,m,U}}
(\D^\mu S_{\kappa^0 +\kappa, \nu_{\kappa^0 +\kappa}^{d,m,D,U}(\nu)}^{d,l} f)
\D^{\lambda -\mu} g_{\kappa^0 +\kappa, \nu}^{d,m} \biggr\|_{L_q(\R^d)}^q = \\
\int_{\R^d} \biggl| \sum_{\nu \in N_{\kappa^0 +\kappa}^{d,m,U}}
(\D^\mu S_{\kappa^0 +\kappa, \nu_{\kappa^0 +\kappa}^{d,m,D,U}(\nu)}^{d,l} f)
\D^{\lambda -\mu} g_{\kappa^0 +\kappa, \nu}^{d,m} \biggr|^q dx = \\
\int_{G_{\kappa^0 +\kappa}^{d,m,U}} \biggl| \sum_{\nu \in N_{\kappa^0 +\kappa}^{d,m,U}}
(\D^\mu S_{\kappa^0 +\kappa, \nu_{\kappa^0 +\kappa}^{d,m,D,U}(\nu)}^{d,l} f)
\D^{\lambda -\mu} g_{\kappa^0 +\kappa, \nu}^{d,m} \biggr|^q dx = \\
\int_{\cup_{n \in \Z^d: Q_{\kappa^0 +\kappa,n}^d \cap
G_{\kappa^0 +\kappa}^{d,m,U} \ne \emptyset} Q_{\kappa^0 +\kappa,n}^d}
\biggl| \sum_{\nu \in N_{\kappa^0 +\kappa}^{d,m,U}}
(\D^\mu S_{\kappa^0 +\kappa, \nu_{\kappa^0 +\kappa}^{d,m,D,U}(\nu)}^{d,l} f)
\D^{\lambda -\mu} g_{\kappa^0 +\kappa, \nu}^{d,m} \biggr|^q dx = \\
\lim_{r \to \infty} \int_{\cup_{\substack{n \in \Z^d: Q_{\kappa^0 +\kappa,n}^d \\ \cap (r B^d)
\cap G_{\kappa^0 +\kappa}^{d,m,U} \ne \emptyset}} Q_{\kappa^0 +\kappa,n}^d}
\biggl| \sum_{\nu \in N_{\kappa^0 +\kappa}^{d,m,U}}
(\D^\mu S_{\kappa^0 +\kappa, \nu_{\kappa^0 +\kappa}^{d,m,D,U}(\nu)}^{d,l} f)
\D^{\lambda -\mu} g_{\kappa^0 +\kappa, \nu}^{d,m} \biggr|^q dx.
\end{multline*}

При оценке правой части (2.1.18) для $ \mu \in \Z_+^d(\lambda), r \in \N $ имеем
\begin{multline*} \tag{2.1.19}
\int_{\cup_{\substack{n \in \Z^d: Q_{\kappa^0 +\kappa,n}^d \\ \cap (r B^d)
\cap G_{\kappa^0 +\kappa}^{d,m,U} \ne \emptyset}} Q_{\kappa^0 +\kappa,n}^d}
\biggl| \sum_{\nu \in N_{\kappa^0 +\kappa}^{d,m,U}}
(\D^\mu S_{\kappa^0 +\kappa, \nu_{\kappa^0 +\kappa}^{d,m,D,U}(\nu)}^{d,l} f)
\D^{\lambda -\mu} g_{\kappa^0 +\kappa, \nu}^{d,m} \biggr|^q dx = \\
\sum_{\substack{n \in \Z^d: Q_{\kappa^0 +\kappa,n}^d \\ \cap (r B^d) \cap
G_{\kappa^0 +\kappa}^{d,m,U} \ne \emptyset}}
\int_{Q_{\kappa^0 +\kappa,n}^d} \biggl| \sum_{\nu \in N_{\kappa^0 +\kappa}^{d,m,U}}
(\D^\mu S_{\kappa^0 +\kappa, \nu_{\kappa^0 +\kappa}^{d,m,D,U}(\nu)}^{d,l} f)
\D^{\lambda -\mu} g_{\kappa^0 +\kappa, \nu}^{d,m} \biggr|^q dx = \\
\sum_{\substack{n \in \Z^d: Q_{\kappa^0 +\kappa,n}^d \\ \cap (r B^d) \cap
G_{\kappa^0 +\kappa}^{d,m,U} \ne \emptyset}}
\int_{Q_{\kappa^0 +\kappa,n}^d} \biggl| \sum_{\substack{\nu \in
N_{\kappa^0 +\kappa}^{d,m,U}: Q_{\kappa^0 +\kappa,n}^d \\ \cap
\supp g_{\kappa^0 +\kappa, \nu}^{d,m} \ne \emptyset}}
(\D^\mu S_{\kappa^0 +\kappa, \nu_{\kappa^0 +\kappa}^{d,m,D,U}(\nu)}^{d,l} f)
\D^{\lambda -\mu} g_{\kappa^0 +\kappa, \nu}^{d,m} \biggr|^q dx = \\
\sum_{\substack{n \in \Z^d: Q_{\kappa^0 +\kappa,n}^d \\ \cap (r B^d) \cap
G_{\kappa^0 +\kappa}^{d,m,U} \ne \emptyset}}
\biggl\| \sum_{\substack{\nu \in N_{\kappa^0 +\kappa}^{d,m,U}: Q_{\kappa^0 +\kappa,n}^d \\
\cap\supp g_{\kappa^0 +\kappa, \nu}^{d,m} \ne \emptyset}}
(\D^\mu S_{\kappa^0 +\kappa, \nu_{\kappa^0 +\kappa}^{d,m,D,U}(\nu)}^{d,l} f)
\D^{\lambda -\mu} g_{\kappa^0 +\kappa, \nu}^{d,m} \biggr\|_{L_q(Q_{\kappa^0 +\kappa,n}^d)}^q \le \\
\sum_{\substack{n \in \Z^d: Q_{\kappa^0 +\kappa,n}^d\\ \cap (r B^d) \cap
G_{\kappa^0 +\kappa}^{d,m,U} \ne \emptyset}}
\biggl(\sum_{\substack{\nu \in N_{\kappa^0 +\kappa}^{d,m,U}: Q_{\kappa^0 +\kappa,n}^d\\ 
\cap\supp g_{\kappa^0 +\kappa, \nu}^{d,m} \ne \emptyset}}
\| (\D^\mu S_{\kappa^0 +\kappa, \nu_{\kappa^0 +\kappa}^{d,m,D,U}(\nu)}^{d,l} f)
\D^{\lambda -\mu} g_{\kappa^0 +\kappa, \nu}^{d,m} \|_{L_q(Q_{\kappa^0 +\kappa,n}^d)}\biggr)^q.
\end{multline*}

Для оценки правой части (2.1.19), используя сначала (1.3.8), а затем
применяя (1.1.3), при $ n \in \Z^d: Q_{\kappa^0 +\kappa,n}^d \cap
G_{\kappa^0 +\kappa}^{d,m,U} \ne \emptyset,
\nu \in N_{\kappa^0 +\kappa}^{d,m,U}: Q_{\kappa^0 +\kappa,n}^d \cap
\supp g_{\kappa^0 +\kappa, \nu}^{d,m} \ne \emptyset, $ выводим
\begin{multline*} \tag{2.1.20}
\| (\D^\mu S_{\kappa^0 +\kappa, \nu_{\kappa^0 +\kappa}^{d,m,D,U}(\nu)}^{d,l} f)
\D^{\lambda -\mu} g_{\kappa^0 +\kappa, \nu}^{d,m} \|_{L_q(Q_{\kappa^0 +\kappa,n}^d)} \le \\
\| \D^{\lambda -\mu} g_{\kappa^0 +\kappa, \nu}^{d,m} \|_{L_\infty(\R^d)}
\| \D^\mu S_{\kappa^0 +\kappa, \nu_{\kappa^0 +\kappa}^{d,m,D,U}(\nu)}^{d,l} f \|_{L_q(Q_{\kappa^0 +\kappa,n}^d)} = \\
c_4 2^{(\kappa^0 +\kappa, \lambda -\mu)}
\| \D^\mu S_{\kappa^0 +\kappa, \nu_{\kappa^0 +\kappa}^{d,m,D,U}(\nu)}^{d,l} f \|_{L_q(Q_{\kappa^0 +\kappa,n}^d)} \le \\
c_5 2^{(\kappa, \lambda -\mu)} 2^{(\kappa, \mu +p^{-1} \e -q^{-1} \e)}
\| S_{\kappa^0 +\kappa, \nu_{\kappa^0 +\kappa}^{d,m,D,U}(\nu)}^{d,l} f\|_{L_p(Q_{\kappa^0 +\kappa,n}^d)} = \\
c_5 2^{(\kappa, \lambda +p^{-1} \e -q^{-1} \e)}
\| S_{\kappa^0 +\kappa, \nu_{\kappa^0 +\kappa}^{d,m,D,U}(\nu)}^{d,l} f\|_{L_p(Q_{\kappa^0 +\kappa,n}^d)}.
\end{multline*}

Далее, принимая во внимание, что в силу (2.1.10) имеет место включение
\begin{multline*}
Q_{\kappa^0 +\kappa,n}^d \subset (2^{-\kappa^0 -\kappa}
\nu_{\kappa^0 +\kappa}^{d,m,D,U}(\nu) +(\gamma^1 +\e) 2^{-\kappa^0 -\kappa} B^d), \\
n \in \Z^d: Q_{\kappa^0 +\kappa,n}^d \cap G_{\kappa^0 +\kappa}^{d,m,U} \ne \emptyset,
\nu \in N_{\kappa^0 +\kappa}^{d,m,U}: \supp g_{\kappa^0 +\kappa, \nu}^{d,m} \cap Q_{\kappa^0 +\kappa,n}^d
\ne \emptyset,
\end{multline*}
на основании (1.1.3), (1.1.5), (2.1.6), (2.1.15) заключаем, что для
$ n \in \Z^d: Q_{\kappa^0 +\kappa,n}^d \cap G_{\kappa^0 +\kappa}^{d,m,U} \ne \emptyset,
\nu \in N_{\kappa^0 +\kappa}^{d,m,U}: \supp g_{\kappa^0 +\kappa, \nu}^{d,m}
\cap Q_{\kappa^0 +\kappa,n}^d \ne \emptyset, $ выполняется неравенство
\begin{multline*} \tag{2.1.21}
\| S_{\kappa^0 +\kappa, \nu_{\kappa^0 +\kappa}^{d,m,D,U}(\nu)}^{d,l} f\|_{L_p(Q_{\kappa^0 +\kappa,n}^d)} \le
c_6 \| S_{\kappa^0 +\kappa, \nu_{\kappa^0 +\kappa}^{d,m,D,U}(\nu)}^{d,l} f
\|_{L_p(Q_{\kappa^0 +\kappa, \nu_{\kappa^0 +\kappa}^{d,m,D,U}(\nu)}^d)} \le \\
c_7 \| f \|_{L_p(Q_{\kappa^0 +\kappa, \nu_{\kappa^0 +\kappa}^{d,m,D,U}(\nu)}^d)} \le c_7
\| f\|_{L_p(D \cap D_{\kappa^0 +\kappa, n}^{d,m,D,U})}.
\end{multline*}

Объединяя (2.1.20) и (2.1.21), находим, что при $ n \in \Z^d: Q_{\kappa^0 +\kappa,n}^d
\cap G_{\kappa^0 +\kappa}^{d,m,U} \ne \emptyset,
\nu \in N_{\kappa^0 +\kappa}^{d,m,U}: Q_{\kappa^0 +\kappa,n}^d \cap \supp g_{\kappa^0 +\kappa, \nu}^{d,m} \ne
\emptyset, $ имеет место неравенство
\begin{equation*}
\| (\D^\mu S_{\kappa^0 +\kappa, \nu_{\kappa^0 +\kappa}^{d,m,D,U}(\nu)}^{d,l} f)
\D^{\lambda -\mu} g_{\kappa^0 +\kappa, \nu}^{d,m} \|_{L_q(Q_{\kappa^0 +\kappa,n}^d)} \le \\
c_8 2^{(\kappa, \lambda +p^{-1} \e -q^{-1} \e)}
\| f\|_{L_p(D \cap D_{\kappa^0 +\kappa,n}^{d,m,D,U})}.
\end{equation*}

Подставляя эту оценку в (2.1.19) и применяя (2.1.16), а затем используя
неравенство (1.1.2) при $ a = p /q \le 1 $ и оценку (2.1.14), приходим к
неравенству
\begin{multline*}
\int_{\cup_{\substack{n \in \Z^d: Q_{\kappa^0 +\kappa,n}^d\\ \cap (r B^d) \cap
G_{\kappa^0 +\kappa}^{d,m,U} \ne \emptyset}} Q_{\kappa^0 +\kappa,n}^d}
\biggl| \sum_{\nu \in N_{\kappa^0 +\kappa}^{d,m,U}}
(\D^\mu S_{\kappa^0 +\kappa, \nu_{\kappa^0 +\kappa}^{d,m,D,U}(\nu)}^{d,l} f)
\D^{\lambda -\mu} g_{\kappa^0 +\kappa, \nu}^{d,m} \biggr|^q dx \le \\
\sum_{\substack{n \in \Z^d: Q_{\kappa^0 +\kappa,n}^d\\ \cap (r B^d) \cap
G_{\kappa^0 +\kappa}^{d,m,U} \ne \emptyset}}
\biggl(\sum_{\substack{\nu \in N_{\kappa^0 +\kappa}^{d,m,U}: \\ Q_{\kappa^0 +\kappa,n}^d \cap
\supp g_{\kappa^0 +\kappa, \nu}^{d,m} \ne \emptyset}}
c_8 2^{(\kappa, \lambda +p^{-1} \e -q^{-1} \e)}
\| f\|_{L_p(D \cap D_{\kappa^0 +\kappa,n}^{d,m,D,U})}\biggr)^q \le \\
(c_8 2^{(\kappa, \lambda +p^{-1} \e -q^{-1} \e)})^q
\sum_{n \in \Z^d: Q_{\kappa^0 +\kappa,n}^d \cap (r B^d) \cap
G_{\kappa^0 +\kappa}^{d,m,U} \ne \emptyset}
\biggl(c_3 \| f\|_{L_p(D \cap D_{\kappa^0 +\kappa,n}^{d,m,D,U})}\biggr)^q = \\
(c_9 2^{(\kappa, \lambda +p^{-1} \e -q^{-1} \e)})^q
\sum_{n \in \Z^d: Q_{\kappa^0 +\kappa,n}^d \cap (r B^d) \cap
G_{\kappa^0 +\kappa}^{d,m,U} \ne \emptyset}
\biggl( \int_{D \cap D_{\kappa^0 +\kappa,n}^{d,m,D,U}} | f(x)|^p dx\biggr)^{q /p} \le \\
(c_9 2^{(\kappa, \lambda +p^{-1} \e -q^{-1} \e)})^q
\biggl(\sum_{n \in \Z^d: Q_{\kappa^0 +\kappa,n}^d \cap (r B^d) \cap
G_{\kappa^0 +\kappa}^{d,m,U} \ne \emptyset}
\int_{D \cap D_{\kappa^0 +\kappa,n}^{d,m,D,U}} | f(x)|^p dx \biggr)^{q /p} = \\
(c_9 2^{(\kappa, \lambda +p^{-1} \e -q^{-1} \e)})^q
\biggl(\sum_{n \in \Z^d: Q_{\kappa^0 +\kappa,n}^d \cap (r B^d) \cap
G_{\kappa^0 +\kappa}^{d,m,U} \ne \emptyset}
\int_D \chi_{ D_{\kappa^0 +\kappa,n}^{d,m,D,U}}(x) | f(x)|^p dx \biggr)^{q /p} = \\
(c_9 2^{(\kappa, \lambda +p^{-1} \e -q^{-1} \e)})^q
\biggl(\int_D \biggl(\sum_{\substack{n \in \Z^d: Q_{\kappa^0 +\kappa,n}^d\\ \cap (r B^d) \cap
G_{\kappa^0 +\kappa}^{d,m,U} \ne \emptyset}}
\chi_{ D_{\kappa^0 +\kappa,n}^{d,m,D,U}}(x)\biggr) | f(x)|^p dx \biggr)^{q /p} \le \\
(c_9 2^{(\kappa, \lambda +p^{-1} \e -q^{-1} \e)})^q
\biggl(\int_D c_2 | f(x)|^p dx \biggr)^{q /p} =\\
(c_{10} 2^{(\kappa, \lambda +p^{-1} \e -q^{-1} \e)} \| f\|_{L_p(D)})^q,
r \in \N, \mu \in \Z_+^d(\lambda).
\end{multline*}

Переходя в последнем неравенстве к пределу при $ r \to \infty, $
в сочетании с (2.1.18) находим, что при $ \mu \in \Z_+^d(\lambda) $ справедливо
неравенство
\begin{multline*} \tag{2.1.22}
\biggl\| \sum_{\nu \in N_{\kappa^0 +\kappa}^{d,m,U}}
(\D^\mu S_{\kappa^0 +\kappa, \nu_{\kappa^0 +\kappa}^{d,m,D,U}(\nu)}^{d,l} f)
\D^{\lambda -\mu} g_{\kappa^0 +\kappa, \nu}^{d,m} \biggr\|_{L_q(\R^d)}^q \le \\
(c_{10} 2^{(\kappa, \lambda +p^{-1} \e -q^{-1} \e)} \|
f\|_{L_p(D)})^q.
\end{multline*}

Соединяя (2.1.17) с (2.1.22). получаем
\begin{multline*}
\| \D^\lambda E_{\kappa^0 +\kappa}^{d,l,m,D,U,\nu_{\kappa^0 +\kappa}} f \|_{L_q(\R^d)} \le \\
\sum_{ \mu \in \Z_+^d(\lambda)} C_\lambda^\mu
c_{10} 2^{(\kappa, \lambda +p^{-1} \e -q^{-1} \e)} \| f\|_{L_p(D)} =
c_1 2^{(\kappa, \lambda +p^{-1} \e -q^{-1} \e)} \| f\|_{L_p(D)},
\end{multline*}
что завершает вывод (2.1.7) для произвольного открытого множества
$ U \subset D $ при $ q \ne \infty. $ Доказательство (2.1.7) при
$ q = \infty $ проводится аналогично.

Если в условиях предложения множество $ U $ -- ограниченно, то при
$ 1 \le q < p \le \infty $ из соблюдения (2.1.7) при $ q = p $ и неравенства
\begin{multline*}
\| \D^\lambda E_{\kappa^0 +\kappa}^{d,l,m,D,U,\nu_{\kappa^0 +\kappa}} f \|_{L_q(\R^d)} =
\| \D^\lambda E_{\kappa^0 +\kappa}^{d,l,m,D,U,\nu_{\kappa^0 +\kappa}} f \|_{L_q(G)} \le \\
(\mes G)^{1 /q -1 /p} \| \D^\lambda E_{\kappa^0 +\kappa}^{d,l,m,D,U,\nu_{\kappa^0 +\kappa}} f \|_{L_p(G)} =
(\mes G)^{1 /q -1 /p} \| \D^\lambda E_{\kappa^0 +\kappa}^{d,l,m,D,U,\nu_{\kappa^0 +\kappa}} f \|_{L_p(\R^d)},
\end{multline*}
где $ G $ -- ограниченная область в $ \R^d $ такая, что $
G_{\kappa^0 +\kappa}^{d,m,U} \subset G $ при $ \kappa \in \Z_+^d,
$ вытекает оценка (2.1.7) в рассматриваемом случае. $ \square $

Обозначим через $ \mathcal I^D $ линейное отображение, которое каждой
функции $ f, $ заданной на множестве $ D \subset \R^d, $ сопоставляет функцию
$ \mathcal I^D f, $ определяемую
на $ \R^d $ равенством
\begin{equation*}
(\mathcal I^D f)(x) = \begin{cases} f(x), \text{ при } x \in D; \\
0, \text{ при } x \in \R^d \setminus D.
\end{cases}
\end{equation*}

Предложение 2.1.3

Пусть $ d \in \N, l \in \N^d, m \in \N^d, 1 \le p < \infty, $
а область $ D \subset \R^d $ и её открытое подмножество $ U \subset D $
удовлетворяют условиям предложения 2.1.2.
Тогда для любой функции $ f \in L_p(D) $ в $ L_p(U) $ имеет место равенство
\begin{equation*} \tag{2.1.23}
f \mid_U = \lim_{\mn(\kappa) \to \infty} (E_{\kappa^0 +\kappa}^{d,l -\e,m,D,U,
\nu_{\kappa^0 +\kappa}} f) \mid_U \text{ (см. } (2.1.3)).
\end{equation*}

Доказательство.

При доказательсве предложения будем использовать объекты и связанные с ними
факты из доказательства предложения 2.1.2.

В условиях предложения определим ещё следующий оператор.
При $ \kappa \in \Z_+^d, n \in \Z^d: Q_{\kappa^0 +\kappa,n}^d \cap
G_{\kappa^0 +\kappa}^{d,m,U} \ne \emptyset, $ введём в рассморение
линейный оператор $ \mathcal S_{\kappa^0 +\kappa,n}^{d,l -\e,m,D,U}, $ полагая
\begin{equation*}
\mathcal S_{\kappa^0 +\kappa,n}^{d,l -\e,m,D,U} =
P_{\delta,x^0}^{d,l -\e} \text{ (см. лемму 1.1.2)}
\end{equation*}
при $ x^0 = x_{\kappa^0 +\kappa,n}^{d,m,D,U},
\delta = \delta_{\kappa^0 +\kappa,n}^{d,m,D,U} \text{ (см. }
(2.1.11), (2.1.12)). $

Пусть теперь в условиях предложения $ f \in L_p(D), \kappa \in \Z_+^d. $
Тогда, полагая $ F = \mathcal I^D f, $ ввиду (2.1.7), (1.3.7), (1.3.9),
(2.1.3), (2.1.6) имеем
\begin{multline*} \tag{2.1.24}
\| f \mid_U  -(E_{\kappa^0 +\kappa}^{d,l -\e,m,D,U,\nu_{\kappa^0 +\kappa}} f) \mid_U
\|_{L_p(U)}^p = \\
\| F \mid_U -(E_{\kappa^0 +\kappa}^{d,l -\e,m,D,U,
\nu_{\kappa^0 +\kappa}} f) \mid_U \|_{L_p(U)}^p = \\
\biggl\| (F (\sum_{ \nu \in \Z^d: \supp g_{\kappa^0 +\kappa, \nu}^{d,m} \cap U \ne \emptyset}
g_{\kappa^0 +\kappa, \nu}^{d,m})) \mid_U -(E_{\kappa^0 +\kappa}^{d,l -\e,m,D,U,
\nu_{\kappa^0 +\kappa}} f) \mid_U \biggr\|_{L_p(U)}^p \le \\
\biggl\| F (\sum_{ \nu \in N_{\kappa^0 +\kappa}^{d,m,U}}
g_{\kappa^0 +\kappa, \nu}^{d,m}) -\sum_{\nu \in N_{\kappa^0 +\kappa}^{d,m,U}}
(S_{\kappa^0 +\kappa, \nu_{\kappa^0 +\kappa}^{d,m,D,U}(\nu)}^{d,l -\e} f)
g_{\kappa^0 +\kappa, \nu}^{d,m} \biggr\|_{L_p(\R^d)}^p = \\
\int_{\R^d} \biggl| \sum_{\nu \in N_{\kappa^0 +\kappa}^{d,m,U}} (F(x) -
(S_{\kappa^0 +\kappa, \nu_{\kappa^0 +\kappa}^{d,m,D,U}(\nu)}^{d,l -\e} F)(x))
g_{\kappa^0 +\kappa, \nu}^{d,m}(x) \biggr|^p dx = \\
\int_{G_{\kappa^0 +\kappa}^{d,m,U}} \biggl| \sum_{\nu \in
N_{\kappa^0 +\kappa}^{d,m,U}} (F(x) -
(S_{\kappa^0 +\kappa, \nu_{\kappa^0 +\kappa}^{d,m,D,U}(\nu)}^{d,l -\e} F)(x))
g_{\kappa^0 +\kappa, \nu}^{d,m}(x) \biggr|^p dx.
\end{multline*}

В силу (2.1.8) выводим
\begin{multline*} \tag{2.1.25}
\int_{G_{\kappa^0 +\kappa}^{d,m,U}} \biggl| \sum_{\nu \in
N_{\kappa^0 +\kappa}^{d,m,U}} (F(x) -
(S_{\kappa^0 +\kappa, \nu_{\kappa^0 +\kappa}^{d,m,D,U}(\nu)}^{d,l -\e} F)(x))
g_{\kappa^0 +\kappa, \nu}^{d,m}(x) \biggr|^p dx = \\
\sum_{n \in \Z^d: Q_{\kappa^0 +\kappa,n}^d \cap
G_{\kappa^0 +\kappa}^{d,m,U} \ne \emptyset}
\int_{Q_{\kappa^0 +\kappa,n}^d} \biggl| \sum_{\nu \in N_{\kappa^0 +\kappa}^{d,m,U}} (F(x) -
(S_{\kappa^0 +\kappa, \nu_{\kappa^0 +\kappa}^{d,m,D,U}(\nu)}^{d,l -\e} F)(x))
g_{\kappa^0 +\kappa, \nu}^{d,m}(x) \biggr|^p dx = \\
\sum_{n \in \Z^d: Q_{\kappa^0 +\kappa,n}^d \cap G_{\kappa^0 +\kappa}^{d,m,U} \ne \emptyset}
\biggl\| \sum_{\substack{\nu \in N_{\kappa^0 +\kappa}^{d,m,U}:
\supp g_{\kappa^0 +\kappa, \nu}^{d,m}\\ \cap Q_{\kappa^0 +\kappa,n}^d \ne \emptyset}} 
(F - (S_{\kappa^0 +\kappa, \nu_{\kappa^0 +\kappa}^{d,m,D,U}(\nu)}^{d,l -\e} F))
g_{\kappa^0 +\kappa, \nu}^{d,m} \biggr\|_{L_p(Q_{\kappa^0 +\kappa,n}^d)}^p \le \\
\sum_{n \in \Z^d: Q_{\kappa^0 +\kappa,n}^d \cap G_{\kappa^0 +\kappa}^{d,m,U} \ne \emptyset}
\biggl( \sum_{\substack{\nu \in N_{\kappa^0 +\kappa}^{d,m,U}:
\supp g_{\kappa^0 +\kappa, \nu}^{d,m}\\ \cap Q_{\kappa^0 +\kappa,n}^d \ne \emptyset}}
\| (F -(S_{\kappa^0 +\kappa, \nu_{\kappa^0 +\kappa}^{d,m,D,U}(\nu)}^{d,l -\e} F))
g_{\kappa^0 +\kappa, \nu}^{d,m} \|_{L_p(Q_{\kappa^0 +\kappa,n}^d)} \biggr)^p.
\end{multline*}

Далее, для $ n \in \Z^d: Q_{\kappa^0 +\kappa,n}^d \cap G_{\kappa^0 +\kappa}^{d,m,U} \ne \emptyset,
\nu \in N_{\kappa^0 +\kappa}^{d,m,U}: \supp g_{\kappa^0 +\kappa, \nu}^{d,m}
\cap Q_{\kappa^0 +\kappa,n}^d \ne \emptyset, $ с учётом (1.3.8) получаем
\begin{multline*} \tag{2.1.26}
\| (F -(S_{\kappa^0 +\kappa, \nu_{\kappa^0 +\kappa}^{d,m,D,U}(\nu)}^{d,l -\e} F))
g_{\kappa^0 +\kappa, \nu}^{d,m} \|_{L_p(Q_{\kappa^0 +\kappa,n}^d)} \le \\
\| F -S_{\kappa^0 +\kappa, \nu_{\kappa^0 +\kappa}^{d,m,D,U}(\nu)}^{d,l -\e} F \|_{L_p(Q_{\kappa^0 +\kappa,n}^d)} \le \\
\| F -\mathcal S_{\kappa^0 +\kappa,n}^{d,l -\e,m,D,U} F \|_{L_p(Q_{\kappa^0 +\kappa,n}^d)} +
\| \mathcal S_{\kappa^0 +\kappa,n}^{d,l -\e,m,D,U} F
-S_{\kappa^0 +\kappa, \nu_{\kappa^0 +\kappa}^{d,m,D,U}(\nu)}^{d, l-\e} F \|_{L_p(Q_{\kappa^0 +\kappa,n}^d)}.
\end{multline*}

Учитывая (2.1.13), на основании (1.1.6), (2.1.12), (2.1.11) заключаем, что
для $ n \in \Z^d: Q_{\kappa^0 +\kappa,n}^d \cap
G_{\kappa^0 +\kappa}^{d,m,U} \ne \emptyset, $ выполняется неравенство
\begin{multline*} \tag{2.1.27}
\| F -\mathcal S_{\kappa^0 +\kappa,n}^{d,l -\e,m,D,U} F
\|_{L_p(Q_{\kappa^0 +\kappa,n}^d)} \le \| F -\mathcal S_{\kappa^0
+\kappa,n}^{d,l -\e,m,D,U} F
\|_{L_p(D_{\kappa^0 +\kappa,n}^{d,m,D,U})} \le \\
c_{11} \sum_{j =1}^d 2^{(\kappa^0_j +\kappa_j) /p} \biggl(\int_{ c_{12} 2^{-\kappa^0_j -\kappa_j} B^1}
\int_{ (D_{\kappa^0 +\kappa,n}^{d,m,D,U})_{l_j \xi e_j}}
|\Delta_{\xi e_j}^{l_j} F(x)|^p dx d\xi \biggr)^{1/p} \le \\
c_{13} \sum_{j =1}^d 2^{\kappa_j /p} \biggl(\int_{ c_{12} 2^{-\kappa_j} B^1}
\int_{ (D_{\kappa^0 +\kappa,n}^{d,m,D,U})_{l_j \xi e_j}}
|\Delta_{\xi e_j}^{l_j} F(x)|^p dx d\xi \biggr)^{1/p}.
\end{multline*}

Принимая во внимание (2.1.13), (2.1.15), (1.1.3), (1.1.4), (1.1.5)
и снова (2.1.15), а также, используя (2.1.27), для $ n \in \Z^d:
Q_{\kappa^0 +\kappa,n}^d \cap G_{\kappa^0 +\kappa}^{d,m,U} \ne
\emptyset, \nu \in N_{\kappa^0 +\kappa}^{d,m,U}:
\supp g_{\kappa^0 +\kappa, \nu}^{d,m} \cap Q_{\kappa^0 +\kappa,n}^d \ne \emptyset, $
находим, что
\begin{multline*} \tag{2.1.28}
\| \mathcal S_{\kappa^0 +\kappa,n}^{d,l -\e,m,D,U} F
-S_{\kappa^0 +\kappa,\nu_{\kappa^0 +\kappa}^{d,m,D,U}(\nu)}^{d,l -\e} F \|_{L_p(Q_{\kappa^0 +\kappa,n}^d)} \le \\
c_{14} \| S_{\kappa^0 +\kappa, \nu_{\kappa^0 +\kappa}^{d,m,D,U}(\nu)}^{d, l -\e}
(\mathcal S_{\kappa^0 +\kappa,n}^{d,l -\e,m,D,U} F -F)
\|_{L_p(Q_{\kappa^0 +\kappa, \nu_{\kappa^0 +\kappa}^{d,m,D,U}(\nu)}^d)} \le \\
c_{15} \| F -\mathcal S_{\kappa^0 +\kappa,n}^{d,l -\e,m,D,U} F \|_{L_p(D_{\kappa^0 +\kappa,n}^{d,m,D,U})} \le \\
c_{16} \sum_{j =1}^d 2^{\kappa_j /p} \biggl(\int_{ c_{12} 2^{-\kappa_j} B^1}
\int_{ (D_{\kappa^0 +\kappa,n}^{d,m,D,U})_{l_j \xi e_j}}
|\Delta_{\xi e_j}^{l_j} F(x)|^p dx d\xi \biggr)^{1/p}.
\end{multline*}

Соединяя (2.1.26), (2.1.27) и (2.1.28), для $ n \in \Z^d: Q_{\kappa^0 +\kappa,n}^d
\cap G_{\kappa^0 +\kappa}^{d,m,U} \ne \emptyset, \nu \in N_{\kappa^0 +\kappa}^{d,m,U}:
\supp g_{\kappa^0 +\kappa, \nu}^{d,m} \cap Q_{\kappa^0 +\kappa,n}^d \ne \emptyset, $ получаем
\begin{multline*}
\| (F -S_{\kappa^0 +\kappa, \nu_{\kappa^0 +\kappa}^{d,m,D,U}(\nu)}^{d, l-\e} F)
g_{\kappa^0 +\kappa,\nu}^{d,m} \|_{L_p(Q_{\kappa^0 +\kappa,n}^d)} \le \\
c_{17} \sum_{j =1}^d 2^{\kappa_j /p} \biggl(\int_{ c_{12} 2^{-\kappa_j} B^1}
\int_{ (D_{\kappa^0 +\kappa, n}^{d,m,D,U})_{l_j \xi e_j}}
|\Delta_{\xi e_j}^{l_j} F(x)|^p dx d\xi \biggr)^{1/p}.
\end{multline*}

Подставляя эту оценку в (2.1.25) и применяя (2.1.16) и неравенство
Гёльдера, а затем используя (2.1.14) и (1.1.2) при $ a =1/p, $
выводим
\begin{multline*}
\int_{G_{\kappa^0 +\kappa}^{d,m,U}} \biggl| \sum_{\nu \in N_{\kappa^0 +\kappa}^{d,m,U}} (F(x) -
(S_{\kappa^0 +\kappa, \nu_{\kappa^0 +\kappa}^{d,m,D,U}(\nu)}^{d,l -\e} F)(x))
g_{\kappa^0 +\kappa, \nu}^{d,m}(x) \biggr|^p dx \le \\
\sum_{n \in \Z^d: Q_{\kappa^0 +\kappa,n}^d \cap G_{\kappa^0 +\kappa}^{d,m,U} \ne \emptyset}
c_{18}^p \biggl(\sum_{j =1}^d 2^{\kappa_j} \int_{ c_{12} 2^{-\kappa_j} B^1}
\int_{ (D_{\kappa^0 +\kappa,n}^{d,m,D,U})_{l_j \xi e_j}}
|\Delta_{\xi e_j}^{l_j} F(x)|^p dx d\xi \biggr) \le \\
c_{18}^p \sum_{j =1}^d 2^{\kappa_j } \int_{ c_{12} 2^{-\kappa_j} B^1}
\biggl(\sum_{n \in \Z^d: Q_{\kappa^0 +\kappa,n}^d \cap G_{\kappa^0 +\kappa}^{d,m,U} \ne \emptyset}
\int_{D_{\kappa^0 +\kappa,n}^{d,m,D,U}} |\Delta_{\xi e_j}^{l_j} F(x)|^p dx \biggr) d\xi = \\
c_{18}^p \sum_{j =1}^d 2^{\kappa_j } \int_{ c_{12} 2^{-\kappa_j} B^1}
\int_{ \R^d} \biggl(\sum_{n \in \Z^d: Q_{\kappa^0 +\kappa,n}^d \cap G_{\kappa^0 +\kappa}^{d,m,U} \ne
\emptyset} \chi_{D_{\kappa^0 +\kappa,n}^{d,m,D,U}}(x) \biggr)
|\Delta_{\xi e_j}^{l_j} F(x)|^p dx d\xi \le \\
\biggl(c_{19} \sum_{j =1}^d 2^{\kappa_j /p} \biggl(\int_{c_{12} 2^{-\kappa_j} B^1}
\int_{ \R^d} |\Delta_{\xi e_j}^{l_j} F(x)|^p dx d\xi \biggr)^{1/p} \biggr)^p.
\end{multline*}
Объединяя последнее неравенство с (2.1.24) и учитывая (1.1.7), для
$ f \in L_p(D) $ при $ \kappa \in \Z_+^d $ приходим к соотношению
\begin{multline*}
\| f \mid_U -(E_{\kappa^0 +\kappa}^{d,l -\e,m,D,U,
\nu_{\kappa^0 +\kappa}} f) \mid_U \|_{L_p(U)} \le \\
c_{20} \sum_{j =1}^d 2^{\kappa_j / p} \biggl(\int_{c_{12} 2^{-\kappa_j} B^1}
\int_{ \R^d} |\Delta_{\xi e_j}^{l_j} F(x)|^p dx d\xi \biggr)^{1/p} \le \\
c_{21} \sum_{j=1}^d \Omega^{l e_j} (\mathcal I^D f, c_{12} 2^{-\kappa_j} )_{L_p(\R^d)} \to 0
\text{ при } \mn(\kappa) \to \infty \text{ (см., например, [20])},
\end{multline*}
что влечёт (2.1.23). $ \square $

Для формулировки следующего утверждения введём обозначение.

Пусть $ d \in \N, D -- $ область в $ \R^d $ и $ \kappa^0 \in \Z_+^d $ таковы,
что выполняется (2.1.1) с $ \kappa^0 $ вместо $ \kappa, $ а $ U \subset D $ --
открытое подмножество $ D. $ Имея в виду замечание после (2.1.3),
при $ m \in \N^d $ рассмотрим некоторое семейство отображений
$ \Nu = \{ \nu_{\kappa^0 +\kappa}^{d,m,D,U},
\kappa \in \Z_+^d \}, $ вида (2.1.2) с $ \kappa^0 +\kappa $ вместо $ \kappa $ и
при $ \kappa, l \in \Z_+^d, $ исходя из (2.1.3) и (1.3.16), определим линейный
оператор $ \mathcal E_{\kappa^0, \kappa}^{d,l,m,D,U,\Nu}: L_1^{\loc \square}(D)
\mapsto \mathcal P_{\kappa^0 +\kappa}^{d,l,m,U}, $ полагая
\begin{equation*} \tag{2.1.29}
\mathcal E_{\kappa^0, \kappa}^{d,l,m,D,U,\Nu} = \sum_{\epsilon
\in \Upsilon^d: \s(\epsilon) \subset \s(\kappa)} (-\e)^\epsilon
H_{\kappa^0 +\kappa, \kappa^0 +\kappa -\epsilon}^{d,l,m,U}
E_{\kappa^0 +\kappa -\epsilon}^{d,l,m,D,U,\nu_{\kappa^0 +\kappa -\epsilon}}.
\end{equation*}
При этом, принимая во внимание (2.1.29), (2.1.3), (1.3.18), (1.3.19), (1.3.20),
имеем
\begin{multline*} \tag{2.1.30}
\mathcal E_{\kappa^0, \kappa}^{d,l,m,D,U,\Nu} f = \\
\sum_{\epsilon \in \Upsilon^d:
\s(\epsilon) \subset \s(\kappa)} (-\e)^\epsilon
\sum_{\nu \in N_{\kappa^0 +\kappa}^{d,m,U}}
\biggl(\sum_{\m^{\epsilon} \in \M_{\epsilon}^m(\nu)} A_{\m^{\epsilon}}^m
S_{\kappa^0 +\kappa -\epsilon, \nu_{\kappa^0 +\kappa -\epsilon}^{d,m,D,U}
(\n_{\epsilon}(\nu,\m^{\epsilon}))}^{d,l} f\biggr)
g_{\kappa^0 +\kappa,\nu}^{d,m} = \\
\sum_{\epsilon \in \Upsilon^d: \s(\epsilon) \subset \s(\kappa)}
\sum_{\nu \in N_{\kappa^0 +\kappa}^{d,m,U}} \biggl((-\e)^\epsilon
\sum_{\m^{\epsilon} \in \M_{\epsilon}^m(\nu)} A_{\m^{\epsilon}}^m
S_{\kappa^0 +\kappa -\epsilon, \nu_{\kappa^0 +\kappa -\epsilon}^{d,m,D,U}
(\n_{\epsilon}(\nu,\m^{\epsilon}))}^{d,l} f\biggr)
g_{\kappa^0 +\kappa,\nu}^{d,m} = \\
\sum_{\nu \in N_{\kappa^0 +\kappa}^{d,m,U}}
\sum_{\epsilon \in \Upsilon^d: \s(\epsilon) \subset \s(\kappa)} \biggl((-\e)^\epsilon
\sum_{\m^{\epsilon} \in \M_{\epsilon}^m(\nu)} A_{\m^{\epsilon}}^m
S_{\kappa^0 +\kappa -\epsilon, \nu_{\kappa^0 +\kappa -\epsilon}^{d,m,D,U}
(\n_{\epsilon}(\nu,\m^{\epsilon}))}^{d,l} f\biggr)
g_{\kappa^0 +\kappa,\nu}^{d,m} = \\
\sum_{\nu \in N_{\kappa^0 +\kappa}^{d,m,U}}
\biggl(\sum_{\epsilon \in \Upsilon^d: \s(\epsilon) \subset \s(\kappa)} (-\e)^\epsilon
\sum_{\m^{\epsilon} \in \M_{\epsilon}^m(\nu)} A_{\m^{\epsilon}}^m
S_{\kappa^0 +\kappa -\epsilon, \nu_{\kappa^0 +\kappa -\epsilon}^{d,m,D,U}
(\n_{\epsilon}(\nu,\m^{\epsilon}))}^{d,l} f\biggr)
g_{\kappa^0 +\kappa,\nu}^{d,m} = \\
\sum_{ \nu \in N_{\kappa^0 +\kappa}^{d,m,U}}
(U_{\kappa^0 +\kappa,\nu}^{d,l,m,D,U,\Nu} f) g_{\kappa^0 +\kappa, \nu}^{d,m},
f \in L_1^{\loc \square}(D),
\end{multline*}
где $ U_{\kappa^0 +\kappa,\nu}^{d,l,m,D,U,\Nu}: L_1^{\loc \square}(D) \mapsto
\mathcal P^{d,l} $ -- линейный оператор, определяемый при $ d \in \N, D \subset
\R^d, U \subset D, l \in \Z_+^d, m \in \N^d, \kappa \in \Z_+^d, \nu \in
N_{\kappa^0 +\kappa}^{d,m,U}, \Nu = \{\nu_{\kappa^0 +\kappa}, \kappa \in \Z_+^d \}, $
удовлетворяющих указанным выше условиям, равенством
\begin{multline*} \tag{2.1.31}
U_{\kappa^0 +\kappa,\nu}^{d,l,m,D,U,\Nu} f =
\sum_{\epsilon \in \Upsilon^d: \s(\epsilon) \subset \s(\kappa)} (-\e)^\epsilon
\sum_{\m^{\epsilon} \in \M_{\epsilon}^m(\nu)} A_{\m^{\epsilon}}^m
S_{\kappa^0 +\kappa -\epsilon, \nu_{\kappa^0 +\kappa -\epsilon}^{d,m,D,U}
(\n_{\epsilon}(\nu,\m^{\epsilon}))}^{d,l} f,\\ f \in L_1^{\loc \square}(D).
\end{multline*}

Определение

При $ d \in \N, m \in \N^d $ будем говорить, что область $ D \subset \R^d $
и её открытое подмножество $ U \subset D $ являются $m$-правильной парой, если
существуют константы $ \Kappa^0 = \Kappa^0(d,m,D,U) \in \Z_+^d, \Gamma^0 =
\Gamma^0(d,m,D,U) \in \R_+^d, $ для которых существуют семейства отображений
$$
\Nu = \Nu^{d,m,D,U} = \{ \nu_{\Kappa^0 +\kappa}^{d,m,D,U}: N_{\Kappa^0 +\kappa}^{d,m,U}
\mapsto \Z^d, \kappa \in \Z_+^d\}, \\
\{n_{\Kappa^0 +\kappa}^{d,m,D,U}: N_{\Kappa^0 +\kappa}^{d,m,U} \mapsto \Z^d,
\kappa \in \Z_+^d\},
$$
обладающих следующими свойствами:

1) при $ \kappa \in \Z_+^d $ для каждого $ \nu \in N_{\Kappa^0 +\kappa}^{d,m,U} $
справедливо включение
\begin{equation*} \tag{2.1.32}
(Q_{\Kappa^0 +\kappa,\nu_{\Kappa^0 +\kappa}^{d,m,D,U}(\nu)}^d \cup
Q_{\Kappa^0 +\kappa,n_{\Kappa^0 +\kappa}^{d,m,D,U}(\nu)}^d) \subset D \cap
(2^{-\Kappa^0 -\kappa} \nu +\Gamma^0 2^{-\Kappa^0 -\kappa} B^d);
\end{equation*}

2) при $ \kappa \in \Z_+^d, \nu \in N_{\Kappa^0 +\kappa}^{d,m,U},
\epsilon \in \Upsilon^d: \s(\epsilon) \subset \s(\kappa), \m^{\epsilon} \in
\M_{\epsilon}^m(\nu) $ для
$$
\mathcal D_{\Kappa^0 +\kappa,\nu,\epsilon,\m^{\epsilon}}^{d,m,D,U} =
\bm x_{\Kappa^0 +\kappa,\nu,\epsilon,\m^{\epsilon}}^{d,m,D,U} +
\bm \delta_{\Kappa^0 +\kappa,\nu,\epsilon,\m^{\epsilon}}^{d,m,D,U} I^d,
$$
где точка $ \bm x_{\Kappa^0 +\kappa,\nu,\epsilon,\m^{\epsilon}}^{d,m,D,U} \in \R^d $
и вектор $ \bm \delta_{\Kappa^0 +\kappa,\nu,\epsilon,\m^{\epsilon}}^{d,m,D,U} \in \R_+^d $
определяются равенствами
\begin{multline*}
(\bm x_{\Kappa^0 +\kappa,\nu,\epsilon,\m^{\epsilon}}^{d,m,D,U})_j =\\
\min(2^{-\Kappa^0_j -\kappa_j} (n_{\Kappa^0 +\kappa}^{d,m,D,U}(\nu))_j,
2^{-\Kappa^0_j -\kappa_j +\epsilon_j} (\nu_{\Kappa^0 +\kappa -\epsilon}^{d,m,D,U}
(\n_{\epsilon}(\nu,\m^{\epsilon})))_j),  j \in \Nu_{1,d}^1; \\
(\bm \delta_{\Kappa^0 +\kappa,\nu,\epsilon,\m^{\epsilon}}^{d,m,D,U})_j =
\max(2^{-\Kappa^0_j -\kappa_j} (n_{\Kappa^0 +\kappa}^{d,m,D,U}(\nu))_j
+2^{-\Kappa^0_j -\kappa_j}, \\
2^{-\Kappa^0_j -\kappa_j +\epsilon_j} (\nu_{\Kappa^0 +\kappa -\epsilon}^{d,m,D,U}
(\n_{\epsilon}(\nu,\m^{\epsilon})))_j +2^{-\Kappa^0_j -\kappa_j +\epsilon_j})
-(\bm x_{\Kappa^0 +\kappa,\nu,\epsilon,\m^{\epsilon}}^{d,m,D,U})_j,
j \in \Nu_{1,d}^1,
\end{multline*}
соблюдается включение
\begin{equation*} \tag{2.1.33}
\mathcal D_{\Kappa^0 +\kappa,\nu,\epsilon,\m^{\epsilon}}^{d,m,D,U} \subset D;
\end{equation*}

3) при $ \kappa \in \Z_+^d, \nu \in N_{\Kappa^0 +\kappa}^{d,m,U} $ для любых
$ \epsilon \in \Upsilon^d: \s(\epsilon) \subset \s(\kappa), $ и $ \m^{\epsilon}
\in \M_{\epsilon}^m(\nu) $ при $ j \in \Nu_{1,d}^1 \setminus \s(\epsilon) $
выполняется равенство
\begin{equation*} \tag{2.1.34}
(\nu_{\Kappa^0 +\kappa -\epsilon}^{d,m,D,U}
(\n_{\epsilon}(\nu,\m^{\epsilon})))_j =
(\nu_{\Kappa^0 +\kappa}^{d,m,D,U}(\nu))_j.
\end{equation*}

Замечание

В условиях приведенного определения  существует константа $ \Gamma^1(d,m,D,U)
\in \R_+^d $ такая, что при $ \kappa \in \Z_+^d, \nu \in N_{\Kappa^0 +\kappa}^{d,m,U},
\epsilon \in \Upsilon^d: \s(\epsilon) \subset \s(\kappa), \m^{\epsilon} \in
\M_{\epsilon}^m(\nu) $ имеют место следующие соотношения:

\begin{equation*} \tag{2.1.35}
2^{-\Kappa^0 -\kappa} \le \bm \delta_{\Kappa^0 +\kappa,\nu,\epsilon,\m^{\epsilon}}^{d,m,D,U}
\le \Gamma^1 2^{-\Kappa^0 -\kappa},
\end{equation*}
\begin{equation*} \tag{2.1.36}
Q_{\Kappa^0 +\kappa, n_{\Kappa^0 +\kappa}^{d,m,D,U}(\nu)}^d \cup
Q_{\Kappa^0 +\kappa -\epsilon, \nu_{\Kappa^0 +\kappa -\epsilon}^{d,m,D,U}
(\n_{\epsilon}(\nu,\m^{\epsilon}))}^d \subset
\mathcal D_{\Kappa^0 +\kappa,\nu,\epsilon,\m^{\epsilon}}^{d,m,D,U}.
\end{equation*}

При проверке (2.1.35) с одной стороны имеем
\begin{multline*}
(\bm \delta_{\Kappa^0 +\kappa,\nu,\epsilon,\m^{\epsilon}}^{d,m,D,U})_j \ge
2^{-\Kappa^0_j -\kappa_j} (n_{\Kappa^0 +\kappa}^{d,m,D,U}(\nu))_j +2^{-\Kappa^0_j -\kappa_j} -\\
2^{-\Kappa^0_j -\kappa_j} (n_{\Kappa^0 +\kappa}^{d,m,D,U}(\nu))_j =
2^{-\Kappa^0_j -\kappa_j}, j \in \Nu_{1,d}^1,
\end{multline*}
т.е. справедливо первое неравенство в (2.1.35).

С другой стороны при $ j \in \Nu_{1,d}^1 $ соблюдаются соотношения
\begin{equation*}
2^{-\Kappa^0_j -\kappa_j} (n_{\Kappa^0 +\kappa}^{d,m,D,U}(\nu))_j +2^{-\Kappa^0_j -\kappa_j} -
2^{-\Kappa^0_j -\kappa_j} (n_{\Kappa^0 +\kappa}^{d,m,D,U}(\nu))_j =
2^{-\Kappa^0_j -\kappa_j};
\end{equation*}
\begin{multline*}
2^{-\Kappa^0_j -\kappa_j +\epsilon_j} (\nu_{\Kappa^0 +\kappa -\epsilon}^{d,m,D,U}
(\n_{\epsilon}(\nu,\m^{\epsilon})))_j +2^{-\Kappa^0_j -\kappa_j +\epsilon_j} - \\
2^{-\Kappa^0_j -\kappa_j +\epsilon_j} (\nu_{\Kappa^0 +\kappa -\epsilon}^{d,m,D,U}
(\n_{\epsilon}(\nu,\m^{\epsilon})))_j = 2^{-\Kappa^0_j -\kappa_j +\epsilon_j};
\end{multline*}
\begin{multline*}
2^{-\Kappa^0_j -\kappa_j +\epsilon_j} (\nu_{\Kappa^0 +\kappa -\epsilon}^{d,m,D,U}
(\n_{\epsilon}(\nu,\m^{\epsilon})))_j +2^{-\Kappa^0_j -\kappa_j +\epsilon_j}
-2^{-\Kappa^0_j -\kappa_j} (n_{\Kappa^0 +\kappa}^{d,m,D,U}(\nu))_j \le \\
| 2^{-\Kappa^0_j -\kappa_j +\epsilon_j} (\nu_{\Kappa^0 +\kappa -\epsilon}^{d,m,D,U}
(\n_{\epsilon}(\nu,\m^{\epsilon})))_j
-2^{-\Kappa^0_j -\kappa_j} (n_{\Kappa^0 +\kappa}^{d,m,D,U}(\nu))_j |
+2^{-\Kappa^0_j -\kappa_j +\epsilon_j};
\end{multline*}
\begin{multline*}
2^{-\Kappa^0_j -\kappa_j} (n_{\Kappa^0 +\kappa}^{d,m,D,U}(\nu))_j
+2^{-\Kappa^0_j -\kappa_j}
-2^{-\Kappa^0_j -\kappa_j +\epsilon_j} (\nu_{\Kappa^0 +\kappa -\epsilon}^{d,m,D,U}
(\n_{\epsilon}(\nu,\m^{\epsilon})))_j \le \\
| 2^{-\Kappa^0_j -\kappa_j +\epsilon_j} (\nu_{\Kappa^0 +\kappa -\epsilon}^{d,m,D,U}
(\n_{\epsilon}(\nu,\m^{\epsilon})))_j
-2^{-\Kappa^0_j -\kappa_j} (n_{\Kappa^0 +\kappa}^{d,m,D,U}(\nu))_j |
+2^{-\Kappa^0_j -\kappa_j};
\end{multline*}
причём, учитывая, что в рассматриваемой ситуации справедливо включение
$ \s(\epsilon) \subset \s(\Kappa^0 +\kappa), $ в силу (1.3.15)
(с $ \Kappa^0 +\kappa $ вместо $ \kappa $) и (2.1.32) имеет место неравенство
\begin{multline*}
| 2^{-\Kappa^0_j -\kappa_j +\epsilon_j} (\nu_{\Kappa^0 +\kappa -\epsilon}^{d,m,D,U}
(\n_{\epsilon}(\nu,\m^{\epsilon})))_j
-2^{-\Kappa^0_j -\kappa_j} (n_{\Kappa^0 +\kappa}^{d,m,D,U}(\nu))_j | = \\
| 2^{-\Kappa^0_j -\kappa_j +\epsilon_j} (\nu_{\Kappa^0 +\kappa -\epsilon}^{d,m,D,U}
(\n_{\epsilon}(\nu,\m^{\epsilon})))_j
-2^{-\Kappa^0_j -\kappa_j +\epsilon_j}
(\n_{\epsilon}(\nu,\m^{\epsilon}))_j +\\
2^{-\Kappa^0_j -\kappa_j +\epsilon_j}
(\n_{\epsilon}(\nu,\m^{\epsilon}))_j
-2^{-\Kappa^0_j -\kappa_j} \nu_j +2^{-\Kappa^0_j -\kappa_j} \nu_j
-2^{-\Kappa^0_j -\kappa_j} (n_{\Kappa^0 +\kappa}^{d,m,D,U}(\nu))_j | \le \\
| 2^{-\Kappa^0_j -\kappa_j +\epsilon_j} (\nu_{\Kappa^0 +\kappa -\epsilon}^{d,m,D,U}
(\n_{\epsilon}(\nu,\m^{\epsilon})))_j
-2^{-\Kappa^0_j -\kappa_j +\epsilon_j}
(\n_{\epsilon}(\nu,\m^{\epsilon}))_j | +\\
| 2^{-\Kappa^0_j -\kappa_j +\epsilon_j}
(\n_{\epsilon}(\nu,\m^{\epsilon}))_j
-2^{-\Kappa^0_j -\kappa_j} \nu_j |
+| 2^{-\Kappa^0_j -\kappa_j} \nu_j
-2^{-\Kappa^0_j -\kappa_j} (n_{\Kappa^0 +\kappa}^{d,m,D,U}(\nu))_j | \le \\
\Gamma^0_j 2^{-\Kappa^0_j -\kappa_j +\epsilon_j}
+| 2^{-\Kappa^0_j -\kappa_j +\epsilon_j} (\n_{\epsilon}(\nu,\m^{\epsilon}))_j
-2^{-\Kappa^0_j -\kappa_j} \nu_j | +\Gamma^0_j 2^{-\Kappa^0_j -\kappa_j},
\end{multline*}
а вследствие (1.3.13) соблюдаются соотношения
\begin{equation*}
| 2^{-\Kappa^0_j -\kappa_j +\epsilon_j} (\n_{\epsilon}(\nu,\m^{\epsilon}))_j
-2^{-\Kappa^0_j -\kappa_j} \nu_j | =
| 2^{-\Kappa^0_j -\kappa_j} \nu_j -2^{-\Kappa^0_j -\kappa_j} \nu_j | = 0 \text{ при }
j \in \Nu_{1,d}^1 \setminus \s(\epsilon),
\end{equation*}
\begin{multline*}
| 2^{-\Kappa^0_j -\kappa_j +\epsilon_j} (\n_{\epsilon}(\nu,\m^{\epsilon}))_j
-2^{-\Kappa^0_j -\kappa_j} \nu_j | =
| 2^{-\Kappa^0_j -\kappa_j +1} (\nu_j -\m_j) /2 -2^{-\Kappa^0_j -\kappa_j} \nu_j | = \\
| 2^{-\Kappa^0_j -\kappa_j} (\nu_j -\m_j) -2^{-\Kappa^0_j -\kappa_j} \nu_j | =
| -\m_j 2^{-\Kappa^0_j -\kappa_j} | \le (m_j +1) 2^{-\Kappa^0_j -\kappa_j}
\text{ при } j \in \s(\epsilon),
\end{multline*}
сопоставление которых приводит к заключению о справедливости второго неравенства
в (2.1.35).

Справедливость (2.1.36) следует из того, что для $ x \in
Q_{\Kappa^0 +\kappa, n_{\Kappa^0 +\kappa}^{d,m,D,U}(\nu)}^d \cup
Q_{\Kappa^0 +\kappa -\epsilon, \nu_{\Kappa^0 +\kappa -\epsilon}^{d,m,D,U}
(\n_{\epsilon}(\nu,\m^{\epsilon}))}^d $ при $ j \in \Nu_{1,d}^1 $ выполняются
неравенства
\begin{multline*}
\min(2^{-\Kappa^0_j -\kappa_j} (n_{\Kappa^0 +\kappa}^{d,m,D,U}(\nu))_j,
2^{-\Kappa^0_j -\kappa_j +\epsilon_j} (\nu_{\Kappa^0 +\kappa -\epsilon}^{d,m,D,U}
(\n_{\epsilon}(\nu,\m^{\epsilon})))_j) < x_j < \\
\max(2^{-\Kappa^0_j -\kappa_j} (n_{\Kappa^0 +\kappa}^{d,m,D,U}(\nu))_j +2^{-\Kappa^0_j -\kappa_j},
2^{-\Kappa^0_j -\kappa_j +\epsilon_j} (\nu_{\Kappa^0 +\kappa -\epsilon}^{d,m,D,U}
(\n_{\epsilon}(\nu,\m^{\epsilon})))_j +2^{-\Kappa^0_j -\kappa_j +\epsilon_j}).
\end{multline*}

Отметим ещё, что при $ d \in \N $ для $ m, \bm m \in \N^d: m \le \bm m, $
всякая $ \bm m $-правильная пара $ (D,U)$ является $m$-правильной парой.

Предложение 2.1.4

Пусть $ d \in \N, m \in \N^d, $ а область $ D \subset \R^d $ и её открытое
подмножество $ U \subset D $ являются $m$-правильной парой. И пусть
$ l \in \N^d, \lambda \in \Z_+^d(m), 1 \le p < \infty, p \le q \le \infty, $ а
также, если множество $ U $ -- ограниченно, пусть ещё $ 1 \le q < p. $
Тогда существуют константы $ c_{22}(d,m,D,U,l,\lambda,p,q) > 0,
c_{23}(d,m,D,U) > 0 $ такие, что при $ \kappa \in \Z_+^d \setminus \{0\} $ для
$ f \in L_p(D) $ соблюдается неравенство
\begin{multline*} \tag{2.1.37}
\| \D^\lambda \mathcal E_{\Kappa^0,\kappa}^{d,l -\e,m,D,U,\Nu} f \|_{L_q(\R^d)}
\le c_{22} 2^{(\kappa, \lambda +(p^{-1} -q^{-1})_+ \e)}
\Omega^{\prime l \chi_{\s(\kappa)}}(f, (c_{23} 2^{-\kappa})^{\s(\kappa)})_{L_p(D)} \\
\text{(см. (2.1.29) с $ \Kappa^0, \Nu $ из определения $m$-правильной пары).}
\end{multline*}

Доказательство.

В условиях предложения пусть $ 1 \le p \le q < \infty, f \in L_p(D), \kappa \in
\Z_+^d \setminus \{0\}. $ Тогда, принимая во внимание (2.1.30), (2.1.5),
(2.1.4), имеем
\begin{multline*} \tag{2.1.38}
\| \D^\lambda \mathcal E_{\Kappa^0,\kappa}^{d,l -\e,m,D,U,\Nu} f \|_{L_q(\R^d)} = \\
\biggl\| \D^\lambda \biggl(\sum_{ \nu \in N_{\Kappa^0 +\kappa}^{d,m,U}}
(U_{\Kappa^0 +\kappa,\nu}^{d,l -\e,m,D,U,\Nu} f) g_{\Kappa^0 +\kappa, \nu}^{d,m}\biggr) \biggr\|_{L_q(\R^d)} = \\
\biggl\| \sum_{ \nu \in N_{\Kappa^0 +\kappa}^{d,m,U}}
\D^\lambda \biggl((U_{\Kappa^0 +\kappa,\nu}^{d,l -\e,m,D,U,\Nu} f)
g_{\Kappa^0 +\kappa, \nu}^{d,m} \biggr) \biggr\|_{L_q(\R^d)} = \\
\biggl\| \sum_{ \nu \in N_{\Kappa^0 +\kappa}^{d,m,U}}
\sum_{ \mu \in \Z_+^d(\lambda)} C_\lambda^\mu
\D^\mu (U_{\Kappa^0 +\kappa,\nu}^{d,l -\e,m,D,U,\Nu} f)
\D^{\lambda -\mu} g_{\Kappa^0 +\kappa, \nu}^{d,m} \biggr\|_{L_q(\R^d)} = \\
\biggl\| \sum_{ \mu \in \Z_+^d(\lambda)} C_\lambda^\mu
\sum_{ \nu \in N_{\Kappa^0 +\kappa}^{d,m,U}}
\D^\mu (U_{\Kappa^0 +\kappa,\nu}^{d,l -\e,m,D,U,\Nu} f)
\D^{\lambda -\mu} g_{\Kappa^0 +\kappa, \nu}^{d,m} \biggr\|_{L_q(\R^d)} \le \\
\sum_{ \mu \in \Z_+^d(\lambda)} C_\lambda^\mu
\biggl\| \sum_{ \nu \in N_{\Kappa^0 +\kappa}^{d,m,U}}
\D^\mu (U_{\Kappa^0 +\kappa,\nu}^{d,l -\e,m,D,U,\Nu} f)
\D^{\lambda -\mu} g_{\Kappa^0 +\kappa, \nu}^{d,m} \biggr\|_{L_q(\R^d)}.
\end{multline*}

Оценивая правую часть (2.1.38), при $ \mu \in \Z_+^d(\lambda) $ с учётом
(2.1.8) получаем
\begin{multline*} \tag{2.1.39}
\biggl\| \sum_{\nu \in N_{\Kappa^0 +\kappa}^{d,m,U}}
\D^\mu (U_{\Kappa^0 +\kappa,\nu}^{d,l -\e,m,D,U,\Nu} f)
\D^{\lambda -\mu} g_{\Kappa^0 +\kappa, \nu}^{d,m} \biggr\|_{L_q(\R^d)}^q = \\
\int_{\R^d} \biggl| \sum_{\nu \in N_{\Kappa^0 +\kappa}^{d,m,U}}
\D^\mu (U_{\Kappa^0 +\kappa,\nu}^{d,l -\e,m,D,U,\Nu} f)
\D^{\lambda -\mu} g_{\Kappa^0 +\kappa, \nu}^{d,m} \biggr|^q dx = \\
\int_{G_{\Kappa^0 +\kappa}^{d,m,U}} \biggl| \sum_{\nu \in N_{\Kappa^0 +\kappa}^{d,m,U}}
\D^\mu (U_{\Kappa^0 +\kappa,\nu}^{d,l -\e,m,D,U,\Nu} f)
\D^{\lambda -\mu} g_{\Kappa^0 +\kappa, \nu}^{d,m} \biggr|^q dx = \\
\int\limits_{\bigcup_{n \in \Z^d: Q_{\Kappa^0 +\kappa,n}^d \cap
G_{\Kappa^0 +\kappa}^{d,m,U} \ne \emptyset} Q_{\Kappa^0 +\kappa,n}^d}
\biggl| \sum_{\nu \in N_{\Kappa^0 +\kappa}^{d,m,U}}
\D^\mu (U_{\Kappa^0 +\kappa,\nu}^{d,l -\e,m,D,U,\Nu} f)
\D^{\lambda -\mu} g_{\Kappa^0 +\kappa, \nu}^{d,m} \biggr|^q dx = \\
\lim_{r \to \infty} \int\limits_{\bigcup{n \in \Z^d: Q_{\Kappa^0 +\kappa,n}^d \cap (r B^d)
\cap G_{\Kappa^0 +\kappa}^{d,m,U} \ne \emptyset} Q_{\Kappa^0 +\kappa,n}^d}
\biggl| \sum_{\nu \in N_{\Kappa^0 +\kappa}^{d,m,U}}
\D^\mu (U_{\Kappa^0 +\kappa,\nu}^{d,l -\e,m,D,U,\Nu} f)
\D^{\lambda -\mu} g_{\Kappa^0 +\kappa, \nu}^{d,m} \biggr|^q dx.
\end{multline*}

Для оценки правой части (2.1.39) при $ \mu \in \Z_+^d(\lambda), r \in \N $ имеем
\begin{multline*} \tag{2.1.40}
\int\limits_{\bigcup_{n \in \Z^d: Q_{\Kappa^0 +\kappa,n}^d \cap (r B^d) \cap
G_{\Kappa^0 +\kappa}^{d,m,U} \ne \emptyset} Q_{\Kappa^0
+\kappa,n}^d} \biggl| \sum_{\nu \in N_{\Kappa^0 +\kappa}^{d,m,U}}
\D^\mu (U_{\Kappa^0 +\kappa,\nu}^{d,l -\e,m,D,U,\Nu} f)
\D^{\lambda -\mu} g_{\Kappa^0 +\kappa, \nu}^{d,m} \biggr|^q dx = \\
\sum_{\substack{n \in \Z^d: Q_{\Kappa^0 +\kappa,n}^d \\ \cap (r B^d) \cap
G_{\Kappa^0 +\kappa}^{d,m,U} \ne \emptyset}}
\int_{Q_{\Kappa^0 +\kappa,n}^d} \biggl| \sum_{\substack{\nu \in
N_{\Kappa^0 +\kappa}^{d,m,U}:\\ Q_{\Kappa^0 +\kappa,n}^d \cap
\supp g_{\Kappa^0 +\kappa, \nu}^{d,m} \ne \emptyset}}
\D^\mu (U_{\Kappa^0 +\kappa,\nu}^{d,l -\e,m,D,U,\Nu} f)
\D^{\lambda -\mu} g_{\Kappa^0 +\kappa, \nu}^{d,m} \biggr|^q dx \le \\
\sum_{\substack{n \in \Z^d: Q_{\Kappa^0 +\kappa,n}^d \\ \cap (r B^d) \cap
G_{\Kappa^0 +\kappa}^{d,m,U} \ne \emptyset}}
\biggl(\sum_{\substack{\nu \in N_{\Kappa^0 +\kappa}^{d,m,U}: \\ Q_{\Kappa^0 +\kappa,n}^d \cap
\supp g_{\Kappa^0 +\kappa, \nu}^{d,m} \ne \emptyset}}
\biggl\| \D^\mu (U_{\Kappa^0 +\kappa,\nu}^{d,l -\e,m,D,U,\Nu} f)
\D^{\lambda -\mu} g_{\Kappa^0 +\kappa, \nu}^{d,m} \biggr\|_{L_q(Q_{\Kappa^0 +\kappa,n}^d)} \biggr)^q.
\end{multline*}

Для проведения оценки правой части (2.1.40) приведём некоторые полезные для
нас факты. Прежде всего, в условиях предложения при $ \kappa \in \Z_+^d $
для $ \nu \in N_{\Kappa^0 +\kappa}^{d,m,U} $ с учётом (1.3.5),
(2.1.32) имеем
\begin{equation*}
(2^{-\Kappa^0 -\kappa} n_{\Kappa^0 +\kappa}^{d,m,D,U}(\nu) +
2^{-\Kappa^0 -\kappa} I^d) \subset
(2^{-\Kappa^0 -\kappa} \nu +\Gamma^0 2^{-\Kappa^0 -\kappa} B^d),
\end{equation*}
а вследствие замкнутости $ B^d $ заключаем, что
\begin{equation*}
(2^{-\Kappa^0 -\kappa} n_{\Kappa^0 +\kappa}^{d,m,D,U}(\nu) +
2^{-\Kappa^0 -\kappa} \overline I^d) \subset
(2^{-\Kappa^0 -\kappa} \nu +\Gamma^0 2^{-\Kappa^0 -\kappa} B^d),
\end{equation*}
и, в частности,
\begin{equation*}
2^{-\Kappa^0 -\kappa} n_{\Kappa^0 +\kappa}^{d,m,D,U}(\nu) \in
(2^{-\Kappa^0 -\kappa} \nu +\Gamma^0 2^{-\Kappa^0 -\kappa} B^d),
\end{equation*}
или
\begin{equation*}
(2^{-\Kappa^0 -\kappa} n_{\Kappa^0 +\kappa}^{d,m,D,U}(\nu) -
2^{-\Kappa^0 -\kappa} \nu) \in \Gamma^0 2^{-\Kappa^0 -\kappa} B^d,
\end{equation*}
откуда ввиду симметричности $ B^d $ относительно $0$ приходим к выводу, что
\begin{equation*} \tag{2.1.41}
(2^{-\Kappa^0 -\kappa} \nu -2^{-\Kappa^0 -\kappa}
n_{\Kappa^0 +\kappa}^{d,m,D,U}(\nu)) \in \Gamma^0 2^{-\Kappa^0 -\kappa} B^d.
\end{equation*}

Кроме того, при $ \kappa \in \Z_+^d $ для $ n \in \Z^d:
Q_{\Kappa^0 +\kappa,n}^d \cap G_{\Kappa^0 +\kappa}^{d,m,U} \ne \emptyset, $
и $ \nu \in N_{\Kappa^0 +\kappa}^{d,m,U}: Q_{\Kappa^0 +\kappa,n}^d
\cap \supp g_{\Kappa^0 +\kappa, \nu}^{d,m} \ne \emptyset, $ учитывая
(2.1.9), (1.3.4), заключаем, что
\begin{equation*}
Q_{\Kappa^0 +\kappa, n}^d \subset 2^{-\Kappa^0 -\kappa} \nu +
2^{-\Kappa^0 -\kappa} (m +\e) \overline I^d \subset 2^{-\Kappa^0 -\kappa} \nu +
2^{-\Kappa^0 -\kappa} (m +\e) B^d
\end{equation*}
или
\begin{equation*}
(Q_{\Kappa^0 +\kappa, n}^d -2^{-\Kappa^0 -\kappa} \nu) \subset
2^{-\Kappa^0 -\kappa} (m +\e) B^d.
\end{equation*}
Из последнего включения и (2.1.41) вытекает, что при $ \kappa \in \Z_+^d $
для $ n \in \Z^d: Q_{\Kappa^0 +\kappa,n}^d \cap
G_{\Kappa^0 +\kappa}^{d,m,U} \ne \emptyset, $
и $ \nu \in N_{\Kappa^0 +\kappa}^{d,m,U}: Q_{\Kappa^0 +\kappa,n}^d
\cap \supp g_{\Kappa^0 +\kappa, \nu}^{d,m} \ne \emptyset, $ справедливо
соотношение
\begin{multline*}
Q_{\Kappa^0 +\kappa, n}^d -2^{-\Kappa^0 -\kappa} n_{\Kappa^0 +\kappa}^{d,m,D,U}(\nu) =\\
(Q_{\Kappa^0 +\kappa, n}^d -2^{-\Kappa^0 -\kappa} \nu +2^{-\Kappa^0 -\kappa} \nu
-2^{-\Kappa^0 -\kappa} n_{\Kappa^0 +\kappa}^{d,m,D,U}(\nu)) \subset \\
2^{-\Kappa^0 -\kappa} (m +\e) B^d +\Gamma^0 2^{-\Kappa^0 -\kappa} B^d \subset \\
(m +\e +\Gamma^0) 2^{-\Kappa^0 -\kappa} B^d =
\Gamma^2(d,m,D,U) 2^{-\Kappa^0 -\kappa} B^d
\end{multline*}
с $ \Gamma^2 = m +\e +\Gamma^0 > \e, $
или
\begin{equation*} \tag{2.1.42}
Q_{\Kappa^0 +\kappa, n}^d \subset
2^{-\Kappa^0 -\kappa} n_{\Kappa^0 +\kappa}^{d,m,D,U}(\nu) +
\Gamma^2 2^{-\Kappa^0 -\kappa} B^d,
\end{equation*}
а также
\begin{equation*}
2^{-\Kappa^0 -\kappa} n +2^{-\Kappa^0 -\kappa} \overline I^d =
\overline Q_{\Kappa^0 +\kappa, n}^d \subset
2^{-\Kappa^0 -\kappa} n_{\Kappa^0 +\kappa}^{d,m,D,U}(\nu) +
\Gamma^2 2^{-\Kappa^0 -\kappa} B^d,
\end{equation*}
в частности,
\begin{equation*}
2^{-\Kappa^0 -\kappa} n \in
2^{-\Kappa^0 -\kappa} n_{\Kappa^0 +\kappa}^{d,m,D,U}(\nu) +
\Gamma^2 2^{-\Kappa^0 -\kappa} B^d,
\end{equation*}
или
\begin{equation*}
(2^{-\Kappa^0 -\kappa} n -
2^{-\Kappa^0 -\kappa} n_{\Kappa^0 +\kappa}^{d,m,D,U}(\nu)) \in
\Gamma^2 2^{-\Kappa^0 -\kappa} B^d,
\end{equation*}
а, значит,
\begin{equation*} \tag{2.1.43}
(2^{-\Kappa^0 -\kappa} n_{\Kappa^0 +\kappa}^{d,m,D,U}(\nu) -
2^{-\Kappa^0 -\kappa} n) \in \Gamma^2 2^{-\Kappa^0 -\kappa} B^d.
\end{equation*}

Далее, при $ \kappa \in \Z_+^d $
для $ \nu \in N_{\Kappa^0 +\kappa}^{d,m,U}, \epsilon \in \Upsilon^d:
\s(\epsilon) \subset \s(\kappa), \m^{\epsilon} \in \M_{\epsilon}^m(\nu) $
с учётом (2.1.35) имеем
\begin{equation*}
\mathcal D_{\Kappa^0 +\kappa,\nu,\epsilon,\m^{\epsilon}}^{d,m,D,U} \subset
\bm x_{\Kappa^0 +\kappa,\nu,\epsilon,\m^{\epsilon}}^{d,m,D,U} +
\Gamma^1 2^{-\Kappa^0 -\kappa} I^d \subset
\bm x_{\Kappa^0 +\kappa,\nu,\epsilon,\m^{\epsilon}}^{d,m,D,U} +
\Gamma^1 2^{-\Kappa^0 -\kappa} B^d,
\end{equation*}
а, следовательно,
\begin{equation*}
\overline {\mathcal D}_{\Kappa^0 +\kappa,\nu,\epsilon,\m^{\epsilon}}^{d,m,D,U} \subset
\bm x_{\Kappa^0 +\kappa,\nu,\epsilon,\m^{\epsilon}}^{d,m,D,U} +
\Gamma^1 2^{-\Kappa^0 -\kappa} B^d,
\end{equation*}
или
\begin{equation*} \tag{2.1.44}
\overline {\mathcal D}_{\Kappa^0 +\kappa,\nu,\epsilon,\m^{\epsilon}}^{d,m,D,U} -
\bm x_{\Kappa^0 +\kappa,\nu,\epsilon,\m^{\epsilon}}^{d,m,D,U} \subset
\Gamma^1 2^{-\Kappa^0 -\kappa} B^d,
\end{equation*}
кроме того, принимая во внимание (1.3.5), (2.1.36), получаем, что
\begin{equation*}
2^{-\Kappa^0 -\kappa} n_{\Kappa^0 +\kappa}^{d,m,D,U}(\nu) \in
\overline Q_{\Kappa^0 +\kappa, n_{\Kappa^0 +\kappa}^{d,m,D,U}(\nu)}^d \subset
\overline {\mathcal D}_{\Kappa^0 +\kappa,\nu,\epsilon,\m^{\epsilon}}^{d,m,D,U},
\end{equation*}
что в силу (2.1.44) приводит к включению
\begin{equation*}
2^{-\Kappa^0 -\kappa} n_{\Kappa^0 +\kappa}^{d,m,D,U}(\nu) -
\bm x_{\Kappa^0 +\kappa,\nu,\epsilon,\m^{\epsilon}}^{d,m,D,U} \in
\Gamma^1 2^{-\Kappa^0 -\kappa} B^d,
\end{equation*}
а, значит,
\begin{equation*} \tag{2.1.45}
\bm x_{\Kappa^0 +\kappa,\nu,\epsilon,\m^{\epsilon}}^{d,m,D,U} -
2^{-\Kappa^0 -\kappa} n_{\Kappa^0 +\kappa}^{d,m,D,U}(\nu) \in
\Gamma^1 2^{-\Kappa^0 -\kappa} B^d,
\end{equation*}
наконец, на основании (2.1.44), (2.1.45) вытекает, что
\begin{multline*} \tag{2.1.46}
\mathcal D_{\Kappa^0 +\kappa,\nu,\epsilon,\m^{\epsilon}}^{d,m,D,U} -
2^{-\Kappa^0 -\kappa} n_{\Kappa^0 +\kappa}^{d,m,D,U}(\nu) \subset
(\overline {\mathcal D}_{\Kappa^0 +\kappa,\nu,\epsilon,\m^{\epsilon}}^{d,m,D,U} -
\bm x_{\Kappa^0 +\kappa,\nu,\epsilon,\m^{\epsilon}}^{d,m,D,U}) +\\
(\bm x_{\Kappa^0 +\kappa,\nu,\epsilon,\m^{\epsilon}}^{d,m,D,U} -
2^{-\Kappa^0 -\kappa} n_{\Kappa^0 +\kappa}^{d,m,D,U}(\nu)) \subset \\
\Gamma^1 2^{-\Kappa^0 -\kappa} B^d +\Gamma^1 2^{-\Kappa^0 -\kappa} B^d
\subset 2 \Gamma^1 2^{-\Kappa^0 -\kappa} B^d.
\end{multline*}
Таким образом, при $ \kappa \in \Z_+^d $
для $ n \in \Z^d: Q_{\Kappa^0 +\kappa,n}^d \cap G_{\Kappa^0 +\kappa}^{d,m,U} \ne \emptyset, $
и $ \nu \in N_{\Kappa^0 +\kappa}^{d,m,U}: Q_{\Kappa^0 +\kappa,n}^d
\cap \supp g_{\Kappa^0 +\kappa, \nu}^{d,m} \ne \emptyset, \epsilon \in
\Upsilon^d: \s(\epsilon) \subset \s(\kappa), \m^{\epsilon} \in
\M_{\epsilon}^m(\nu) $ согласно (2.1.46), (2.1.43) соблюдается включение
\begin{multline*}
\mathcal D_{\Kappa^0 +\kappa,\nu,\epsilon,\m^{\epsilon}}^{d,m,D,U} -
2^{-\Kappa^0 -\kappa} n =
(\mathcal D_{\Kappa^0 +\kappa,\nu,\epsilon,\m^{\epsilon}}^{d,m,D,U} -
2^{-\Kappa^0 -\kappa} n_{\Kappa^0 +\kappa}^{d,m,D,U}(\nu)) +\\
(2^{-\Kappa^0 -\kappa} n_{\Kappa^0 +\kappa}^{d,m,D,U}(\nu) -
2^{-\Kappa^0 -\kappa} n) \subset \\
2 \Gamma^1 2^{-\Kappa^0 -\kappa} B^d +\Gamma^2 2^{-\Kappa^0 -\kappa} B^d \subset
\Gamma^3(d,m,D,U) 2^{-\Kappa^0 -\kappa} B^d
\end{multline*}
с $ \Gamma^3 = 2 \Gamma^1 +\Gamma^2 > \e, $ которое в соединении с (2.1.33) даёт
\begin{equation*}
\mathcal D_{\Kappa^0 +\kappa,\nu,\epsilon,\m^{\epsilon}}^{d,m,D,U} \subset
D \cap (2^{-\Kappa^0 -\kappa} n +\Gamma^3 2^{-\Kappa^0 -\kappa} B^d)
\end{equation*}
или
\begin{equation*} \tag{2.1.47}
\mathcal D_{\Kappa^0 +\kappa,\nu,\epsilon,\m^{\epsilon}}^{d,m,D,U} \subset
D \cap D_{\Kappa^0 +\kappa,n}^{\prime d,m,D,U},
\end{equation*}

где
$$
D_{\Kappa^0 +\kappa,n}^{\prime d,m,D,U} = 2^{-\Kappa^0 -\kappa} n +\Gamma^3
2^{-\Kappa^0 -\kappa} B^d,
\kappa \in \Z_+^d, n \in \Z^d: Q_{\Kappa^0 +\kappa,n}^d \cap
G_{\Kappa^0 +\kappa}^{d,m,U} \ne \emptyset.
$$

Из приведенных определений с учётом того, что $ \Gamma^3 > \e, $ видно, что
при $ \kappa \in \Z_+^d, n \in \Z^d: Q_{\Kappa^0 +\kappa,n}^d \cap
G_{\Kappa^0 +\kappa}^{d,m,U} \ne \emptyset, $ справедливо включение
\begin{equation*} \tag{2.1.48}
Q_{\Kappa^0 +\kappa, n}^d \subset D_{\Kappa^0 +\kappa, n}^{\prime d,m,D,U}.
\end{equation*}

Используя (2.1.48), легко проверить, что существует константа
$ c_{24}(d,m,D,U) >0 $ такая, что при $ \kappa \in \Z_+^d $ для каждого
$ x \in \R^d $ число
\begin{equation*} \tag{2.1.49}
\card \{ n \in \Z^d: Q_{\Kappa^0 +\kappa,n}^d \cap G_{\Kappa^0 +\kappa}^{d,m,U} \ne \emptyset,
x \in D_{\Kappa^0 +\kappa,n}^{\prime d,m,D,U} \} \le c_{24}.
\end{equation*}

Перейдём к оценке правой части (2.1.40). При $ n \in \Z^d:
Q_{\Kappa^0 +\kappa,n}^d \cap G_{\Kappa^0 +\kappa}^{d,m,U} \ne \emptyset,
\nu \in N_{\Kappa^0 +\kappa}^{d,m,U}:  Q_{\Kappa^0 +\kappa,n}^d \cap
\supp g_{\Kappa^0 +\kappa, \nu}^{d,m} \ne \emptyset, \mu \in \Z_+^d(\lambda), $
используя (1.3.8), а также благодаря (2.1.31), (2.1.42), применяя (1.1.3), выводим
\begin{multline*} \tag{2.1.50}
\biggl\| \D^\mu (U_{\Kappa^0 +\kappa,\nu}^{d,l -\e,m,D,U,\Nu} f)
\D^{\lambda -\mu} g_{\Kappa^0 +\kappa, \nu}^{d,m} \biggr\|_{L_q(Q_{\Kappa^0 +\kappa,n}^d)} \le \\
\| \D^{\lambda -\mu} g_{\Kappa^0 +\kappa, \nu}^{d,m} \|_{L_\infty(\R^d)}
\biggl\| \D^\mu (U_{\Kappa^0 +\kappa,\nu}^{d,l -\e,m,D,U,\Nu} f) \biggr\|_{L_q(Q_{\Kappa^0 +\kappa,n}^d)} \le \\
c_{25}(d,m,D,U) 2^{(\kappa, \lambda -\mu)}
\biggl\| \D^\mu (U_{\Kappa^0 +\kappa,\nu}^{d,l -\e,m,D,U,\Nu} f) \biggr\|_{L_q(Q_{\Kappa^0 +\kappa,n}^d)} \le \\
c_{25} 2^{(\kappa, \lambda -\mu)} c_{26} (2^{-\kappa})^{-\mu -p^{-1} \e +q^{-1} \e}
\biggl\| U_{\Kappa^0 +\kappa,\nu}^{d,l -\e,m,D,U,\Nu} f
\biggr\|_{L_p(Q_{\Kappa^0 +\kappa, n_{\Kappa^0 +\kappa}^{d,m,D,U}(\nu)}^d)} = \\
c_{27} 2^{(\kappa, \lambda +p^{-1} \e -q^{-1} \e)}
\biggl\| U_{\Kappa^0 +\kappa,\nu}^{d,l -\e,m,D,U,\Nu} f
\biggr\|_{L_p(Q_{\Kappa^0 +\kappa, n_{\Kappa^0 +\kappa}^{d,m,D,U}(\nu)}^d)}.
\end{multline*}

В условиях предложения преобразуем выражение (2.1.31) (с $ l -\e $ вместо
$ l $ ) к виду, подходящему для получения интересующих нас оценок.
Для этого отметим, что в силу равенства (1.2.12) при $ d \in \N, l \in \Z_+^d,
\Delta \in \R_+^d, X^0 \in \R^d $ для области $ D \subset \R^d,
\kappa \in \Z_+^d, \nu \in \Z^d: Q_{\kappa,\nu}^d \subset D, 1 \le p < \infty $
и $ f \in L_p(D) $ почти для всех $ x \in \R^d $ имеет место равенство

\begin{equation*}
\chi_{X^0 +\Delta I^d}(x) (S_{\kappa, \nu}^{d,l} (f \mid_{Q_{\kappa,\nu}^d}))(x)
= ((\prod_{j=1}^d V_j(M_{\chi_{X_j^0 +\Delta_j I}}
S_{\kappa_j, \nu_j}^{1,l_j})) F)(x), \text{ где } F = \mathcal I^D f.
\end{equation*}

Далее, в условиях предложения для $ \kappa \in \Z_+^d, \nu \in N_{\Kappa^0 +\kappa}^{d,m,U} $
выберем $ X_{\Kappa^0 +\kappa,\nu}^{d,m,D,U} \in \R^d $ и
$ \Delta_{\Kappa^0 +\kappa,\nu}^{d,m,D,U} \in \R_+^d $ так, чтобы
при $ \epsilon \in \Upsilon^d: \s(\epsilon) \subset \s(\kappa),
\m^{\epsilon} \in \M_{\epsilon}^m(\nu) $ выполнялось включение
\begin{equation*} \tag{2.1.51}
\mathcal D_{\Kappa^0 +\kappa,\nu,\epsilon,\m^{\epsilon}}^{d,m,D,U} \subset
X_{\Kappa^0 +\kappa,\nu}^{d,m,D,U} +\Delta_{\Kappa^0 +\kappa,\nu}^{d,m,D,U} I^d.
\end{equation*}
Тогда в силу сказанного, а также благодаря п. 2 леммы 1.2.1 и (1.2.2),
в условиях предложения для $ f \in L_p(D) $ при $ \kappa \in \Z_+^d
\setminus \{0\}, \nu \in N_{\Kappa^0 +\kappa}^{d,m,U}, \epsilon \in \Upsilon^d:
\s(\epsilon) \subset \s(\kappa), \m^{\epsilon} \in \M_{\epsilon}^m(\nu),
X^0 = X_{\Kappa^0 +\kappa,\nu}^{d,m,D,U},
\Delta = \Delta_{\Kappa^0 +\kappa,\nu}^{d,m,D,U} $ почти для всех
$ x \in \R^d $ выполняется равенство
\begin{multline*} \tag{2.1.52}
\chi_{X^0 +\Delta I^d}(x) (S_{\Kappa^0 +\kappa -\epsilon,
\nu_{\Kappa^0 +\kappa -\epsilon}^{d,m,D,U}
(\n_{\epsilon}(\nu,\m^{\epsilon}))}^{d, l -\e} f)(x) = \\
\biggl((\prod_{j=1}^d V_j(M_{\chi_{X_j^0 +\Delta_j I}}
S_{\Kappa^0_j +\kappa_j -\epsilon_j,
(\nu_{\Kappa^0 +\kappa -\epsilon}^{d,m,D,U}
(\n_{\epsilon}(\nu,\m^{\epsilon})))_j}^{1,l_j -1})) F\biggr)(x) = \\
\biggl((\prod_{j=1}^d (E -V_j(E -M_{\chi_{X_j^0 +\Delta_j I}}
S_{\Kappa^0_j +\kappa_j -\epsilon_j,
(\nu_{\Kappa^0 +\kappa -\epsilon}^{d,m,D,U}
(\n_{\epsilon}(\nu,\m^{\epsilon})))_j}^{1,l_j -1}))) F\biggr)(x) = \\
\sum_{ \gamma \in \Upsilon^d} (-\e)^\gamma \biggl((\prod_{j \in \s(\gamma)}
V_j(E -M_{\chi_{X_j^0 +\Delta_j I}} S_{\Kappa^0_j +\kappa_j -\epsilon_j,
(\nu_{\Kappa^0 +\kappa -\epsilon}^{d,m,D,U}
(\n_{\epsilon}(\nu,\m^{\epsilon})))_j}^{1,l_j -1})) F\biggr)(x).
\end{multline*}

Исходя из (2.1.31) (с $ l -\e $ вместо $ l $ ), пользуясь (2.1.52), находим, что
в условиях предложения при $ \kappa \in \Z_+^d \setminus \{0\}, \nu \in
N_{\Kappa^0 +\kappa}^{d,m,U}, X^0 = X_{\Kappa^0 +\kappa,\nu}^{d,m,D,U},
\Delta = \Delta_{\Kappa^0 +\kappa,\nu}^{d,m,D,U} $ для $ f \in L_p(D) $
имеет место равенство
\begin{multline*} \tag{2.1.53}
(U_{\Kappa^0 +\kappa,\nu}^{d,l -\e,m,D,U,\Nu} f) \mid_{X^0 +\Delta I^d} =\\
\biggl(\sum_{\epsilon \in \Upsilon^d: \s(\epsilon) \subset \s(\kappa)} (-\e)^\epsilon
(\sum_{\m^{\epsilon} \in \M_{\epsilon}^m(\nu)} A_{\m^{\epsilon}}^m
S_{\Kappa^0 +\kappa -\epsilon, \nu_{\Kappa^0 +\kappa -\epsilon}^{d,m,D,U}
(\n_{\epsilon}(\nu,\m^{\epsilon}))}^{d,l -\e} f) \biggr) \mid_{X^0 +\Delta I^d} = \\
\sum_{\epsilon \in \Upsilon^d: \s(\epsilon) \subset \s(\kappa)} (-\e)^\epsilon
(\sum_{\m^{\epsilon} \in \M_{\epsilon}^m(\nu)} A_{\m^{\epsilon}}^m
(\chi_{X^0 +\Delta I^d} (S_{\Kappa^0 +\kappa -\epsilon, \nu_{\Kappa^0 +\kappa -\epsilon}^{d,m,D,U}
(\n_{\epsilon}(\nu,\m^{\epsilon}))}^{d,l -\e} f)) \mid_{X^0 +\Delta I^d}) = \\
\sum_{\epsilon \in \Upsilon^d: \s(\epsilon) \subset \s(\kappa)} (-\e)^\epsilon
(\sum_{\m^{\epsilon} \in \M_{\epsilon}^m(\nu)} A_{\m^{\epsilon}}^m \\
\biggl(\sum_{ \gamma \in \Upsilon^d} (-\e)^\gamma ((\prod_{j \in \s(\gamma)}
V_j(E -M_{\chi_{X_j^0 +\Delta_j I}} S_{\Kappa^0_j +\kappa_j -\epsilon_j,
(\nu_{\Kappa^0 +\kappa -\epsilon}^{d,m,D,U}
(\n_{\epsilon}(\nu,\m^{\epsilon})))_j}^{1,l_j -1})) F)\biggr) \mid_{X^0 +\Delta I^d}) = \\
\biggl(\sum_{\epsilon \in \Upsilon^d: \s(\epsilon) \subset \s(\kappa)} (-\e)^\epsilon
(\sum_{\m^{\epsilon} \in \M_{\epsilon}^m(\nu)} A_{\m^{\epsilon}}^m
(\sum_{ \gamma \in \Upsilon^d} (-\e)^\gamma \\
\times \biggl((\prod_{j \in \s(\gamma)}
V_j(E -M_{\chi_{X_j^0 +\Delta_j I}} S_{\Kappa^0_j +\kappa_j -\epsilon_j,
(\nu_{\Kappa^0 +\kappa -\epsilon}^{d,m,D,U}
(\n_{\epsilon}(\nu,\m^{\epsilon})))_j}^{1,l_j -1})) F \biggr) )) \biggr) \mid_{X^0 +\Delta I^d} = \\
\biggl(\sum_{ \gamma \in \Upsilon^d} (-\e)^\gamma
(\sum_{\epsilon \in \Upsilon^d: \s(\epsilon) \subset \s(\kappa)} (-\e)^\epsilon
(\sum_{\m^{\epsilon} \in \M_{\epsilon}^m(\nu)} A_{\m^{\epsilon}}^m \\
\biggl((\prod_{j \in \s(\gamma)}
V_j(E -M_{\chi_{X_j^0 +\Delta_j I}} S_{\Kappa^0_j +\kappa_j -\epsilon_j,
(\nu_{\Kappa^0 +\kappa -\epsilon}^{d, m,D,U}
(\n_{\epsilon}(\nu,\m^{\epsilon})))_j}^{1,l_j -1})) F\biggr) ))\biggr) \mid_{X^0 +\Delta I^d} = \\
\biggl(\sum_{ \gamma \in \Upsilon^d: \s(\kappa) \setminus \s(\gamma) = \emptyset} (-\e)^\gamma
(\sum_{\epsilon \in \Upsilon^d: \s(\epsilon) \subset \s(\kappa)} (-\e)^\epsilon
(\sum_{\m^{\epsilon} \in \M_{\epsilon}^m(\nu)} A_{\m^{\epsilon}}^m \\
\biggl((\prod_{j \in \s(\gamma)}
V_j(E -M_{\chi_{X_j^0 +\Delta_j I}} S_{\Kappa^0_j +\kappa_j -\epsilon_j,
(\nu_{\Kappa^0 +\kappa -\epsilon}^{d, m,D,U}
(\n_{\epsilon}(\nu,\m^{\epsilon})))_j}^{1,l_j -1})) F \biggr) )) + \\
\sum_{ \gamma \in \Upsilon^d: \s(\kappa) \setminus \s(\gamma) \ne \emptyset} (-\e)^\gamma
(\sum_{\epsilon \in \Upsilon^d: \s(\epsilon) \subset \s(\kappa)} (-\e)^\epsilon
(\sum_{\m^{\epsilon} \in \M_{\epsilon}^m(\nu)} A_{\m^{\epsilon}}^m \\
\biggl((\prod_{j \in \s(\gamma)}
V_j(E -M_{\chi_{X_j^0 +\Delta_j I}} S_{\Kappa^0_j +\kappa_j -\epsilon_j,
(\nu_{\Kappa^0 +\kappa -\epsilon}^{d, m,D,U}
(\n_{\epsilon}(\nu,\m^{\epsilon})))_j}^{1,l_j -1})) F\biggr) )) \biggr) \mid_{X^0 +\Delta I^d}.
\end{multline*}

При $ \kappa \in \Z_+^d \setminus \{0\}, \nu \in N_{\Kappa^0 +\kappa}^{d,m,U},
X^0 = X_{\Kappa^0 +\kappa,\nu}^{d,m,D,U},
\Delta = \Delta_{\Kappa^0 +\kappa,\nu}^{d,m,D,U} $
для $ \gamma \in \Upsilon^d: \s(\kappa) \setminus \s(\gamma) \ne \emptyset, $
ввиду (1.3.12), (1.3.21) соблюдается равенство
\begin{multline*} \tag{2.1.54}
\sum_{\epsilon \in \Upsilon^d: \s(\epsilon) \subset \s(\kappa)} (-\e)^\epsilon
\biggl(\sum_{\m^{\epsilon} \in \M_{\epsilon}^m(\nu)} A_{\m^{\epsilon}}^m \\
\biggl(\prod_{j \in \s(\gamma)}
V_j(E -M_{\chi_{X_j^0 +\Delta_j I}} S_{\Kappa^0_j +\kappa_j -\epsilon_j,
(\nu_{\Kappa^0 +\kappa -\epsilon}^{d, m,D,U}
(\n_{\epsilon}(\nu,\m^{\epsilon})))_j}^{1,l_j -1})\biggr) \biggr) = \\
\sum_{\epsilon, \epsilon^\prime \in \Upsilon^d: \s(\epsilon) \subset
(\s(\kappa) \cap \s(\gamma)), \s(\epsilon^\prime) \subset
(\s(\kappa) \setminus \s(\gamma))} (-\e)^{\epsilon +\epsilon^\prime}
\biggl(\sum_{\m^{\epsilon +\epsilon^\prime} \in \M_{\epsilon +\epsilon^\prime}^m(\nu)} A_{\m^{\epsilon +\epsilon^\prime}}^m \\
\biggl(\prod_{j \in \s(\gamma)}
V_j(E -M_{\chi_{X_j^0 +\Delta_j I}} S_{\Kappa^0_j +\kappa_j -\epsilon_j -\epsilon_j^\prime,
(\nu_{\Kappa^0 +\kappa -\epsilon -\epsilon^\prime}^{d, m,D,U}
(\n_{\epsilon +\epsilon^\prime}(\nu,\m^{\epsilon +\epsilon^\prime})))_j}^{1,l_j -1})\biggr) \biggr) = \\
\sum_{\substack{\epsilon \in \Upsilon^d:\\ \s(\epsilon) \subset
(\s(\kappa) \cap \s(\gamma))}} (-\e)^\epsilon
\biggl(\sum_{\substack{\epsilon^\prime \in \Upsilon^d:\\
\s(\epsilon^\prime) \subset (\s(\kappa) \setminus \s(\gamma))}}
(-\e)^{\epsilon^\prime} \biggl(\sum_{\substack{ \m^{\epsilon +\epsilon^\prime} =
(\m^{\epsilon}, \m^{\epsilon^\prime}):\\ \m^{\epsilon} \in \M_{\epsilon}^m(\nu),
\m^{\epsilon^\prime} \in \M_{\epsilon^\prime}^m(\nu)}} (\prod_{j \in
\s(\epsilon +\epsilon^\prime)} a_{\m_j}^{m_j}) \\ \times
\biggl(\prod_{j \in \s(\gamma)}
V_j(E -M_{\chi_{X_j^0 +\Delta_j I}} S_{\Kappa^0_j +\kappa_j -\epsilon_j -\epsilon_j^\prime,
(\nu_{\Kappa^0 +\kappa -\epsilon -\epsilon^\prime}^{d, m,D,U}
(\n_{\epsilon +\epsilon^\prime}(\nu,\m^{\epsilon +\epsilon^\prime})))_j}^{1,l_j -1})\biggr) \biggr)\biggr) = \\
\sum_{\substack{\epsilon \in \Upsilon^d:\\ \s(\epsilon) \subset
(\s(\kappa) \cap \s(\gamma))}} (-\e)^\epsilon
\biggl(\sum_{\substack{\epsilon^\prime \in \Upsilon^d:\\
\s(\epsilon^\prime) \subset (\s(\kappa) \setminus \s(\gamma))}}
(-\e)^{\epsilon^\prime} \biggl(\sum_{\substack{\m^{\epsilon} \in \M_{\epsilon}^m(\nu),
\m^{\epsilon^\prime} \in \M_{\epsilon^\prime}^m(\nu),\\ \m^{\epsilon +\epsilon^\prime} =
(\m^{\epsilon}, \m^{\epsilon^\prime})}} (\prod_{j \in \s(\epsilon)} a_{\m_j}^{m_j})
(\prod_{j \in \s(\epsilon^\prime)} a_{\m_j}^{m_j}) \\ \times
\biggl(\prod_{j \in \s(\gamma)}
V_j(E -M_{\chi_{X_j^0 +\Delta_j I}} S_{\Kappa^0_j +\kappa_j -\epsilon_j -\epsilon_j^\prime,
(\nu_{\Kappa^0 +\kappa -\epsilon -\epsilon^\prime}^{d, m,D,U}
(\n_{\epsilon +\epsilon^\prime}(\nu,\m^{\epsilon +\epsilon^\prime})))_j}^{1,l_j -1}) \biggr) \biggr)\biggr).
\end{multline*}

Далее, заметим, что при $ \kappa \in \Z_+^d \setminus \{0\} $ для $ \nu \in
N_{\Kappa^0 +\kappa}^{d,m,U}, \epsilon, \epsilon^\prime \in \Upsilon^d:
\s(\epsilon) \cap \s(\epsilon^\prime) = \emptyset,
\s(\epsilon +\epsilon^\prime) \subset \s(\kappa), $ и любых
$ \m^{\epsilon +\epsilon^\prime} \in \M_{\epsilon +\epsilon^\prime}^m(\nu),
\m^{\epsilon} \in \M_{\epsilon}^m(\nu),
\m^{\epsilon^\prime} \in \M_{\epsilon^\prime}^m(\nu):
(\m^{\epsilon +\epsilon^\prime})_j = (\m^{\epsilon})_j, j \in \s(\epsilon),
(\m^{\epsilon +\epsilon^\prime})_j = (\m^{\epsilon^\prime})_j, j \in
\s(\epsilon^\prime), $ согласно (1.3.14), (1.3.15) и (2.1.34)
при $ j \in \Nu_{1,d}^1 \setminus \s(\epsilon^\prime) $ имеет место равенство
\begin{equation*}
(\nu_{\Kappa^0 +\kappa -\epsilon -\epsilon^\prime}^{d,m,D,U}
(\n_{\epsilon +\epsilon^\prime}(\nu,\m^{\epsilon +\epsilon^\prime})))_j =
(\nu_{\Kappa^0 +\kappa -\epsilon -\epsilon^\prime}^{d,m,D,U}
(\n_{\epsilon^\prime}(\n_{\epsilon}(\nu, \m^{\epsilon}), \m^{\epsilon^\prime})))_j = \\
(\nu_{\Kappa^0 +\kappa -\epsilon}^{d,m,D,U}(\n_{\epsilon}(\nu,\m^{\epsilon})))_j.
\end{equation*}
Принимая во внимание это обстоятельство,
при $ \kappa \in \Z_+^d \setminus \{0\}, \nu \in N_{\Kappa^0 +\kappa}^{d,m,U},
\gamma \in \Upsilon^d: \s(\kappa) \setminus \s(\gamma) \ne \emptyset,
\epsilon \in \Upsilon^d: \s(\epsilon) \subset (\s(\kappa) \cap \s(\gamma)),
\epsilon^\prime \in \Upsilon^d: \s(\epsilon^\prime) \subset (\s(\kappa)
\setminus \s(\gamma)), \m^{\epsilon} \in \M_{\epsilon}^m(\nu),
\m^{\epsilon^\prime} \in \M_{\epsilon^\prime}^m(\nu), \m^{\epsilon +\epsilon^\prime}: $
\begin{equation*}
(\m^{\epsilon +\epsilon^\prime})_j = \begin{cases} (\m^{\epsilon})_j,
\text{ при } j \in \s(\epsilon); \\
(\m^{\epsilon^\prime})_j, \text{ при } j \in \s(\epsilon^\prime),
\end{cases}
\end{equation*}
$ X^0 = X_{\Kappa^0 +\kappa,\nu}^{d,m,D,U},
\Delta = \Delta_{\Kappa^0 +\kappa,\nu}^{d,m,D,U} $
с учётом (1.3.21), а также включения $ \s(\gamma) \subset
\Nu_{1,d}^1 \setminus \s(\epsilon^\prime) $ выводим
\begin{multline*} \tag{2.1.55}
(\prod_{j \in \s(\epsilon)} a_{\m_j}^{m_j})
(\prod_{j \in \s(\epsilon^\prime)} a_{\m_j}^{m_j}) \\ \times
(\prod_{j \in \s(\gamma)}
V_j(E -M_{\chi_{X_j^0 +\Delta_j I}} S_{\Kappa^0_j +\kappa_j -\epsilon_j -\epsilon_j^\prime,
(\nu_{\Kappa^0 +\kappa -\epsilon -\epsilon^\prime}^{d,m,D,U}
(\n_{\epsilon +\epsilon^\prime}(\nu,\m^{\epsilon +\epsilon^\prime})))_j}^{1,l_j -1})) = \\
A_{\m^{\epsilon}}^m A_{\m^{\epsilon^\prime}}^m
(\prod_{j \in \s(\gamma)}
V_j(E -M_{\chi_{X_j^0 +\Delta_j I}} S_{\Kappa^0_j +\kappa_j -\epsilon_j,
(\nu_{\Kappa^0 +\kappa -\epsilon}^{d,m,D,U}
(\n_{\epsilon}(\nu,\m^{\epsilon})))_j}^{1,l_j -1})).
\end{multline*}

Кроме того, используя (2.1.55), (1.3.22), при $ \kappa \in \Z_+^d
\setminus \{0\}, \nu \in N_{\Kappa^0 +\kappa}^{d,m,U}, \gamma \in \Upsilon^d:
\s(\kappa) \setminus \s(\gamma) \ne \emptyset, \epsilon \in \Upsilon^d:
\s(\epsilon) \subset (\s(\kappa) \cap \s(\gamma)), \epsilon^\prime \in \Upsilon^d:
\s(\epsilon^\prime) \subset (\s(\kappa) \setminus \s(\gamma)),
X^0 = X_{\Kappa^0 +\kappa,\nu}^{d,m,D,U},
\Delta = \Delta_{\Kappa^0 +\kappa,\nu}^{d,m,D,U} $ находим, что
\begin{multline*} \tag{2.1.56}
\sum_{ \m^{\epsilon} \in \M_{\epsilon}^m(\nu), \m^{\epsilon^\prime}
\in \M_{\epsilon^\prime}^m(\nu), \m^{\epsilon +\epsilon^\prime} =
(\m^{\epsilon}, \m^{\epsilon^\prime})}
(\prod_{j \in \s(\epsilon)} a_{\m_j}^{m_j})
(\prod_{j \in \s(\epsilon^\prime)} a_{\m_j}^{m_j}) \\ \times
(\prod_{j \in \s(\gamma)}
V_j(E -M_{\chi_{X_j^0 +\Delta_j I}} S_{\Kappa^0_j +\kappa_j -\epsilon_j -\epsilon_j^\prime,
(\nu_{\Kappa^0 +\kappa -\epsilon -\epsilon^\prime}^{d, m,D,U}
(\n_{\epsilon +\epsilon^\prime}(\nu,\m^{\epsilon +\epsilon^\prime})))_j}^{1,l_j -1})) = \\
\sum_{\substack{ (\m^{\epsilon}, \m^{\epsilon^\prime}): \m^{\epsilon} \in \M_{\epsilon}^m(\nu),\\
\m^{\epsilon^\prime} \in \M_{\epsilon^\prime}^m(\nu)}}
A_{\m^{\epsilon}}^m A_{\m^{\epsilon^\prime}}^m
(\prod_{j \in \s(\gamma)}
V_j(E -M_{\chi_{X_j^0 +\Delta_j I}} S_{\Kappa^0_j +\kappa_j -\epsilon_j,
(\nu_{\Kappa^0 +\kappa -\epsilon}^{d, m,D,U}
(\n_{\epsilon}(\nu,\m^{\epsilon})))_j}^{1,l_j -1})) = \\
\sum_{\m^{\epsilon^\prime} \in \M_{\epsilon^\prime}^m(\nu)} A_{\m^{\epsilon^\prime}}^m
(\sum_{ \m^{\epsilon} \in \M_{\epsilon}^m(\nu)} A_{\m^{\epsilon}}^m
(\prod_{j \in \s(\gamma)}
V_j(E -M_{\chi_{X_j^0 +\Delta_j I}} S_{\Kappa^0_j +\kappa_j -\epsilon_j,
(\nu_{\Kappa^0 +\kappa -\epsilon}^{d, m,D,U}
(\n_{\epsilon}(\nu,\m^{\epsilon})))_j}^{1,l_j -1}))) = \\
(\sum_{\m^{\epsilon^\prime} \in \M_{\epsilon^\prime}^m(\nu)} A_{\m^{\epsilon^\prime}}^m)
(\sum_{ \m^{\epsilon} \in \M_{\epsilon}^m(\nu)} A_{\m^{\epsilon}}^m
(\prod_{j \in \s(\gamma)}
V_j(E -M_{\chi_{X_j^0 +\Delta_j I}} S_{\Kappa^0_j +\kappa_j -\epsilon_j,
(\nu_{\Kappa^0 +\kappa -\epsilon}^{d, m,D,U}
(\n_{\epsilon}(\nu,\m^{\epsilon})))_j}^{1,l_j -1}))) = \\
\sum_{ \m^{\epsilon} \in \M_{\epsilon}^m(\nu)} A_{\m^{\epsilon}}^m
(\prod_{j \in \s(\gamma)}
V_j(E -M_{\chi_{X_j^0 +\Delta_j I}} S_{\Kappa^0_j +\kappa_j -\epsilon_j,
(\nu_{\Kappa^0 +\kappa -\epsilon}^{d, m,D,U}
(\n_{\epsilon}(\nu,\m^{\epsilon})))_j}^{1,l_j -1})).
\end{multline*}

Подставляя (2.1.56) в (2.1.54), приходим к равенству
\begin{multline*} \tag{2.1.57}
\sum_{\epsilon \in \Upsilon^d: \s(\epsilon) \subset \s(\kappa)} (-\e)^\epsilon
(\sum_{\m^{\epsilon} \in \M_{\epsilon}^m(\nu)} A_{\m^{\epsilon}}^m \\
(\prod_{j \in \s(\gamma)}
V_j(E -M_{\chi_{X_j^0 +\Delta_j I}} S_{\Kappa^0_j +\kappa_j -\epsilon_j,
(\nu_{\Kappa^0 +\kappa -\epsilon}^{d, m,D,U}
(\n_{\epsilon}(\nu,\m^{\epsilon})))_j}^{1,l_j -1}))) = \\
\sum_{\substack{\epsilon \in \Upsilon^d:\\ \s(\epsilon) \subset
(\s(\kappa) \cap \s(\gamma))}} (-\e)^\epsilon
(\sum_{\substack{\epsilon^\prime \in \Upsilon^d:\\
\s(\epsilon^\prime) \subset (\s(\kappa) \setminus \s(\gamma))}}
(-\e)^{\epsilon^\prime} (\sum_{ \m^{\epsilon} \in \M_{\epsilon}^m(\nu)} A_{\m^{\epsilon}}^m \\
(\prod_{j \in \s(\gamma)}
V_j(E -M_{\chi_{X_j^0 +\Delta_j I}} S_{\Kappa^0_j +\kappa_j -\epsilon_j,
(\nu_{\Kappa^0 +\kappa -\epsilon}^{d, m,D,U}
(\n_{\epsilon}(\nu,\m^{\epsilon})))_j}^{1,l_j -1})))) = \\
\sum_{\substack{\epsilon^\prime \in \Upsilon^d:\\
\s(\epsilon^\prime) \subset (\s(\kappa) \setminus \s(\gamma))}} (-\e)^{\epsilon^\prime}
(\sum_{\substack{\epsilon \in \Upsilon^d:\\ \s(\epsilon) \subset
(\s(\kappa) \cap \s(\gamma))}} (-\e)^\epsilon
(\sum_{ \m^{\epsilon} \in \M_{\epsilon}^m(\nu)} A_{\m^{\epsilon}}^m\\
(\prod_{j \in \s(\gamma)}
V_j(E -M_{\chi_{X_j^0 +\Delta_j I}} S_{\Kappa^0_j +\kappa_j -\epsilon_j,
(\nu_{\Kappa^0 +\kappa -\epsilon}^{d, m,D,U}
(\n_{\epsilon}(\nu,\m^{\epsilon})))_j}^{1,l_j -1})))) = \\
(\sum_{\substack{\epsilon^\prime \in \Upsilon^d: \\
\s(\epsilon^\prime) \subset (\s(\kappa) \setminus \s(\gamma))}} (-\e)^{\epsilon^\prime})
(\sum_{\substack{\epsilon \in \Upsilon^d: \\ \s(\epsilon) \subset
(\s(\kappa) \cap \s(\gamma))}} (-\e)^\epsilon
(\sum_{ \m^{\epsilon} \in \M_{\epsilon}^m(\nu)} A_{\m^{\epsilon}}^m\\
(\prod_{j \in \s(\gamma)}
V_j(E -M_{\chi_{X_j^0 +\Delta_j I}} S_{\Kappa^0_j +\kappa_j -\epsilon_j,
(\nu_{\Kappa^0 +\kappa -\epsilon}^{d, m,D,U}
(\n_{\epsilon}(\nu,\m^{\epsilon})))_j}^{1,l_j -1})))) = 0,
\end{multline*}
при $ \kappa \in \Z_+^d \setminus \{0\}, \nu \in N_{\Kappa^0 +\kappa}^{d,m,U},
X^0 = X_{\Kappa^0 +\kappa,\nu}^{d,m,D,U},
\Delta = \Delta_{\Kappa^0 +\kappa,\nu}^{d,m,D,U},
\gamma \in \Upsilon^d: \s(\kappa) \setminus \s(\gamma) \ne \emptyset, $
ибо в этом случае
$$
\sum_{\epsilon^\prime \in \Upsilon^d: \s(\epsilon^\prime) \subset
(\s(\kappa) \setminus \s(\gamma))} (-\e)^{\epsilon^\prime} =0.
$$

Соединяя (2.1.57) с (2.1.53), заключаем, что в условиях предложения
при $ \kappa \in \Z_+^d \setminus \{0\}, \nu \in N_{\Kappa^0 +\kappa}^{d,m,U},
X^0 = X_{\Kappa^0 +\kappa,\nu}^{d,m,D,U},
\Delta = \Delta_{\Kappa^0 +\kappa,\nu}^{d,m,D,U}, $ для $ f \in L_p(D) $
имеет место равенство
\begin{multline*} \tag{2.1.58}
(U_{\Kappa^0 +\kappa,\nu}^{d,l -\e,m,D,U,\Nu} f) \mid_{X^0 +\Delta I^d} = \\
(\sum_{ \gamma \in \Upsilon^d: \s(\kappa) \setminus \s(\gamma) = \emptyset} (-\e)^\gamma
(\sum_{\epsilon \in \Upsilon^d: \s(\epsilon) \subset \s(\kappa)} (-\e)^\epsilon
(\sum_{\m^{\epsilon} \in \M_{\epsilon}^m(\nu)} A_{\m^{\epsilon}}^m \\
((\prod_{j \in \s(\gamma)}
V_j(E -M_{\chi_{X_j^0 +\Delta_j I}} S_{\Kappa^0_j +\kappa_j -\epsilon_j,
(\nu_{\Kappa^0 +\kappa -\epsilon}^{d, m,D,U}
(\n_{\epsilon}(\nu,\m^{\epsilon})))_j}^{1,l_j -1})) F) ))) \mid_{X^0 +\Delta I^d} = \\
(\sum_{ \gamma \in \Upsilon^d: \s(\kappa) \subset \s(\gamma)} (-\e)^\gamma
(\sum_{\epsilon \in \Upsilon^d: \s(\epsilon) \subset \s(\kappa)} (-\e)^\epsilon
(\sum_{\m^{\epsilon} \in \M_{\epsilon}^m(\nu)} A_{\m^{\epsilon}}^m \\
((\prod_{j \in \s(\gamma)}
V_j(E -M_{\chi_{X_j^0 +\Delta_j I}} S_{\Kappa^0_j +\kappa_j -\epsilon_j,
(\nu_{\Kappa^0 +\kappa -\epsilon}^{d, m,D,U}
(\n_{\epsilon}(\nu,\m^{\epsilon})))_j}^{1,l_j -1})) F) ))) \mid_{X^0 +\Delta I^d}.
\end{multline*}

Возвращаясь к оценке правой части (2.1.50), при $ \kappa \in \Z_+^d
\setminus \{0\}, \nu \in N_{\Kappa^0 +\kappa}^{d,m,U},
X^0 = X_{\Kappa^0 +\kappa,\nu}^{d,m,D,U},
\Delta = \Delta_{\Kappa^0 +\kappa,\nu}^{d,m,D,U} $ для $ f \in L_p(D), $
опираясь на (2.1.58) с учётом (2.1.36), и (2.1.51), получаем
\begin{multline*} \tag{2.1.59}
\biggl\| U_{\Kappa^0 +\kappa,\nu}^{d,l -\e,m,D,U,\Nu} f
\biggr\|_{L_p(Q_{\Kappa^0 +\kappa, n_{\Kappa^0 +\kappa}^{d,m,D,U}(\nu)}^d)} = \\
\biggl\| \sum_{ \gamma \in \Upsilon^d: \s(\kappa) \subset \s(\gamma)} (-\e)^\gamma
(\sum_{\epsilon \in \Upsilon^d: \s(\epsilon) \subset \s(\kappa)} (-\e)^\epsilon
(\sum_{\m^{\epsilon} \in \M_{\epsilon}^m(\nu)} A_{\m^{\epsilon}}^m \\
((\prod_{j \in \s(\gamma)}
V_j(E -M_{\chi_{X_j^0 +\Delta_j I}} S_{\Kappa^0_j +\kappa_j -\epsilon_j,
(\nu_{\Kappa^0 +\kappa -\epsilon}^{d, m,D,U}
(\n_{\epsilon}(\nu,\m^{\epsilon})))_j}^{1,l_j -1})) F)))
\biggr\|_{L_p(Q_{\Kappa^0 +\kappa, n_{\Kappa^0 +\kappa}^{d,m,D,U}(\nu)}^d)} \le \\
\sum_{ \gamma \in \Upsilon^d: \s(\kappa) \subset \s(\gamma)}
\sum_{\epsilon \in \Upsilon^d: \s(\epsilon) \subset \s(\kappa)}
\sum_{\m^{\epsilon} \in \M_{\epsilon}^m(\nu)} A_{\m^{\epsilon}}^m \\
\biggl\| (\prod_{j \in \s(\gamma)}
V_j(E -M_{\chi_{X_j^0 +\Delta_j I}} S_{\Kappa^0_j +\kappa_j -\epsilon_j,
(\nu_{\Kappa^0 +\kappa -\epsilon}^{d, m,D,U}
(\n_{\epsilon}(\nu,\m^{\epsilon})))_j}^{1,l_j -1})) F
\biggr\|_{L_p(Q_{\Kappa^0 +\kappa, n_{\Kappa^0 +\kappa}^{d,m,D,U}(\nu)}^d)} \le \\
\sum_{ \gamma \in \Upsilon^d: \s(\kappa) \subset \s(\gamma)}
\sum_{\epsilon \in \Upsilon^d: \s(\epsilon) \subset \s(\kappa)}
\sum_{\m^{\epsilon} \in \M_{\epsilon}^m(\nu)} A_{\m^{\epsilon}}^m \\
\biggl\| (\prod_{j \in \s(\gamma)}
V_j(E -M_{\chi_{X_j^0 +\Delta_j I}} S_{\Kappa^0_j +\kappa_j -\epsilon_j,
(\nu_{\Kappa^0 +\kappa -\epsilon}^{d, m,D,U}
(\n_{\epsilon}(\nu,\m^{\epsilon})))_j}^{1,l_j -1})) F
\biggr\|_{L_p(\mathcal D_{\Kappa^0 +\kappa,\nu,\epsilon,\m^{\epsilon}}^{d,m,D,U})}.
\end{multline*}
Оценивая правую часть (2.1.59), при $ \kappa \in \Z_+^d \setminus \{0\},
\nu \in N_{\Kappa^0 +\kappa}^{d,m,U}, \gamma \in \Upsilon^d: \s(\kappa) \subset \s(\gamma),
\epsilon \in \Upsilon^d: \s(\epsilon) \subset \s(\kappa),
\m^{\epsilon} \in \M_{\epsilon}^m(\nu),
X^0 = X_{\Kappa^0 +\kappa,\nu}^{d,m,D,U},
\Delta = \Delta_{\Kappa^0 +\kappa,\nu}^{d,m,D,U}, $ используя (1.2.2), а
затем с учётом (2.1.36), (2.1.35), (2.1.51) применяя (1.2.8), заключаем, что
\begin{multline*} \tag{2.1.60}
\biggl\| (\prod_{j \in \s(\gamma)}
V_j(E -M_{\chi_{X_j^0 +\Delta_j I}} S_{\Kappa^0_j +\kappa_j -\epsilon_j,
(\nu_{\Kappa^0 +\kappa -\epsilon}^{d, m,D,U}
(\n_{\epsilon}(\nu,\m^{\epsilon})))_j}^{1,l_j -1})) F
\biggr\|_{L_p(\mathcal D_{\Kappa^0 +\kappa,\nu,\epsilon,\m^{\epsilon}}^{d,m,D,U})} = \\
\biggl\| (\prod_{j \in \s(\gamma) \setminus \s(\kappa)}
V_j(E -M_{\chi_{X_j^0 +\Delta_j I}} S_{\Kappa^0_j +\kappa_j -\epsilon_j,
(\nu_{\Kappa^0 +\kappa -\epsilon}^{d, m,D,U}
(\n_{\epsilon}(\nu,\m^{\epsilon})))_j}^{1,l_j -1})) \\
((\prod_{j \in \s(\kappa)}
V_j(E -M_{\chi_{X_j^0 +\Delta_j I}} S_{\Kappa^0_j +\kappa_j -\epsilon_j,
(\nu_{\Kappa^0 +\kappa -\epsilon}^{d, m,D,U}
(\n_{\epsilon}(\nu,\m^{\epsilon})))_j}^{1,l_j -1})) F)
\biggr\|_{L_p(\mathcal D_{\Kappa^0 +\kappa,\nu,\epsilon,\m^{\epsilon}}^{d,m,D,U})} \le \\
c_{28} \biggl\| (\prod_{j \in \s(\kappa)}
V_j(E -M_{\chi_{X_j^0 +\Delta_j I}} S_{\Kappa^0_j +\kappa_j -\epsilon_j,
(\nu_{\Kappa^0 +\kappa -\epsilon}^{d, m,D,U}
(\n_{\epsilon}(\nu,\m^{\epsilon})))_j}^{1,l_j -1})) F
\biggr\|_{L_p(\mathcal D_{\Kappa^0 +\kappa,\nu,\epsilon,\m^{\epsilon}}^{d,m,D,U})}.
\end{multline*}

Для проведения оценки правой части (2.1.60) при $ \kappa \in \Z_+^d
\setminus \{0\}, \nu \in N_{\Kappa^0 +\kappa}^{d,m,U}, \epsilon \in \Upsilon^d:
\s(\epsilon) \subset \s(\kappa), \m^{\epsilon} \in \M_{\epsilon}^m(\nu),
j =1,\ldots,d $ обозначим через
$ \mathcal S_{\Kappa^0 +\kappa,\nu,\epsilon,\m^{\epsilon}}^{j,d,l -\e,m,D,U}:
L_1((\bm x_{\Kappa^0 +\kappa,\nu,\epsilon,\m^{\epsilon}}^{d,m,D,U})_j +
(\bm \delta_{\Kappa^0 +\kappa,\nu,\epsilon,\m^{\epsilon}}^{d,m,D,U})_j I) \mapsto
\mathcal P^{1, l_j -1}, $ линейный оператор, определяемый равенством
$$
\mathcal S_{\Kappa^0 +\kappa,\nu,\epsilon,\m^{\epsilon}}^{j,d,l -\e,m,D,U} =
P_{\delta, x^0}^{1, l_j -1}
$$
при $ \delta = (\bm \delta_{\Kappa^0 +\kappa,\nu,\epsilon,\m^{\epsilon}}^{d,m,D,U})_j,
x^0 = (\bm x_{\Kappa^0 +\kappa,\nu,\epsilon,\m^{\epsilon}}^{d,m,D,U})_j $
(см. лемму 1.1.2).

Далее, пользуясь тем, что при $ \kappa \in \Z_+^d \setminus \{0\}, \nu \in
N_{\Kappa^0 +\kappa}^{d,m,U}, \epsilon \in \Upsilon^d: \s(\epsilon) \subset
\s(\kappa), \m^{\epsilon} \in \M_{\epsilon}^m(\nu),
X^0 = X_{\Kappa^0 +\kappa,\nu}^{d,m,D,U},
\Delta = \Delta_{\Kappa^0 +\kappa,\nu}^{d,m,D,U}, j =1,\ldots,d $ ввиду
(2.1.36), (2.1.51), (1.1.4) в $ \mathcal B(L_p(\R), L_p(\R)) $ (при $ 1 \le p < \infty $)
соблюдено равенство
\begin{multline*}
(E -M_{\chi_{X_j^0 +\Delta_j I}} S_{\Kappa^0_j +\kappa_j -\epsilon_j,
(\nu_{\Kappa^0 +\kappa -\epsilon}^{d, m,D,U}
(\n_{\epsilon}(\nu,\m^{\epsilon})))_j}^{1,l_j -1})
(E -M_{\chi_{X_j^0 +\Delta_j I}}
\mathcal S_{\Kappa^0 +\kappa,\nu,\epsilon,\m^{\epsilon}}^{j,d,l -\e,m,D,U}) = \\
E -M_{\chi_{X_j^0 +\Delta_j I}} S_{\Kappa^0_j +\kappa_j -\epsilon_j,
(\nu_{\Kappa^0 +\kappa -\epsilon}^{d, m,D,U}
(\n_{\epsilon}(\nu,\m^{\epsilon})))_j}^{1,l_j -1}
-M_{\chi_{X_j^0 +\Delta_j I}}
\mathcal S_{\Kappa^0 +\kappa,\nu,\epsilon,\m^{\epsilon}}^{j,d,l -\e,m,D,U} + \\
M_{\chi_{X_j^0 +\Delta_j I}} S_{\Kappa^0_j +\kappa_j -\epsilon_j,
(\nu_{\Kappa^0 +\kappa -\epsilon}^{d, m,D,U}
(\n_{\epsilon}(\nu,\m^{\epsilon})))_j}^{1,l_j -1}
(M_{\chi_{X_j^0 +\Delta_j I}}
\mathcal S_{\Kappa^0 +\kappa,\nu,\epsilon,\m^{\epsilon}}^{j,d,l -\e,m,D,U}) = \\
E -M_{\chi_{X_j^0 +\Delta_j I}} S_{\Kappa^0_j +\kappa_j -\epsilon_j,
(\nu_{\Kappa^0 +\kappa -\epsilon}^{d, m,D,U}
(\n_{\epsilon}(\nu,\m^{\epsilon})))_j}^{1,l_j -1}
-M_{\chi_{X_j^0 +\Delta_j I}}
\mathcal S_{\Kappa^0 +\kappa,\nu,\epsilon,\m^{\epsilon}}^{j,d,l -\e,m,D,U} + \\
M_{\chi_{X_j^0 +\Delta_j I}} S_{\Kappa^0_j +\kappa_j -\epsilon_j,
(\nu_{\Kappa^0 +\kappa -\epsilon}^{d, m,D,U}
(\n_{\epsilon}(\nu,\m^{\epsilon})))_j}^{1,l_j -1}
(\mathcal S_{\Kappa^0 +\kappa,\nu,\epsilon,\m^{\epsilon}}^{j,d,l -\e,m,D,U}) = \\
E -M_{\chi_{X_j^0 +\Delta_j I}} S_{\Kappa^0_j +\kappa_j -\epsilon_j,
(\nu_{\Kappa^0 +\kappa -\epsilon}^{d, m,D,U}
(\n_{\epsilon}(\nu,\m^{\epsilon})))_j}^{1,l_j -1}
-M_{\chi_{X_j^0 +\Delta_j I}}
\mathcal S_{\Kappa^0 +\kappa,\nu,\epsilon,\m^{\epsilon}}^{j,d,l -\e,m,D,U} + \\
M_{\chi_{X_j^0 +\Delta_j I}}
\mathcal S_{\Kappa^0 +\kappa,\nu,\epsilon,\m^{\epsilon}}^{j,d,l -\e,m,D,U} =
E -M_{\chi_{X_j^0 +\Delta_j I}} S_{\Kappa^0_j +\kappa_j -\epsilon_j,
(\nu_{\Kappa^0 +\kappa -\epsilon}^{d, m,D,U}
(\n_{\epsilon}(\nu,\m^{\epsilon})))_j}^{1,l_j -1},
\end{multline*}
в силу п. 2 леммы 1.2.1 и (1.2.2), а также благодаря (2.1.36), (2.1.35),
(2.1.51) используя (1.2.8), выводим
\begin{multline*} \tag{2.1.61}
\biggl\| (\prod_{j \in \s(\kappa)}
V_j(E -M_{\chi_{X_j^0 +\Delta_j I}} S_{\Kappa^0_j +\kappa_j -\epsilon_j,
(\nu_{\Kappa^0 +\kappa -\epsilon}^{d, m,D,U}
(\n_{\epsilon}(\nu,\m^{\epsilon})))_j}^{1,l_j -1})) F
\biggr\|_{L_p(\mathcal D_{\Kappa^0 +\kappa,\nu,\epsilon,\m^{\epsilon}}^{d,m,D,U})} = \\
\biggl\| (\prod_{j \in \s(\kappa)}
V_j((E -M_{\chi_{X_j^0 +\Delta_j I}} S_{\Kappa^0_j +\kappa_j -\epsilon_j,
(\nu_{\Kappa^0 +\kappa -\epsilon}^{d, m,D,U}
(\n_{\epsilon}(\nu,\m^{\epsilon})))_j}^{1,l_j -1}) \\
(E -M_{\chi_{X_j^0 +\Delta_j I}}
\mathcal S_{\Kappa^0 +\kappa,\nu,\epsilon,\m^{\epsilon}}^{j,d,l -\e,m,D,U}))) F
\biggr\|_{L_p(\mathcal D_{\Kappa^0 +\kappa,\nu,\epsilon,\m^{\epsilon}}^{d,m,D,U})} = \\
\biggl\| (\prod_{j \in \s(\kappa)}
V_j(E -M_{\chi_{X_j^0 +\Delta_j I}} S_{\Kappa^0_j +\kappa_j -\epsilon_j,
(\nu_{\Kappa^0 +\kappa -\epsilon}^{d, m,D,U}
(\n_{\epsilon}(\nu,\m^{\epsilon})))_j}^{1,l_j -1})) \\
((\prod_{j \in \s(\kappa)}
V_j(E -M_{\chi_{X_j^0 +\Delta_j I}}
\mathcal S_{\Kappa^0 +\kappa,\nu,\epsilon,\m^{\epsilon}}^{j,d,l -\e,m,D,U})) F)
\biggr\|_{L_p(\mathcal D_{\Kappa^0 +\kappa,\nu,\epsilon,\m^{\epsilon}}^{d,m,D,U})} \le \\
c_{29} \biggl\| (\prod_{j \in \s(\kappa)} V_j(E -M_{\chi_{X_j^0 +\Delta_j I}}
\mathcal S_{\Kappa^0 +\kappa,\nu,\epsilon,\m^{\epsilon}}^{j,d,l -\e,m,D,U})) F
\biggr\|_{L_p(\mathcal D_{\Kappa^0 +\kappa,\nu,\epsilon,\m^{\epsilon}}^{d,m,D,U})}.
\end{multline*}

Подставляя (2.1.61) в (2.1.60) и применяя с учётом (2.1.51) неравенство
(1.2.3) (при $ J = \emptyset, J^\prime = \s(\kappa)$), а также принимая во
внимание (2.1.35), (2.1.33), находим, что при $ \kappa \in \Z_+^d \setminus \{0\},
\nu \in N_{\Kappa^0 +\kappa}^{d,m,U}, \gamma \in \Upsilon^d: \s(\kappa) \subset
\s(\gamma), \epsilon \in \Upsilon^d: \s(\epsilon) \subset \s(\kappa),
\m^{\epsilon} \in \M_{\epsilon}^m(\nu),
X^0 = X_{\Kappa^0 +\kappa,\nu}^{d,m,D,U},
\Delta = \Delta_{\Kappa^0 +\kappa,\nu}^{d,m,D,U} $
выполняется неравенство
\begin{multline*}
\biggl\| (\prod_{j \in \s(\gamma)}
V_j(E -M_{\chi_{X_j^0 +\Delta_j I}} S_{\Kappa^0_j +\kappa_j -\epsilon_j,
(\nu_{\Kappa^0 +\kappa -\epsilon}^{d, m,D,U}
(\n_{\epsilon}(\nu,\m^{\epsilon})))_j}^{1,l_j -1})) F
\biggr\|_{L_p(\mathcal D_{\Kappa^0 +\kappa,\nu,\epsilon,\m^{\epsilon}}^{d,m,D,U})} \le  \\
c_{30} \biggl\| (\prod_{j \in \s(\kappa)} V_j(E -M_{\chi_{X_j^0 +\Delta_j I}}
\mathcal S_{\Kappa^0 +\kappa,\nu,\epsilon,\m^{\epsilon}}^{j,d,l -\e,m,D,U})) F
\biggr\|_{L_p(\mathcal D_{\Kappa^0 +\kappa,\nu,\epsilon,\m^{\epsilon}}^{d,m,D,U})} \le \\
c_{31} (\prod_{j \in \s(\kappa)}
(\bm \delta_{\Kappa^0 +\kappa,\nu,\epsilon,\m^{\epsilon}}^{d,m,D,U})_j^{-1 /p})
\biggl(\int\limits_{ (\bm \delta_{\Kappa^0 +\kappa,\nu,\epsilon,\m^{\epsilon}}^{d,m,D,U} B^d)^{\s(\kappa)}}
\int\limits_{ (\mathcal D_{\Kappa^0 +\kappa,\nu,\epsilon,\m^{\epsilon}}^{d,m,D,U})_\xi^{l \chi_{\s(\kappa)}}}
|\Delta_\xi^{l \chi_{\s(\kappa)}} F(x)|^p dx d\xi^{\s(\kappa)}\biggr)^{1/p} \le \\
c_{32} (\prod_{j \in \s(\kappa)} 2^{\kappa_j /p})
\biggl(\int_{ (c_{23} 2^{-\kappa} B^d)^{\s(\kappa)}}
\int_{ (\mathcal D_{\Kappa^0 +\kappa,\nu,\epsilon,\m^{\epsilon}}^{d,m,D,U})_\xi^{l \chi_{\s(\kappa)}}}
|\Delta_\xi^{l \chi_{\s(\kappa)}} F(x)|^p dx d\xi^{\s(\kappa)}\biggr)^{1/p} = \\
c_{32} (\prod_{j \in \s(\kappa)} 2^{\kappa_j /p})
\biggl(\int_{ (c_{23} 2^{-\kappa} B^d)^{\s(\kappa)}}
\int_{ (\mathcal D_{\Kappa^0 +\kappa,\nu,\epsilon,\m^{\epsilon}}^{d,m,D,U})_\xi^{l \chi_{\s(\kappa)}}}
|\Delta_\xi^{l \chi_{\s(\kappa)}} f(x)|^p dx d\xi^{\s(\kappa)}\biggr)^{1/p},
\end{multline*}
которое в соединении с (2.1.59) и (1.3.22) приводит к оценке
\begin{multline*}
\biggl\| U_{\Kappa^0 +\kappa,\nu}^{d,l -\e,m,D,U,\Nu} f
\biggr\|_{L_p(Q_{\Kappa^0 +\kappa, n_{\Kappa^0 +\kappa}^{d,m,D,U}(\nu)}^d)} \le
\sum_{ \gamma \in \Upsilon^d: \s(\kappa) \subset \s(\gamma)}
\sum_{\epsilon \in \Upsilon^d: \s(\epsilon) \subset \s(\kappa)}
\sum_{\m^{\epsilon} \in \M_{\epsilon}^m(\nu)} A_{\m^{\epsilon}}^m \\
c_{32} (\prod_{j \in \s(\kappa)} 2^{\kappa_j /p})
\biggl(\int_{ (c_{23} 2^{-\kappa} B^d)^{\s(\kappa)}}
\int_{ (\mathcal D_{\Kappa^0 +\kappa,\nu,\epsilon,\m^{\epsilon}}^{d,m,D,U})_\xi^{l \chi_{\s(\kappa)}}}
|\Delta_\xi^{l \chi_{\s(\kappa)}} f(x)|^p dx d\xi^{\s(\kappa)}\biggr)^{1/p} \le \\
c_{33} (\prod_{j \in \s(\kappa)} 2^{\kappa_j /p})
\sum_{\epsilon \in \Upsilon^d: \s(\epsilon) \subset \s(\kappa)}
\sum_{\m^{\epsilon} \in \M_{\epsilon}^m(\nu)}
\biggl(\int\limits_{ (c_{23} 2^{-\kappa} B^d)^{\s(\kappa)}}
\int\limits_{ (\mathcal D_{\Kappa^0 +\kappa,\nu,\epsilon,\m^{\epsilon}}^{d,m,D,U})_\xi^{l \chi_{\s(\kappa)}}}
|\Delta_\xi^{l \chi_{\s(\kappa)}} f(x)|^p dx d\xi^{\s(\kappa)}\biggr)^{1/p}, \\
\kappa \in \Z_+^d \setminus \{0\}, \nu \in N_{\Kappa^0 +\kappa}^{d,m,U}, f \in L_p(D).
\end{multline*}
Подставляя эту оценку в (2.1.50) и используя (2.1.47) (см. также (1.3.12)),
при $ n \in \Z^d: Q_{\Kappa^0 +\kappa,n}^d \cap
G_{\Kappa^0 +\kappa}^{d,m,U} \ne \emptyset,
\nu \in N_{\Kappa^0 +\kappa}^{d,m,U}:  Q_{\Kappa^0 +\kappa,n}^d \cap
\supp g_{\Kappa^0 +\kappa, \nu}^{d,m} \ne \emptyset, \mu \in \Z_+^d(\lambda), $
имеем
\begin{multline*} \tag{2.1.62}
\biggl\| \D^\mu (U_{\Kappa^0 +\kappa,\nu}^{d,l -\e,m,D,U,\Nu} f)
\D^{\lambda -\mu} g_{\Kappa^0 +\kappa, \nu}^{d,m} \biggr\|_{L_q(Q_{\Kappa^0 +\kappa,n}^d)} \le \\
c_{27} 2^{(\kappa, \lambda +p^{-1} \e -q^{-1} \e)}
c_{33} (\prod_{j \in \s(\kappa)} 2^{\kappa_j /p}) \\
\sum_{\epsilon \in \Upsilon^d: \s(\epsilon) \subset \s(\kappa)}
\sum_{\m^{\epsilon} \in \M_{\epsilon}^m(\nu)}
\biggl(\int_{ (c_{23} 2^{-\kappa} B^d)^{\s(\kappa)}}
\int_{ (\mathcal D_{\Kappa^0 +\kappa,\nu,\epsilon,\m^{\epsilon}}^{d,m,D,U})_\xi^{l \chi_{\s(\kappa)}}}
|\Delta_\xi^{l \chi_{\s(\kappa)}} f(x)|^p dx d\xi^{\s(\kappa)}\biggr)^{1/p} = \\
c_{34} 2^{(\kappa, \lambda +p^{-1} \e -q^{-1} \e)}
(\prod_{j \in \s(\kappa)} 2^{\kappa_j /p}) \\
\sum_{\epsilon \in \Upsilon^d: \s(\epsilon) \subset \s(\kappa)}
\sum_{\m^{\epsilon} \in \M_{\epsilon}^m(\nu)}
\biggl(\int_{ (c_{23} 2^{-\kappa} B^d)^{\s(\kappa)}}
\int_{ (\mathcal D_{\Kappa^0 +\kappa,\nu,\epsilon,\m^{\epsilon}}^{d,m,D,U})_\xi^{l \chi_{\s(\kappa)}}}
|\Delta_\xi^{l \chi_{\s(\kappa)}} f(x)|^p dx d\xi^{\s(\kappa)}\biggr)^{1/p} \le \\
c_{34} 2^{(\kappa, \lambda +p^{-1} \e -q^{-1} \e)}
(\prod_{j \in \s(\kappa)} 2^{\kappa_j /p}) \\
\sum_{\epsilon \in \Upsilon^d: \s(\epsilon) \subset \s(\kappa)}
\sum_{\m^{\epsilon} \in \M_{\epsilon}^m(\nu)}
\biggl(\int_{ (c_{23} 2^{-\kappa} B^d)^{\s(\kappa)}}
\int_{ (D \cap D_{\Kappa^0 +\kappa,n}^{\prime d,m,D,U})_\xi^{l \chi_{\s(\kappa)}}}
|\Delta_\xi^{l \chi_{\s(\kappa)}} f(x)|^p dx d\xi^{\s(\kappa)}\biggr)^{1/p} \le \\
c_{34} 2^{(\kappa, \lambda +p^{-1} \e -q^{-1} \e)}
(\prod_{j \in \s(\kappa)} 2^{\kappa_j /p})
(\card \Upsilon^d) (\card \Nu_{0, m +\e}^d) \\
\biggl(\int_{ (c_{23} 2^{-\kappa} B^d)^{\s(\kappa)}}
\int_{ (D \cap D_{\Kappa^0 +\kappa,n}^{\prime d,m,D,U})_\xi^{l \chi_{\s(\kappa)}}}
|\Delta_\xi^{l \chi_{\s(\kappa)}} f(x)|^p dx d\xi^{\s(\kappa)}\biggr)^{1/p} \le \\
c_{35} 2^{(\kappa, \lambda +p^{-1} \e -q^{-1} \e)}
(\prod_{j \in \s(\kappa)} 2^{\kappa_j /p})
\biggl(\int_{ (c_{23} 2^{-\kappa} B^d)^{\s(\kappa)}}
\int_{ D_\xi^{l \chi_{\s(\kappa)}} \cap D_{\Kappa^0 +\kappa,n}^{\prime d,m,D,U}}
|\Delta_\xi^{l \chi_{\s(\kappa)}} f(x)|^p dx d\xi^{\s(\kappa)}\biggr)^{1/p}.
\end{multline*}
При $ \mu \in \Z_+^d(\lambda), r \in \N, $ объединяя (2.1.40) с (2.1.62) и
используя (2.1.16), а затем применяя (1.1.2) при $ a = p /q \le 1, $ и,
наконец, пользуясь (2.1.49), приходим к оценке
\begin{multline*} \tag{2.1.63}
\int_{\cup_{n \in \Z^d: Q_{\Kappa^0 +\kappa,n}^d \cap (r B^d) \cap
G_{\Kappa^0 +\kappa}^{d,m,U} \ne \emptyset} Q_{\Kappa^0 +\kappa,n}^d} \\
\biggl| \sum_{\nu \in N_{\Kappa^0 +\kappa}^{d,m,U}}
\D^\mu (U_{\Kappa^0 +\kappa,\nu}^{d,l -\e,m,D,U,\Nu} f)
\D^{\lambda -\mu} g_{\Kappa^0 +\kappa, \nu}^{d,m} \biggr|^q dx \le \\
\sum_{n \in \Z^d: Q_{\Kappa^0 +\kappa,n}^d \cap (r B^d) \cap
G_{\Kappa^0 +\kappa}^{d,m,U} \ne \emptyset}
\biggl(\sum_{\substack{\nu \in N_{\Kappa^0 +\kappa}^{d,m,U}: \\ Q_{\Kappa^0 +\kappa,n}^d \cap
\supp g_{\Kappa^0 +\kappa, \nu}^{d,m} \ne \emptyset}}
c_{35} 2^{(\kappa, \lambda +p^{-1} \e -q^{-1} \e)}
(\prod_{j \in \s(\kappa)} 2^{\kappa_j /p})\\
\biggl(\int_{ (c_{23} 2^{-\kappa} B^d)^{\s(\kappa)}}
\int_{ D_\xi^{l \chi_{\s(\kappa)}} \cap D_{\Kappa^0 +\kappa,n}^{\prime d,m,D,U}}
|\Delta_\xi^{l \chi_{\s(\kappa)}} f(x)|^p dx d\xi^{\s(\kappa)}\biggr)^{1/p} \biggr)^q \le \\
\sum_{n \in \Z^d: Q_{\Kappa^0 +\kappa,n}^d \cap (r B^d) \cap
G_{\Kappa^0 +\kappa}^{d,m,U} \ne \emptyset}
\biggl(c_3 c_{35} 2^{(\kappa, \lambda +p^{-1} \e -q^{-1} \e)}
(\prod_{j \in \s(\kappa)} 2^{\kappa_j /p}) \\
\biggl(\int_{ (c_{23} 2^{-\kappa} B^d)^{\s(\kappa)}}
\int_{ D_\xi^{l \chi_{\s(\kappa)}} \cap D_{\Kappa^0 +\kappa,n}^{\prime d,m,D,U}}
|\Delta_\xi^{l \chi_{\s(\kappa)}} f(x)|^p dx d\xi^{\s(\kappa)}\biggr)^{1/p} \biggr)^q = \\
(c_{36} 2^{(\kappa, \lambda +p^{-1} \e -q^{-1} \e)}
(\prod_{j \in \s(\kappa)} 2^{\kappa_j /p}) )^q \\
\sum_{\substack{n \in \Z^d: Q_{\Kappa^0 +\kappa,n}^d\\ \cap (r B^d) \cap
G_{\Kappa^0 +\kappa}^{d,m,U} \ne \emptyset}}
\biggl(\int_{ (c_{23} 2^{-\kappa} B^d)^{\s(\kappa)}}
\int_{ D_\xi^{l \chi_{\s(\kappa)}} \cap D_{\Kappa^0 +\kappa,n}^{\prime d,m,D,U}}
|\Delta_\xi^{l \chi_{\s(\kappa)}} f(x)|^p dx d\xi^{\s(\kappa)} \biggr)^{q /p} \le \\
(c_{36} 2^{(\kappa, \lambda +p^{-1} \e -q^{-1} \e)}
(\prod_{j \in \s(\kappa)} 2^{\kappa_j /p}) )^q \\
\biggl(\sum_{\substack{n \in \Z^d: Q_{\Kappa^0 +\kappa,n}^d\\ \cap (r B^d) \cap
G_{\Kappa^0 +\kappa}^{d,m,U} \ne \emptyset}}
\int_{ (c_{23} 2^{-\kappa} B^d)^{\s(\kappa)}}
\int_{ D_\xi^{l \chi_{\s(\kappa)}} \cap D_{\Kappa^0 +\kappa,n}^{\prime d,m,D,U}}
|\Delta_\xi^{l \chi_{\s(\kappa)}} f(x)|^p dx d\xi^{\s(\kappa)} \biggr)^{q /p} = \\
(c_{36} 2^{(\kappa, \lambda +p^{-1} \e -q^{-1} \e)}
(\prod_{j \in \s(\kappa)} 2^{\kappa_j /p}) )^q \\
\biggl(\int\limits_{ (c_{23} 2^{-\kappa} B^d)^{\s(\kappa)}}
\int\limits_{ D_\xi^{l \chi_{\s(\kappa)}}}
\biggl(\sum_{\substack{n \in \Z^d: Q_{\Kappa^0 +\kappa,n}^d\\ \cap (r B^d) \cap
G_{\Kappa^0 +\kappa}^{d,m,U} \ne \emptyset}}
\chi_{D_{\Kappa^0 +\kappa,n}^{\prime d,m,D,U}}(x)\biggr)
|\Delta_\xi^{l \chi_{\s(\kappa)}} f(x)|^p dx d\xi^{\s(\kappa)} \biggr)^{q /p} \le \\
(c_{36} 2^{(\kappa, \lambda +p^{-1} \e -q^{-1} \e)} (\prod_{j \in
\s(\kappa)} 2^{\kappa_j /p}) )^q \biggl(\int_{ (c_{23} 2^{-\kappa}
B^d)^{\s(\kappa)}} \int_{ D_\xi^{l \chi_{\s(\kappa)}}} c_{24}
|\Delta_\xi^{l \chi_{\s(\kappa)}} f(x)|^p dx d\xi^{\s(\kappa)} \biggr)^{q /p} = \\
\biggl( c_{37} 2^{(\kappa, \lambda +p^{-1} \e -q^{-1} \e)}
(\prod_{j \in \s(\kappa)} 2^{\kappa_j /p}) \biggl(\int_{ (c_{23}
2^{-\kappa} B^d)^{\s(\kappa)}} \int_{ D_\xi^{l \chi_{\s(\kappa)}}}
|\Delta_\xi^{l \chi_{\s(\kappa)}} f(x)|^p dx d\xi^{\s(\kappa)}
\biggr)^{1 /p} \biggr)^q.
\end{multline*}

Переходя к пределу при $ r \to \infty $ в (2.1.63), и учитывая (2.1.39),
получаем, что при $ \mu \in \Z_+^d(\lambda) $ имеет место оценка
\begin{multline*}
\biggl\| \sum_{ \nu \in N_{\Kappa^0 +\kappa}^{d,m,U}}
\D^\mu (U_{\Kappa^0 +\kappa,\nu}^{d,l -\e,m,D,U,\Nu} f)
\D^{\lambda -\mu} g_{\Kappa^0 +\kappa, \nu}^{d,m} \biggr\|_{L_q(\R^d)}^q \le \\
\biggl( c_{37} 2^{(\kappa, \lambda +p^{-1} \e -q^{-1} \e)}
(\prod_{j \in \s(\kappa)} 2^{\kappa_j /p})
\biggl(\int_{ (c_{23} 2^{-\kappa} B^d)^{\s(\kappa)}}
\int_{ D_\xi^{l \chi_{\s(\kappa)}}}
|\Delta_\xi^{l \chi_{\s(\kappa)}} f(x)|^p dx d\xi^{\s(\kappa)} \biggr)^{1 /p} \biggr)^q.
\end{multline*}

Соединяя эту оценку с (2.1.38), заключаем, что при $ \kappa \in \Z_+^d
\setminus \{0\}, 1 \le p \le q < \infty $ для $ f \in L_p(D) $ выполняется неравенство
\begin{multline*}
\| \D^\lambda \mathcal E_{\Kappa^0, \kappa}^{d,l -\e,m,D,U,\Nu} f \|_{L_q(\R^d)} \le 
\sum_{ \mu \in \Z_+^d(\lambda)} C_\lambda^\mu
c_{37} 2^{(\kappa, \lambda +p^{-1} \e -q^{-1} \e)}\times\\
(\prod_{j \in \s(\kappa)} 2^{\kappa_j /p})
\biggl(\int_{ (c_{23} 2^{-\kappa} B^d)^{\s(\kappa)}}
\int_{ D_\xi^{l \chi_{\s(\kappa)}}}
|\Delta_\xi^{l \chi_{\s(\kappa)}} f(x)|^p dx d\xi^{\s(\kappa)} \biggr)^{1 /p} \le \\
c_{22} 2^{(\kappa, \lambda +(p^{-1} -q^{-1}) \e)}
\Omega^{\prime l \chi_{\s(\kappa)}}(f, (c_{23} 2^{-\kappa})^{\s(\kappa)})_{L_p(D)},
\end{multline*}
которое совпадает с (2.1.37). Вывод (2.1.37) при $ q = \infty $ проводится
аналогично.

Если в условиях предложения множество $ U $ -- ограниченно, то при
$ 1 \le q < p $ в силу соображений, приведенных в конце доказательства
предложения 2.1.2, имеет место оценка
\begin{equation*}
\| \D^\lambda \mathcal E_{\Kappa^0,\kappa}^{d,l -\e,m,D,U,\Nu} f \|_{L_q(\R^d)}
\le c_{22} 2^{(\kappa, \lambda)}
\Omega^{\prime l \chi_{\s(\kappa)}}(f, (c_{23} 2^{-\kappa})^{\s(\kappa)})_{L_p(D)},
\kappa \in \Z_+^d \setminus \{0\}, f \in L_p(D),
\end{equation*}
совпадающая с (2.1.37) в рассматриваемой ситуации. $ \square $

Предложение 2.1.5

Пусть выполнены условия предложения 2.1.4.
Тогда если для функции $ f \in L_p(D) $ и любого непустого множества $ J \subset
\Nu_{1,d}^1 $ функция
\begin{equation*} \tag{2.1.64}
\biggl(\prod_{j \in J} t_j^{-\lambda_j -(p^{-1} -q^{-1})_+ -1}\biggr)
\Omega^{\prime l \chi_J}(f, c_{23} t^J)_{L_p(D)} \in L_1((I^d)^J), \
\end{equation*}
то в $ L_q(U) $ имеет место равенство
\begin{equation*} \tag{2.1.65}
\D^\lambda (f \mid_U) = \sum_{\kappa \in \Z_+^d}
(\D^\lambda (\mathcal E_{\Kappa^0, \kappa}^{d, l -\e,m,D,U,\Nu} f)) \mid_U.
\end{equation*}

Доказательство.

Прежде всего отметим, что в условиях теоремы (см. также (2.1.32)) ввиду
предложения 2.1.3 в $ L_p(U) $ справедливо равенство (2.1.23)
(с $ \Kappa^0 $ вместо $ \kappa^0 $). Поэтому на основании леммы 1.3.1 и
соотношений (1.3.17), (2.1.29) заключаем, что в $ L_p(U) $ имеет место
равенство
\begin{multline*} \tag{2.1.66}
f \mid_U = \sum_{\kappa \in \Z_+^d} (\sum_{\epsilon \in \Upsilon^d:
\s(\epsilon) \subset \s(\kappa)} (-\e)^\epsilon
(E_{\Kappa^0 +\kappa -\epsilon}^{d,l -\e,m,D,U,\nu_{\Kappa^0 +\kappa -\epsilon}} f) \mid_U) = \\
\sum_{\kappa \in \Z_+^d} (\sum_{\epsilon \in \Upsilon^d:
\s(\epsilon) \subset \s(\kappa)} (-\e)^\epsilon
(H_{\Kappa^0 +\kappa, \Kappa^0 +\kappa -\epsilon}^{d,l -\e,m,U}
E_{\Kappa^0 +\kappa -\epsilon}^{d,l -\e,m,D,U,\nu_{\Kappa^0 +\kappa -\epsilon}} f) \mid_U) = \\
\sum_{\kappa \in \Z_+^d} (\sum_{\epsilon \in \Upsilon^d:
\s(\epsilon) \subset \s(\kappa)} (-\e)^\epsilon
(H_{\Kappa^0 +\kappa, \Kappa^0 +\kappa -\epsilon}^{d,l -\e,m,U}
E_{\Kappa^0 +\kappa -\epsilon}^{d,l -\e,m,D,U,\nu_{\Kappa^0 +\kappa -\epsilon}} f)) \mid_U = \\
\sum_{ \kappa \in \Z_+^d} (\mathcal E_{\Kappa^0, \kappa}^{d,l -\e,m,D,U,\Nu} f) \mid_U.
\end{multline*}

Далее, для любой функции $ f \in L_p(D) $ и каждого множества $ J \subset
\Nu_{1,d}^1: J \ne \emptyset, $ при $ \kappa \in \Z_+^d: \s(\kappa) = J, $
ввиду (2.1.37) выполняется неравенство
\begin{multline*} \tag{2.1.67}
\| \D^\lambda ((\mathcal E_{\Kappa^0, \kappa}^{d, l -\e,m,D,U,\Nu} f) \mid_U) \|_{L_q(U)} = \\
\| (\D^\lambda \mathcal E_{\Kappa^0, \kappa}^{d, l -\e,m,D,U,\Nu} f) \mid_U \|_{L_q(U)} \le \\
c_{38} \biggl(\prod_{j \in J} 2^{\kappa_j (\lambda_j +p^{-1}
+(p^{-1} -q^{-1})_+)}\biggr) \biggl(\int_{ (c_{23} 2^{-\kappa} B^d)^J}
\int_{ D_\xi^{l \chi_J}} |\Delta_\xi^{l \chi_J} f(x)|^p dx d\xi^J\biggr)^{1/p} \le \\
c_{39} \int\limits_{ (2^{-\kappa} +2^{-\kappa} I^d)^J}
\biggl(\prod_{j \in J} t_j^{-\lambda_j -p^{-1} -(p^{-1} -q^{-1})_+ -1}\biggr)
\biggl(\int_{ (c_{23} t B^d)^J} \int_{ D_\xi^{l \chi_J}}
|\Delta_\xi^{l \chi_J} f(x)|^p dx d\xi^J\biggr)^{1/p} dt^J \le \\
c_{40} \int_{ (2^{-\kappa} +2^{-\kappa} I^d)^J}
\biggl(\prod_{j \in J} t_j^{-\lambda_j -(p^{-1} -q^{-1})_+ -1} \biggr)
\Omega^{\prime l \chi_J}(f, c_{23} t^J)_{L_p(D)} dt^J.
\end{multline*}

Поскольку множества $ (2^{-\kappa} +2^{-\kappa} I^d)^J, \kappa \in
\Z_+^d: \s(\kappa) = J, $ попарно не пересекаются и $ \cup_{ \kappa
\in \Z_+^d: \s(\kappa) = J} (2^{-\kappa} +2^{-\kappa} I^d)^J
\subset (I^d)^J, $ то из (2.1.67) вытекает, что для функции $ f \in
L_p(D), $ подчинённой условию (2.1.64), ряд
$$
\sum_{ \kappa \in \Z_+^d: \s(\kappa) = J} \| \D^\lambda ((\mathcal E_{\Kappa^0,
\kappa}^{d,l -\e,m,D,U,\Nu} f) \mid_U) \|_{L_q(U)}
$$
сходится для $ J \subset \Nu_{1,d}^1: J \ne \emptyset, $
а, следовательно, и ряд $ \sum_{ \kappa \in \Z_+^d}
\| \D^\lambda ((\mathcal E_{\Kappa^0,\kappa}^{d, l -\e,m,D,U,\Nu} f) \mid_U) \|_{L_q(U)} $ сходится, и,
значит, ряд
$$
\sum_{ \kappa \in \Z_+^d} \D^\lambda ((\mathcal E_{\Kappa^0,\kappa}^{d, l -\e,m,D,U,\Nu} f) \mid_U)
$$
сходится в $ L_q(U). $

Принимая во внимание это обстоятельство и равенство (2.1.66), для
любой функции $ \phi \in C_0^\infty (U) $ имеем
\begin{multline*}
\langle \D^\lambda (f \mid_U), \phi \rangle =
\int_{U} (\sum_{ \kappa \in \Z_+^d} (\mathcal E_{\Kappa^0,\kappa}^{d,l -\e,m,D,U,\Nu} f) \mid_U)
(-1)^{(\lambda,\e)} \D^\lambda \phi dx = \\
\sum_{ \kappa \in \Z_+^d}
\int_{U} (\mathcal E_{\Kappa^0,\kappa}^{d,l -\e,m,D,U,\Nu} f) \mid_U
(-1)^{(\lambda,\e)} \D^\lambda \phi dx = \\
\sum_{ \kappa \in \Z_+^d}
\int_{U} \D^\lambda ((\mathcal E_{\Kappa^0,\kappa}^{d, l -\e,m,D,U,\Nu} f) \mid_U) \phi dx = \\
\int_{U} (\sum_{ \kappa \in \Z_+^d}
\D^\lambda ((\mathcal E_{\Kappa^0,\kappa}^{d, l -\e,m,D,U,\Nu} f) \mid_U)) \phi dx
= \int_{U} (\sum_{ \kappa \in \Z_+^d}
(\D^\lambda (\mathcal E_{\Kappa^0,\kappa}^{d, l -\e,m,D,U,\Nu} f)) \mid_U) \phi dx.
\end{multline*}

А это значит, что в $ L_q(U) $ верно равенство (2.1.65). $ \square $

Предложение 2.1.6

Пусть $ d \in \N, m \in \N^d, $ область $ D \subset \R^d $ и её открытое
подмножество $ U \subset D $ являются $m$-правильной парой. Пусть также
$ \alpha \in \R_+^d, 1 \le p < \infty, p \le q \le \infty $ и если множество
$ U $ -- ограниченно, пусть ещё $ 1 \le q < p, \lambda \in \Z_+^d(m) $ таковы,
что соблюдается условие
\begin{equation*} \tag{2.1.68}
\alpha -\lambda -(p^{-1} -q^{-1})_+ \e >0.
\end{equation*}
Тогда

1) для любой функции $ f \in (S_p^\alpha H)^\prime(D) $ при $ l = l(\alpha) $
в $ L_q(U) $ имеет место равенство (2.1.65);

2) существует константа $ c_{41}(d,\alpha,p,q,\lambda,m,D,U) > 0 $ такая,
что для $ f \in (S_p^\alpha H)^\prime(D) $ выполняется неравенство
\begin{equation*} \tag{2.1.69}
\| \D^\lambda (f \mid_U) \|_{L_q(U)} \le c_{41} \| f \|_{(S_p^\alpha H)^\prime(D)}.
\end{equation*}

Доказательство.

Для доказательства п. 1) предложения достаточно заметить, что для $ f \in
(S_p^\alpha H)^\prime(D), l = l(\alpha)  $ и любого непустого множества $ J
\subset \Nu_{1,d}^1 $ при $ t^J \in (\R_+^d)^J $ имеет место неравенство
\begin{multline*}
\Omega^{\prime l \chi_J}(f, t^J)_{L_p(D)} = (t^J)^{\alpha^J}
(t^J)^{-\alpha^J},\Omega^{\prime l \chi_J}(f, t^J)_{L_p(D)} \le \\
(t^J)^{\alpha^J} \sup_{\tau^J \in (\R_+^d)^J}
(\tau^J)^{-\alpha^J},\Omega^{\prime l \chi_J}(f, \tau^J)_{L_p(D)} \le
\| f \|_{(S_p^\alpha H)^\prime(D)} (t^J)^{\alpha^J},
\end{multline*}
а, значит, ввиду (2.1.68) соблюдается (2.1.64), и применить предложение 2.1.5.

Теперь установим справедливость п. 2). При $ l = l(\alpha), $ в соответствии с
п. 1) предложения, используя (2.1.65), а также применяя (2.1.7) (см. (2.1.29),
(1.3.16)) и (2.1.37), для $ f \in (S_p^\alpha H)^\prime(D) $ имеем
\begin{multline*}
\| \D^\lambda (f \mid_U) \|_{L_q(U)} = \| \sum_{\kappa \in \Z_+^d}
(\D^\lambda (\mathcal E_{\Kappa^0, \kappa}^{d, l -\e,m,D,U,\Nu} f)) \mid_U \|_{L_q(U)} \le \\
\sum_{\kappa \in \Z_+^d} \| (\D^\lambda (\mathcal E_{\Kappa^0, \kappa}^{d, l -\e,m,D,U,\Nu} f)) \mid_U \|_{L_q(U)} \le
\sum_{\kappa \in \Z_+^d} \| \D^\lambda \mathcal E_{\Kappa^0, \kappa}^{d, l -\e,m,D,U,\Nu} f \|_{L_q(\R^d)} = \\
\| \D^\lambda \mathcal E_{\Kappa^0,0}^{d, l -\e,m,D,U,\Nu} f \|_{L_q(\R^d)} +
\sum_{\kappa \in \Z_+^d: \kappa \ne 0}
\| \D^\lambda \mathcal E_{\Kappa^0, \kappa}^{d, l -\e,m,D,U,\Nu} f \|_{L_q(\R^d)} = \\
\| \D^\lambda E_{\Kappa^0}^{d, l -\e,m,D,U,\nu_{\Kappa^0}} f \|_{L_q(\R^d)} +
\sum_{J \subset \Nu_{1,d}^1: J \ne \emptyset}
\sum_{\kappa \in \Z_+^d: \s(\kappa) = J}
\| \D^\lambda \mathcal E_{\Kappa^0, \kappa}^{d, l -\e,m,D,U,\Nu} f \|_{L_q(\R^d)} \le \\
c_1 \| f\|_{L_p(D)} +\sum_{J \subset \Nu_{1,d}^1: J \ne \emptyset}
\sum_{\kappa \in \Z_+^d: \s(\kappa) = J} c_{22} 2^{(\kappa^J, \lambda^J +(p^{-1} -q^{-1})_+ \e^J)}
\Omega^{\prime l \chi_{J}}(f, (c_{23} 2^{-\kappa})^{J})_{L_p(D)} \le \\
c_1 \| f\|_{L_p(D)} +\sum_{J \subset \Nu_{1,d}^1: J \ne \emptyset}
\sum_{\kappa \in \Z_+^d: \s(\kappa) = J} c_{22} 2^{(\kappa^J, \lambda^J +(p^{-1} -q^{-1})_+ \e^J)} \\
((c_{23} 2^{-\kappa})^{J})^{\alpha^J} \| f\|_{(S_p^\alpha H)^\prime(D)} \le \\
c_1 \| f\|_{(S_p^\alpha H)^\prime(D)} +\sum_{J \subset \Nu_{1,d}^1: J \ne \emptyset}
\sum_{\kappa \in \Z_+^d: \s(\kappa) = J} c_{42} 2^{-(\kappa^J, \alpha^J -\lambda^J -(p^{-1} -q^{-1})_+ \e^J)}
\| f\|_{(S_p^\alpha H)^\prime(D)} = \\
c_1 \| f\|_{(S_p^\alpha H)^\prime(D)} +c_{42}
\sum_{J \subset \Nu_{1,d}^1: J \ne \emptyset}
\sum_{\kappa \in \Z_+^d: \s(\kappa) = J} 2^{-(\kappa, \alpha -\lambda -(p^{-1} -q^{-1})_+ \e)}
\| f\|_{(S_p^\alpha H)^\prime(D)} = \\
(c_1 +c_{42} \sum_{\kappa \in \Z_+^d \setminus \{0\}}
2^{-(\kappa, \alpha -\lambda -(p^{-1} -q^{-1})_+ \e)})
\| f\|_{(S_p^\alpha H)^\prime(D)} = c_{41} \| f\|_{(S_p^\alpha H)^\prime(D)},
\end{multline*}
что совпадает с (2.1.69). $\square $
\bigskip

2.2. В этом пункте будут построены средства приближения функций из
рассматриваемых нами пространств, на которые опирается вывод основных
результатов работы.

Лемма 2.2.1

При $ d \in \N $ пусть область $ D \subset \R^d $ и её открытое подмножество
$ U \subset D $ таковы, что существует $ \delta \in \R_+^d, $ для которого имеет
место включение $ (U +\delta I^d) \subset D. $ Тогда $ (D,U) $ является
$m$-правильной парой для любого $ m \in \N^d. $

Доказательство.

В условиях леммы при $ m \in \N^d $ и $ \Kappa^0, \kappa \in \Z_+^d $ для $ \nu
\in N_{\Kappa^0 +\kappa}^{d,m,U}, $ учитывая открытость $ U $ и (1.3.4),
выберем $ x \in U \cap (2^{-\Kappa^0 -\kappa} \nu +2^{-\Kappa^0 -\kappa} (m +\e) I^d). $
Тогда $ x \in U $ и при $ j = 1,\ldots,d $ выполняются соотношения
\begin{equation*} \tag{2.2.1}
2^{-\Kappa^0_j -\kappa_j} \nu_j < x_j < 2^{-\Kappa^0_j -\kappa_j} \nu_j +
2^{-\Kappa^0_j -\kappa_j} (m_j +1),
\end{equation*}
а для $ y \in (2^{-\Kappa^0 -\kappa} (\nu +m +\e) +2^{-\Kappa^0 -\kappa} I^d) $
при $ j = 1,\ldots,d $ соблюдаются неравенства
$$
2^{-\Kappa^0_j -\kappa_j} (\nu_j +m_j +1) < y_j < 2^{-\Kappa^0_j -\kappa_j}
(\nu_j +m_j +1) +2^{-\Kappa^0_j -\kappa_j},
$$
и, значит,
\begin{equation*} \tag{2.2.2}
x_j < 2^{-\Kappa^0_j -\kappa_j} (\nu_j +m_j +1) < y_j <
2^{-\Kappa^0_j -\kappa_j} \nu_j +2^{-\Kappa^0_j -\kappa_j} (m_j +2) < x_j
+2^{-\Kappa^0_j -\kappa_j} (m_j +2),
\end{equation*}
т.е.
\begin{equation*} \tag{2.2.3}
Q_{\Kappa^0 +\kappa, \nu +m +\e}^d \subset x +2^{-\Kappa^0 -\kappa} (m +2 \e) I^d
\subset (U +2^{-\Kappa^0 -\kappa} (m +2 \e) I^d).
\end{equation*}

Теперь в условиях леммы при $ m \in \N^d $ фиксируем $ \Kappa^0 =
\Kappa^0(d,m,D,U) \in \Z_+^d $ такое, что $ 2^{-\Kappa^0} (m +2 \e) < \delta, $
и определим при $ \kappa \in \Z_+^d $ отображения
$ \nu_{\Kappa^0 +\kappa}^{d,m,D,U}: N_{\Kappa^0 +\kappa}^{d,m,U} \mapsto \Z^d,
n_{\Kappa^0 +\kappa}^{d,m,D,U}: N_{\Kappa^0 +\kappa}^{d,m,U} \mapsto \Z^d, $
полагая для $ \nu \in N_{\Kappa^0 +\kappa}^{d,m,U} $ значение
$$
\nu_{\Kappa^0 +\kappa}^{d,m,D,U}(\nu) = n_{\Kappa^0 +\kappa}^{d,m,D,U}(\nu) =
\nu +m +\e.
$$
Проверим, что для рассматриваемых объектов соблюдаются условия определения
$m$-правильной пары. Для этого при $ \kappa \in \Z_+^d $ для $ \nu \in
N_{\Kappa^0 +\kappa}^{d,m,U} $ в силу (2.2.3) имеем
\begin{multline*}
Q_{\Kappa^0 +\kappa,\nu_{\Kappa^0 +\kappa}^{d,m,D,U}(\nu)}^d =
Q_{\Kappa^0 +\kappa,n_{\Kappa^0 +\kappa}^{d,m,D,U}(\nu)}^d =
Q_{\Kappa^0 +\kappa, \nu +m +\e}^d \subset\\
 (U +2^{-\Kappa^0 -\kappa} (m +2 \e) I^d)
\subset (U +2^{-\Kappa^0} (m +2 \e) I^d) \subset (U +\delta I^d) \subset D
\end{multline*}
и
\begin{multline*}
Q_{\Kappa^0 +\kappa,\nu_{\Kappa^0 +\kappa}^{d,m,D,U}(\nu)}^d =
Q_{\Kappa^0 +\kappa,n_{\Kappa^0 +\kappa}^{d,m,D,U}(\nu)}^d =
Q_{\Kappa^0 +\kappa, \nu +m +\e}^d =\\
(2^{-\Kappa^0 -\kappa} (\nu +m +\e) +2^{-\Kappa^0 -\kappa} I^d) \subset
(2^{-\Kappa^0 -\kappa} \nu +2^{-\Kappa^0 -\kappa} (m +\e +I^d)) \subset\\
(2^{-\Kappa^0 -\kappa} \nu +2^{-\Kappa^0 -\kappa} (m +2 \e) \overline I^d) \subset
(2^{-\Kappa^0 -\kappa} \nu +2^{-\Kappa^0 -\kappa} (m +2 \e) B^d),
\end{multline*}
т.е. выполняется (2.1.32) с $ \Gamma^0 = m +2 \e. $

Далее, при $ \kappa \in \Z_+^d, \nu \in N_{\Kappa^0 +\kappa}^{d,m,U},
\epsilon \in \Upsilon^d: \s(\epsilon) \subset \s(\kappa), \m^{\epsilon} \in
\M_{\epsilon}^m(\nu), $ выбирая $ x \in (U \cap (2^{-\Kappa^0 -\kappa} \nu +2^{-\Kappa^0 -\kappa} (m +\e) I^d)), $
ввиду (2.2.2), (2.2.1) для $ j =1.\ldots,d $ имеем
\begin{multline*} \tag{2.2.4}
x_j < 2^{-\Kappa^0_j -\kappa_j} (\nu_{\Kappa^0 +\kappa}^{d,m,D,U}(\nu))_j <
2^{-\Kappa^0_j -\kappa_j} (\nu_{\Kappa^0 +\kappa}^{d,m,D,U}(\nu))_j
+2^{-\Kappa^0_j -\kappa_j} =\\
 2^{-\Kappa^0_j -\kappa_j} \nu_j + 2^{-\Kappa^0_j -\kappa_j} (m_j +2) <
x_j +2^{-\Kappa^0_j -\kappa_j} (m_j +2) < \\
x_j +2^{-\Kappa^0_j} (m_j +2) < x_j +\delta_j.
\end{multline*}
Принимая во внимание, что в описанных условиях имеет место включение
$$
x \in (U \cap (2^{-\Kappa^0 -\kappa +\epsilon} \n_{\epsilon}(\nu,\m^{\epsilon})
+2^{-\Kappa^0 -\kappa +\epsilon} (m +\e) I^d)), \text{ (см. вывод (1.3.15))},
$$
заключаем, что для $ j = 1,\ldots,d $ соблюдается (2.2.4) при
$ \Kappa^0 +\kappa -\epsilon $ вместо $ \Kappa^0 +\kappa $ и
$ \n_{\epsilon}(\nu,\m^{\epsilon}) \in N_{\Kappa^0 +\kappa -\epsilon}^{d,m,U} $
вместо $ \nu, $ т.е.
\begin{multline*} \tag{2.2.5}
x_j < 2^{-\Kappa^0_j -\kappa_j +\epsilon_j}
(\nu_{\Kappa^0 +\kappa -\epsilon}^{d,m,D,U}(\n_{\epsilon}(\nu,\m^{\epsilon})))_j <\\
2^{-\Kappa^0_j -\kappa_j +\epsilon_j}
(\nu_{\Kappa^0 +\kappa -\epsilon}^{d,m,D,U}(\n_{\epsilon}(\nu,\m^{\epsilon})))_j +
2^{-\Kappa^0_j -\kappa_j +\epsilon_j} =\\
2^{-\Kappa^0_j -\kappa_j +\epsilon_j} (\n_{\epsilon}(\nu,\m^{\epsilon}))_j +
2^{-\Kappa^0_j -\kappa_j +\epsilon_j} (m_j +2) <\\
x_j +2^{-\Kappa^0_j -\kappa_j +\epsilon_j} (m_j +2) \le
x_j +2^{-\Kappa^0_j} (m_j +2) < x_j +\delta_j.
\end{multline*}
Из (2.2.4) и (2.2.5) с учётом определений получаем, что при
$ j =1,\ldots,d $ справедливо соотношение
\begin{multline*}
x_j < \min(2^{-\Kappa^0_j -\kappa_j} (n_{\Kappa^0 +\kappa}^{d,m,D,U}(\nu))_j,
2^{-\Kappa^0_j -\kappa_j +\epsilon_j} (\nu_{\Kappa^0 +\kappa -\epsilon}^{d,m,D,U}
(\n_{\epsilon}(\nu,\m^{\epsilon})))_j) = \\
(\bm x_{\Kappa^0 +\kappa,\nu,\epsilon,\m^{\epsilon}}^{d,m,D,U})_j <
(\bm x_{\Kappa^0 +\kappa,\nu,\epsilon,\m^{\epsilon}}^{d,m,D,U})_j +
(\bm \delta_{\Kappa^0 +\kappa,\nu,\epsilon,\m^{\epsilon}}^{d,m,D,U})_j =\\
\max(2^{-\Kappa^0_j -\kappa_j} (n_{\Kappa^0 +\kappa}^{d,m,D,U}(\nu))_j
+2^{-\Kappa^0_j -\kappa_j}, \\
2^{-\Kappa^0_j -\kappa_j +\epsilon_j} (\nu_{\Kappa^0 +\kappa -\epsilon}^{d,m,D,U}
(\n_{\epsilon}(\nu,\m^{\epsilon})))_j +2^{-\Kappa^0_j -\kappa_j +\epsilon_j}) <
x_j +\delta_j,
\end{multline*}
а, значит,
\begin{equation*}
\mathcal D_{\Kappa^0 +\kappa,\nu,\epsilon,\m^{\epsilon}}^{d,m,D,U} =
\bm x_{\Kappa^0 +\kappa,\nu,\epsilon,\m^{\epsilon}}^{d,m,D,U} +
\bm \delta_{\Kappa^0 +\kappa,\nu,\epsilon,\m^{\epsilon}}^{d,m,D,U} I^d
\subset x +\delta I^d \subset (U +\delta I^d) \subset D,
\end{equation*}
что приводит к (2.1.33).

Наконец, при $ \kappa \in \Z_+^d, \nu \in N_{\Kappa^0 +\kappa}^{d,m,U} $ для
любых $ \epsilon \in \Upsilon^d: \s(\epsilon) \subset \s(\kappa), $ и
$ \m^{\epsilon} \in \M_{\epsilon}^m(\nu) $ при $ j \in \Nu_{1,d}^1 \setminus \s(\epsilon) $
в силу определений и (1.3.13) имеем
\begin{multline*}
(\nu_{\Kappa^0 +\kappa -\epsilon}^{d,m,D,U}
(\n_{\epsilon}(\nu,\m^{\epsilon})))_j = (\n_{\epsilon}(\nu,\m^{\epsilon}) +m +\e)_j =
(\n_{\epsilon}(\nu,\m^{\epsilon}))_j  +m_j +1 = \\
\nu_j +m_j +1 = (\nu +m +\e)_j =
(\nu_{\Kappa^0 +\kappa}^{d,m,D,U}(\nu))_j,
\end{multline*}
что совпадает с (2.1.34). $ \square $

При $ d \in \N, $ определяя
$$
\Sigma^d = \{ \sigma \in \Z^d: \sigma_j \in \{-1,1\}, j =1,\ldots,d\},
$$
для $ \sigma \in \Sigma^d $ обозначим через $ \bm h_\sigma $ отображение,
которое каждой функции $ f, $ заданной на некотором множестве $ S \subset \R^d, $
ставит в соответствие функцию $ \bm h_\sigma f, $ определяемую на множестве
$ \{ x \in \R^d: \sigma x \in S\} = \sigma^{-1} S = \sigma S $ равенством
$ (\bm h_\sigma f)(x) = f(\sigma x). $
Так как для $ \sigma \in \Sigma^d $ отображение
$ \R^d \ni x \mapsto \sigma x \in \R^d $ --- взаимно однозначно, то отображение
$ \bm h_\sigma $ является биекцией на себя множества всех функций с
областью определения в $ \R^d. $
При этом обратное отображение $ \bm h_\sigma^{-1} $ для $ f: S \mapsto \R $
задается равенством
\begin{equation*} \tag{2.2.6}
(\bm h_\sigma^{-1} f)(x) = f(\sigma^{-1} x) = f(\sigma x) =
(\bm h_\sigma f)(x), x \in \sigma S.
\end{equation*}

Отметим некоторые полезные для нас свойства отображений $ \bm h_\sigma. $
При $ d \in \N, \sigma \in \Sigma^d $ для любых множеств $ S \subset S^\prime
\subset \R^d $ и любой функции
$ f: S^\prime \mapsto \R $ верно равенство
\begin{equation*} \tag{2.2.7}
\bm h_\sigma (f \mid_S) = (\bm h_\sigma f) \mid_{\sigma^{-1} S}.
\end{equation*}
В самом деле, для $ x \in \sigma^{-1} S $ имеем
\begin{equation*}
(\bm h_\sigma (f \mid_S))(x) = (f \mid_S)(\sigma x) = f(\sigma x) =
(\bm h_\sigma f)(x) = ((\bm h_\sigma f) \mid_{\sigma^{-1} S})(x).
\end{equation*}

При $ d \in \N, \sigma \in \Sigma^d $ для открытого множества $ D \subset \R^d,
1 \le p \le \infty $ и $ f \in L_p(D) $ имеет место равенство
\begin{equation*} \tag{2.2.8}
\| \bm h_\sigma f \|_{L_p(\sigma^{-1} D)} = \| f \|_{L_p(D)},
\text{ а, значит, } \bm h_\sigma \in \mathcal B(L_p(D),L_p(\sigma^{-1} D)).
\end{equation*}

Действительно, делая замену переменных в интеграле, получаем
\begin{multline*}
\| f \|_{L_p(D)}^p = \int_D |f(y)|^p dy = \int_{x: \sigma x \in D}
| f(\sigma x)|^p | \sigma^\e| dx = \int_{\sigma^{-1} D} | f(\sigma x)|^p dx =\\
\int_{\sigma^{-1} D} | (\bm h_\sigma f)(x)|^p dx =
\| \bm h_\sigma f \|_{L_p(\sigma^{-1} D)}^p,
\end{multline*}
откуда следует (2.2.8) при $ 1 \le p < \infty. $
При $ p = \infty $ имеем
\begin{multline*}
\| \bm h_\sigma f \|_{L_\infty(\sigma^{-1} D)} =
\supvrai_{x \in \sigma^{-1} D}| (\bm h_\sigma f)(x)| = 
\supvrai_{x \in \sigma^{-1} D}| f(\sigma x)| =\\
\supvrai_{y \in D}| f(\sigma \sigma^{-1} y)| =
\supvrai_{y \in D}| f(y)| = \| f \|_{L_\infty(D)}.
\end{multline*}

Отметим ещё, что при $ d \in \N, \sigma \in \Sigma^d $ вследствие (2.2.8) и
равенства
\begin{equation*}
\sigma Q_{\kappa,\nu}^d = Q_{\kappa,\sigma(\nu +\chi_J)}^d, \\
\text{ где } J = \{ j \in \Nu_{1,d}^1: \sigma_j = -1\}, \kappa \in \Z_+^d,
\nu \in \Z^d,
\end{equation*}
имеет место равенство
\begin{equation*}
\bm h_\sigma(L_1^{\loc \square}(D)) = L_1^{\loc \square}(\sigma^{-1} D), \\
D  \text{ -- произвольное открытое множество в } \R^d.
\end{equation*}

Лемма 2.2.2

Пусть $ d \in \N, \sigma \in \Sigma^d, D $ -- открытое множество в $ \R^d,
1 \le p, q \le \infty, \lambda \in \Z_+^d $ и $ f \in L_p(D),
\D^\lambda f \in L_q(D). $ Тогда соблюдается равенство
\begin{equation*} \tag{2.2.9}
\D^\lambda (\bm h_\sigma f) = \sigma^\lambda \bm h_\sigma (\D^\lambda f).
\end{equation*}

Доказательство.

В условиях леммы для $ \phi \in C_0^\infty(\sigma^{-1} D), $ учитывая (2.2.8) и
делая замену переменных в интегралах, имеем
\begin{multline*}
\langle \D^\lambda (\bm h_\sigma f), \phi \rangle =
(-1)^{(\lambda, \e)} \langle (\bm h_\sigma f), \D^\lambda
\phi \rangle = (-1)^{(\lambda, \e)} \int_{\sigma^{-1} D} (\bm
h_\sigma f)(y) (\D^\lambda \phi)(y) dy = \\
(-1)^{(\lambda, \e)}
\int_{\sigma^{-1} D} f(\sigma y) (\D^\lambda \phi)(y) dy =\\
(-1)^{(\lambda, \e)} \int_{\sigma \sigma^{-1} D} f(\sigma
\sigma^{-1} x) (\D^\lambda \phi)(\sigma^{-1} x) |
(\sigma^{-1})^\e| dx =\\
 (-1)^{(\lambda, \e)} \int_{D} f(x)
(\frac{\D^{|\lambda|} \phi} {\D y_1^{\lambda_1} \ldots \D
y_d^{\lambda_d}}) (\sigma^{-1} x) dx =\\
 (-1)^{(\lambda, \e)} \int_D
f(x) \sigma^{\lambda} \frac{\D^{|\lambda|}} {\D x_1^{\lambda_1}
\ldots \D x_d^{\lambda_d}} (\phi(\sigma^{-1} x)) dx =\\
(-1)^{(\lambda, \e)} \int_D f(x) \sigma^{\lambda} \D^\lambda
(\phi(\sigma^{-1} x)) dx = \sigma^{\lambda} \int_D \D^\lambda f(x)
\phi(\sigma^{-1} x) dx =\\
 \sigma^{\lambda} \int_{\sigma^{-1} D}
(\D^\lambda f)(\sigma y) \phi(\sigma^{-1} \sigma y) | \sigma^\e|
dy = \sigma^{\lambda} \int_{\sigma^{-1} D} (\bm
h_\sigma(\D^\lambda f))(y) \phi(y) dy,
\end{multline*}
что влечёт (2.2.9). $ \square $

Лемма 2.2.3

Пусть $ d \in \N, D $ -- открытое множество в $ \R^d, 1 \le p < \infty,
\sigma \in \Sigma^d. $ Тогда для
$ f \in L_p(D), l \in \N^d, J \subset \Nu_{1,d}^1: J \ne \emptyset, $ при $ t^J
\in (\R_+^d)^J $ выполняется равенство
\begin{equation*} \tag{2.2.10}
\Omega^{\prime l \chi_{J}} (\bm h_\sigma f, t^{J})_{L_p(\sigma^{-1} D)} =
\Omega^{\prime l \chi_{J}} (f t^{J})_{L_p(D)},
\end{equation*}
и
\begin{multline*} \tag{2.2.11}
\bm h_\sigma ((S_p^\alpha H)^\prime(D)) \subset (S_p^\alpha H)^\prime(\sigma^{-1} D),
\| \bm h_\sigma f \|_{(S_p^\alpha H)^\prime(\sigma^{-1} D)} =
\| f \|_{(S_p^\alpha H)^\prime(D)},\\ f \in (S_p^\alpha H)^\prime(D),
\alpha \in \R_+^d.
\end{multline*}

Доказательство.

В самом деле, в условиях леммы, с учётом (2.2.8) делая замену переменных в
интеграле, имеем
\begin{multline*}
\int_{ (t B^d)^{J}} \int_{(\sigma^{-1} D)_\eta^{l \chi_{J}}}
| (\Delta_\eta^{l \chi_{J}} (\bm h_\sigma f))(y)|^p dy d \eta^{J} = \\
\int_{\{(\eta^J,y): \eta^J \in (t B^d)^J, y \in (\sigma^{-1} D)_\eta^{l \chi_J}\}}
\biggl| \sum_{k^J \in (\Z_+^d(l))^J} C_{l^J}^{k^J} (-\e^J)^{l^J -k^J}
(\bm h_\sigma f)(y +k \eta \chi_J)\biggr|^p d\eta^J dy = \\
\int_{\{(\eta^J,y): \eta^J \in (t B^d)^J, y \in (\sigma^{-1} D)_\eta^{l \chi_J}\}}
\biggl| \sum_{k^J \in (\Z_+^d(l))^J} C_{l^J}^{k^J} (-\e^J)^{l^J -k^J}
f(\sigma y +\sigma k \eta \chi_J)\biggr|^p d\eta^J dy = \\
\int_{\substack{\{ (\xi^J,x): (\sigma^{-1} \xi)^J \in (t B^d)^J, \\
\sigma^{-1} x \in (\sigma^{-1} D)_{\sigma^{-1} \xi}^{l \chi_J} \}}}
\biggl| \sum_{k^J \in (\Z_+^d(l))^J} C_{l^J}^{k^J} (-\e^J)^{l^J -k^J}
f(\sigma \sigma^{-1} x +\sigma k \sigma^{-1} \xi \chi_J)\biggr|^p \times\\
| (\sigma^{-1})^{\e}| \cdot |((\sigma^{-1})^J)^{\e^J}| d\xi^J dx = \\
\int_{\substack{\{ (\xi^J,x): \xi^J \in (t B^d)^J, \\ x \in D_{\xi}^{l \chi_J} \}}}
\biggl| \sum_{k^J \in (\Z_+^d(l))^J} C_{l^J}^{k^J} (-\e^J)^{l^J -k^J}
f(x +k \xi \chi_J)\biggr|^p d\xi^J dx = \\
\int_{(t B^d)^J} \int_{ D_\xi^{l \chi_J}}
| (\Delta_\xi^{l \chi_J} f)(x)|^p dx d\xi^J.
\end{multline*}
Откуда выводим
\begin{multline*}
\Omega^{\prime l \chi_J}((\bm h_\sigma f),t^J)_{L_p(\sigma^{-1} D)} =
((2 t^J)^{-\e^J} \int_{ (t B^d)^J} \int_{(\sigma^{-1} D)_\eta^{l \chi_J}}
| (\Delta_\eta^{l \chi_J} (\bm h_\sigma f))(y)|^p dy d\eta^J)^{1 /p} = \\
((2 t^J)^{-\e^J} \int_{(t B^d)^J} \int_{ D_\xi^{l \chi_J}}
| (\Delta_\xi^{l \chi_J} f)(x)|^p dx d\xi^J)^{1 /p} = \\
\Omega^{\prime l \chi_J}(f, t^J)_{L_p(D)},
f \in L_p(D), J \subset \{1,\ldots,d\}: J \ne \emptyset, t^J \in (\R_+^d)^J,
\end{multline*}
что совпадает с (2.2.10).

Включение (2.2.11) следует из (2.2.8) и (2.2.10). $ square $

Теперь можно установить следующее предложение.

Предложение 2.2.4

При $ d \in \N $ пусть область $ D \subset \R^d $ и её открытое подмножество
$ U \subset D $ таковы, что существуют $ \delta \in \R_+^d $ и $ \sigma \in
\Sigma^d, $ для которых выполняется
включение $ (U +\sigma \delta I^d) \subset D. $ Тогда при любых
$ l \in \Z_+^d, m \in \N^d $ существует семейство линейных операторов
$$
\mathfrak E_\kappa^{d,l,m,D,U}: \bm h_\sigma^{-1}(L_1^{\loc \square}(\sigma^{-1} D)) =
L_1^{\loc \square}(D) \mapsto L_1^{\loc}(\R^d), \kappa \in \Z_+^d,
$$
обладающих следующими свойствами:

1) при $ l \in \Z_+^d, m \in \N^d, \lambda \in \Z_+^d(m) $ и $ 1 \le p \le q
\le \infty, $ а также, если $ U $ -- ограниченно, ещё и при $ 1 \le q < p \le \infty, $
существует константа $ c_1(d,l,m,D,U,\lambda,p,q) > 0 $ такая, что для $ f \in
L_p(D) $ соблюдается неравенство
\begin{equation*} \tag{2.2.12}
\| \D^\lambda \mathfrak E_0^{d,l,m,D,U} f \|_{L_q(\R^d)} \le c_1 \| f \|_{L_p(D)};
\end{equation*}

2) при $ l \in \N^d, m \in \N^d, \lambda \in \Z_+^d(m), 1 \le p < \infty,
p \le q \le \infty, $ а также, если множество $ U $ -- ограниченно, ещё и
при $ 1 \le q < p, $ существуют константы $ c_2(d,l,m,D,U,\lambda,p,q) > 0,
c_3(d,m,D,U) > 0 $ такие, что при $ \kappa \in \Z_+^d \setminus \{0\} $ для
$ f \in L_p(D) $ выполняется неравенство
\begin{multline*} \tag{2.2.13}
\| \D^\lambda \mathfrak E_{\kappa}^{d,l -\e,m,D,U} f \|_{L_q(\R^d)}
\le c_2 2^{(\kappa, \lambda +(p^{-1} -q^{-1})_+ \e)}
\Omega^{\prime l \chi_{\s(\kappa)}}(f, (c_3 2^{-\kappa})^{\s(\kappa)})_{L_p(D)};
\end{multline*}

3) при $ m \in \N^d, \alpha \in \R_+^d, 1 \le p < \infty, p \le q \le \infty, $
а также, если $ U $ -- ограниченное множество, ещё и при $ 1 \le q < p,
\lambda \in \Z_+^d(m), $ удовлетворяющих условию (2.1.68), для любой функции
$ f \in (S_p^\alpha H)^\prime(D) $ и $ l = l(\alpha) $ в $ L_q(U) $ имеет место
равенство
\begin{equation*} \tag{2.2.14}
\D^\lambda (f \mid_U) = \sum_{\kappa \in \Z_+^d} (\D^\lambda (\mathfrak E_{\kappa}^{d, l -\e,m,D,U} f)) \mid_U,
\end{equation*}

а также

4) при $ \alpha \in \R_+^d, 1 \le p < \infty, p \le q \le \infty, $ а также,
если $ U $ -- ограниченное множество, ещё и при $ 1 \le q < p, \lambda \in
\Z_+^d, $ удовлетворяющих условию (2.1.68), существует константа
$ c_4(d,\alpha,p,q,\lambda,D,U) > 0 $ такая, что для любой функции
$ f \in (S_p^\alpha H)^\prime(D) $ соблюдается неравенство
\begin{equation*} \tag{2.2.15}
\| \D^\lambda (f \mid_U) \|_{L_q(U)} \le c_4 \| f \|_{(S_p^\alpha H)^\prime(D)}.
\end{equation*}

Доказательство.

В условиях предложения, учитывая, что $ \sigma^{-1} U +\delta I^d =
\sigma^{-1} U +\sigma^{-1} \sigma \delta I^d = \sigma^{-1}(U +\sigma \delta I^d)
\subset \sigma^{-1} D, $ видим,
что в силу леммы 2.2.1 при любом $ m \in \N^d $ область $ \sigma^{-1} D $
и её открытое подмножество $ \sigma^{-1} U \subset \sigma^{-1} D $ является
$m$-правильной парой. Принимая во внимание это обстоятельство, а также
замечание перед леммой 2.2.2, определим при $ l \in \Z_+^d, m \in \N^d $
семейство операторов
$$
\mathfrak E_\kappa^{d,l,m,D,U}:
\bm h_\sigma^{-1}(L_1^{\loc \square}(\sigma^{-1} D)) =
L_1^{\loc \square}(D) \mapsto L_1^{\loc}(\R^d), \kappa \in \Z_+^d,
$$
полагая (см. (2.2.6), (2.1.3), (2.1.29))
\begin{multline*} \tag{2.2.16}
\mathfrak E_0^{d,l,m,D,U} = \bm h_\sigma^{-1}
\mathcal E_{\Kappa^0,0}^{d,l,m,\sigma^{-1} D,\sigma^{-1} U,\Nu} \bm h_\sigma =
\bm h_\sigma^{-1}
E_{\Kappa^0}^{d,l,m,\sigma^{-1} D,\sigma^{-1} U,\nu_{\Kappa^0}} \bm h_\sigma =\\
\bm h_\sigma E_{\Kappa^0}^{d,l,m,\sigma D,\sigma U,\nu_{\Kappa^0}} \bm h_\sigma, 
\mathfrak E_\kappa^{d,l,m,D,U} = \bm h_\sigma^{-1}
\mathcal E_{\Kappa^0, \kappa}^{d,l,m,\sigma^{-1} D,\sigma^{-1} U,\Nu} \bm h_\sigma =\\
\bm h_\sigma \mathcal E_{\Kappa^0, \kappa}^{d,l,m,\sigma D,\sigma U,\Nu}
\bm h_\sigma, \kappa \in \Z_+^d \setminus \{0\},
\end{multline*}
где
\begin{multline*}
\Kappa^0 = \Kappa^0(d,m,\sigma^{-1} D, \sigma^{-1} U) \in \Z_+^d, 
\Nu = \Nu^{d,m,\sigma^{-1} D,\sigma^{-1} U} =\\ 
\{ \nu_{\Kappa^0
+\kappa}^{d,m,\sigma^{-1} D,\sigma^{-1} U}: N_{\Kappa^0
+\kappa}^{d,m,\sigma^{-1} U} \mapsto \Z^d, \kappa \in \Z_+^d\}
\end{multline*}
-- объекты из определения $m$-правильной пары.

Проверим соблюдение соотношений (2.2.12) -- (2.2.14).
В условиях п. 1), используя (2.2.16), (2.2.9), (2.2.8), (2.1.7) и снова
(2.2.8), для $ f \in L_p(D) $ выводим
\begin{multline*}
\| \D^\lambda \mathfrak E_0^{d,l,m,D,U} f \|_{L_q(\R^d)} = \| \D^\lambda
(\bm h_\sigma E_{\Kappa^0}^{d,l,m,\sigma D,\sigma U,\nu_{\Kappa^0}}
\bm h_\sigma f) \|_{L_q(\R^d)} = \\
\| \sigma^\lambda \bm h_\sigma
(\D^\lambda E_{\Kappa^0}^{d,l,m,\sigma D,\sigma U,\nu_{\Kappa^0}}
\bm h_\sigma f) \|_{L_q(\R^d)} = \\
| \sigma^\lambda| \cdot \| \bm h_\sigma
(\D^\lambda E_{\Kappa^0}^{d,l,m,\sigma D,\sigma U,\nu_{\Kappa^0}}
\bm h_\sigma f) \|_{L_q(\R^d)} = \\
\| \D^\lambda (E_{\Kappa^0}^{d,l,m,\sigma^{-1} D,\sigma^{-1} U,\nu_{\Kappa^0}}
(\bm h_\sigma f)) \|_{L_q(\R^d)} \le c_1 \| \bm h_\sigma f \|_{L_p(\sigma^{-1} D)} =
c_1 \| f \|_{L_p(D)},
\end{multline*}
т.е. верно (2.2.12).

Далее, в условиях п. 2), применяя (2.2.16), (2.2.9), (2.2.8), (2.1.37),
(2.2.10), при $ \kappa \in \Z_+^d \setminus \{0\} $ для $ f \in L_p(D) $
получаем
\begin{multline*}
\| \D^\lambda \mathfrak E_{\kappa}^{d,l -\e,m,D,U} f \|_{L_q(\R^d)} =
\| \D^\lambda (\bm h_\sigma \mathcal E_{\Kappa^0, \kappa}^{d,l -\e,m,\sigma D,\sigma U,\Nu}
\bm h_\sigma f)\|_{L_q(\R^d)} = \\
\| \sigma^\lambda \bm h_\sigma
(\D^\lambda \mathcal E_{\Kappa^0, \kappa}^{d,l -\e,m,\sigma D,\sigma U,\Nu}
\bm h_\sigma f)\|_{L_q(\R^d)} = \\
| \sigma^\lambda| \cdot \| \bm h_\sigma
(\D^\lambda \mathcal E_{\Kappa^0, \kappa}^{d,l -\e,m,\sigma D,\sigma U,\Nu}
\bm h_\sigma f)\|_{L_q(\R^d)} = \\
\| \D^\lambda (\mathcal E_{\Kappa^0, \kappa}^{d,l -\e,m,\sigma^{-1} D,\sigma^{-1} U,\Nu}
\bm h_\sigma f)\|_{L_q(\R^d)} \le \\
c_2 2^{(\kappa, \lambda +(p^{-1} -q^{-1})_+ \e)}
\Omega^{\prime l \chi_{\s(\kappa)}}(\bm h_\sigma f,
(c_3 2^{-\kappa})^{\s(\kappa)})_{L_p(\sigma^{-1} D)} = \\
c_2 2^{(\kappa, \lambda +(p^{-1} -q^{-1})_+ \e)}
\Omega^{\prime l \chi_{\s(\kappa)}}(f, (c_3 2^{-\kappa})^{\s(\kappa)})_{L_p(D)},
\end{multline*}
что совпадает с (2.2.13).

Наконец, в условиях п. 3), ввиду (2.2.11) согласно п. 1) предложения
2.1.6 для $ f \in (S_p^\alpha H)^\prime(D) $ и $ l = l(\alpha) $ в
$ L_q(\sigma^{-1} U) $ имеет место равенство
\begin{equation*}
\D^\lambda ((\bm h_\sigma f) \mid_{\sigma^{-1} U}) = \sum_{\kappa \in \Z_+^d}
(\D^\lambda (\mathcal E_{\Kappa^0, \kappa}^{d, l -\e,m,\sigma^{-1} D,\sigma^{-1} U,\Nu}
(\bm h_\sigma f))) \mid_{\sigma^{-1} U},
\end{equation*}
откуда в силу (2.2.8) заключаем, что в $ l_q(\sigma^{-1} \sigma^{-1} U) =
L_q(U) $ верно равенство
\begin{multline*} \tag{2.2.17}
\bm h_\sigma (\D^\lambda ((\bm h_\sigma f) \mid_{\sigma^{-1} U})) =
\bm h_\sigma(\sum_{\kappa \in \Z_+^d}
(\D^\lambda (\mathcal E_{\Kappa^0, \kappa}^{d, l -\e,m,\sigma^{-1} D,\sigma^{-1} U,\Nu}
(\bm h_\sigma f))) \mid_{\sigma^{-1} U}) = \\
\sum_{\kappa \in \Z_+^d} \bm h_\sigma
((\D^\lambda (\mathcal E_{\Kappa^0, \kappa}^{d, l -\e,m,\sigma^{-1} D,\sigma^{-1} U,\Nu}
(\bm h_\sigma f))) \mid_{\sigma^{-1} U}).
\end{multline*}

Используя (2.2.9), (2.2.7), (2.2.6), имеем
\begin{multline*} \tag{2.2.18}
\bm h_\sigma (\D^\lambda ((\bm h_\sigma f) \mid_{\sigma^{-1} U})) =
\sigma^{-\lambda} \D^\lambda (\bm h_\sigma ((\bm h_\sigma f) \mid_{\sigma^{-1} U})) = \\
\sigma^{\lambda} \D^\lambda (\bm h_\sigma (\bm h_\sigma (f \mid_{U}))) =
\sigma^{\lambda} \D^\lambda (\bm h_\sigma^{-1} (\bm h_\sigma (f \mid_{U}))) = \\
\sigma^{\lambda} \D^\lambda (f \mid_{U}).
\end{multline*}
Применяя (2.2.7), (2.2.9), (2.2.16),  получаем
\begin{multline*} \tag{2.2.19}
\bm h_\sigma
((\D^\lambda (\mathcal E_{\Kappa^0, \kappa}^{d, l -\e,m,\sigma^{-1} D,\sigma^{-1} U,\Nu}
(\bm h_\sigma f))) \mid_{\sigma^{-1} U}) = \\
(\bm h_\sigma
(\D^\lambda (\mathcal E_{\Kappa^0, \kappa}^{d, l -\e,m,\sigma^{-1} D,\sigma^{-1} U,\Nu}
(\bm h_\sigma f)))) \mid_{\sigma^{-1} \sigma^{-1} U} = \\
(\bm h_\sigma
(\D^\lambda (\mathcal E_{\Kappa^0, \kappa}^{d, l -\e,m,\sigma D,\sigma U,\Nu}
(\bm h_\sigma f)))) \mid_{U} = \\
(\sigma^{-\lambda} \D^\lambda (\bm h_\sigma
(\mathcal E_{\Kappa^0, \kappa}^{d, l -\e,m,\sigma D,\sigma U,\Nu}
(\bm h_\sigma f)))) \mid_{U} = \\
\sigma^{\lambda} (\D^\lambda (\bm h_\sigma
(\mathcal E_{\Kappa^0, \kappa}^{d, l -\e,m,\sigma D,\sigma U,\Nu}
(\bm h_\sigma f)))) \mid_{U} = \\
\sigma^{\lambda} (\D^\lambda ((\bm h_\sigma
\mathcal E_{\Kappa^0, \kappa}^{d, l -\e,m,\sigma D,\sigma U,\Nu}
\bm h_\sigma) f)) \mid_{U} = \\
\sigma^{\lambda} (\D^\lambda (\mathfrak E_{\kappa}^{d, l -\e,m,D,U} f)) \mid_U.
\end{multline*}
Вывод (2.2.19) проведен при $ \kappa \ne 0. $ Для получения (2.2.19) при
$ \kappa = 0 $ в проведенных выкладках следует заменить
$ \mathcal E_{\Kappa^0, \kappa}^{d,l -\e,m,\sigma D,\sigma U,\Nu} $ на
$ E_{\Kappa^0}^{d,l -\e,m,\sigma D,\sigma U,\nu_{\Kappa^0}}. $
Подставляя (2.2.18) и (2.2.19) в (2.2.17), в $ L_q(U) $ приходим к
равенству
\begin{equation*}
\sigma^\lambda \D^\lambda (f \mid_U) = \sum_{\kappa \in \Z_+^d}
\sigma^\lambda (\D^\lambda (\mathfrak E_{\kappa}^{d, l -\e,m,D,U} f)) \mid_U,
\end{equation*}
из которого следует (2.2.14).

 Остаётся убедиться в справедливости (2.2.15). Для этого в
 условиях п. 4) предложения, фиксируем $ m
\in \N^d $ так, чтобы $ \lambda \in \Z_+^d(m). $ Учитывая, что $
(\sigma^{-1} U +\delta I^d) \subset \sigma^{-1} D, $ в силу леммы
2.2.1 $ (\sigma^{-1} D, \sigma^{-1} U) $ является $m$-правильной
парой. Принимая во внимание эти обстоятельства, а также (2.2.11),
на основании (2.1.69) заключаем, что для $ f \in (S_p^\alpha
H)^\prime(D) $ имеет место соотношение
\begin{equation*}
\| \D^\lambda ((\bm h_\sigma f) \mid_{\sigma^{-1} U}) \|_{L_q(\sigma^{-1} U)}
\le c_4 \| \bm h_\sigma f \|_{(S_p^\alpha H)^\prime(\sigma^{-1} D)}.
\end{equation*}
Отсюда, учитывая (2.2.11), а также тот факт, что вследствие (2.2.8), (2.2.18)
выполняется соотношение
\begin{multline*}
\| \D^\lambda ((\bm h_\sigma f) \mid_{\sigma^{-1} U}) \|_{L_q(\sigma^{-1} U)} =
\| \bm h_\sigma (\D^\lambda ((\bm h_\sigma f) \mid_{\sigma^{-1} U}))
\|_{L_q(\sigma^{-1} \sigma^{-1} U)} = \\
\| \bm h_\sigma (\D^\lambda ((\bm h_\sigma f)
\mid_{\sigma^{-1} U})) \|_{L_q(U)} = \\
\| \sigma^{\lambda} \D^\lambda (f \mid_{U}) \|_{L_q(U)} =
\| \D^\lambda (f \mid_{U}) \|_{L_q(U)},
\end{multline*}
получаем (2.2.15). $ \square $

Теорема 2.2.5

При $ d \in \N $ пусть $ D $ -- область в $ \R^d, $ для которой существует
система открытых подмножеств $ \{ U_i \subset D, i =1,\ldots,\mathcal I\} $

такая, что при $ i =1,\ldots,\mathcal I $ существуют $ \delta^i \in \R_+^d $ и
$ \sigma^i \in \Sigma^d, $ для которых имеет место включение
$ (U_i +\sigma^i \delta^i I^d) \subset D, $ и $ D = \cup_{i =1}^{\mathcal I} U_i. $
Тогда при $ \alpha \in \R_+^d, 1 \le p < \infty, p \le q \le \infty, $ а также,
если $ D $ -- ограниченная область, ещё и при $ 1 \le q < p, \lambda \in
\Z_+^d, $ удовлетворяющих условию (2.1.68), существует константа
$ c_5(d,\alpha,p,q,\lambda,D) > 0 $ такая, что для любой функции
$ f \in (S_p^\alpha H)^\prime(D) $ верно неравенство
\begin{equation*} \tag{2.2.20}
\| \D^\lambda f \|_{L_q(D)} \le c_5 \| f \|_{(S_p^\alpha H)^\prime(D)}.
\end{equation*}

Доказательство.

Сначала покажем, что в условиях теоремы для $ f \in (S_p^\alpha H)^\prime(D) $
(обобщённая) производная $ \D^\lambda f \in L_q(D). $ Для этого при
$ i =1,\ldots,\mathcal I, $ имея в виду
(2.2.15), возьмём функцию $ \bm f_i \in L_q(U_i), $ для которой выполняется равенство
\begin{equation*}
\langle \D^\lambda (f \mid_{U_i}), \phi \rangle = \int_{U_i} \bm f_i \phi dx
\end{equation*}
для всех $ \phi \in C_0^\infty(U_i). $

Тогда при $ i,j =1,\ldots,\mathcal I: U_i \cap U_j \ne \emptyset, $ почти для
всех $ x \in U_i \cap U_j $ соблюдается равенство
\begin{equation*} \tag{2.2.21}
\bm f_i(x) = \bm f_j(x).
\end{equation*}
в самом деле, в описанных условиях для $ \phi \in C_0^\infty(U_i \cap U_j) =
C_0^\infty(U_i) \cap C_0^\infty(U_j) \subset C_0^\infty(D) $ имем
\begin{multline*}
\int_{U_i \cap U_j} \bm f_i \phi dx =
\int_{U_i} \bm f_i \phi dx =
\langle \D^\lambda (f \mid_{U_i}), \phi \rangle =
\langle (\D^\lambda f) \mid_{U_i}, \phi \rangle = \\
\langle \D^\lambda f, \phi \rangle =
\langle (\D^\lambda f) \mid_{U_j}, \phi \rangle =
\langle \D^\lambda (f \mid_{U_j}), \phi \rangle =
\int_{U_j} \bm f_j \phi dx =
\int_{U_i \cap U_j} \bm f_j \phi dx,
\end{multline*}
откуда в силу леммы Дюбуа-Реймона следует (2.2.21).

Учитывая, что $ D = \cup_{i =1}^{\mathcal I} U_i, $ определим функцию
$ \bm f: D \mapsto \R, $ полагая $ \bm f(x) = \bm f_i(x) $ почти
для всех $ x \in U_i, i =1,\ldots,\mathcal I. $
Ввиду (2.2.21) функция $ \bm f $ определена корректно на $ D. $ При этом
понятно, что $ \bm f \in L_q(D) $ и выполняется неравенство
$$
| \bm f | \le (\sum_{i =1}^{\mathcal I} \chi_{U_i}) | \bm f | =
| (\sum_{i =1}^{\mathcal I} \chi_{U_i}) \bm f | =
| \sum_{i =1}^{\mathcal I} (\chi_{U_i} \bm f) |,
$$
откуда
\begin{equation*} \tag{2.2.22}
\| \bm f \|_{L_q(D)} \le \| \sum_{i =1}^{\mathcal I} (\chi_{U_i} \bm f)\|_{L_q(D)} \le
\sum_{i =1}^{\mathcal I} \| \chi_{U_i} \bm f\|_{L_q(D)} = \\
\sum_{i =1}^{\mathcal I} \| \bm f \mid_{U_i}\|_{L_q(U_i)} =
\sum_{i =1}^{\mathcal I} \| \bm f_i\|_{L_q(U_i)}.
\end{equation*}

Теперь для $ \phi \in C_0^\infty(D) $ построим систему функций
$ \{ g_i \in C_0^\infty(U_i), i =1,\ldots,\mathcal I\}, $ обладающую тем
свойством, что для $ x \in \supp \phi $ соблюдается равенство
$ \sum_{i =1}^{\mathcal I} g_i(x) =1. $
Тогда, поскольку $ (g_i \phi) \in C_0^\infty(U_i), i =1,\ldots,\mathcal I, $
имеем
\begin{multline*}
\langle \D^\lambda f, \phi \rangle =
\langle \D^\lambda f, ( \sum_{i =1}^{\mathcal I} g_i) \phi \rangle =
\langle \D^\lambda f, \sum_{i =1}^{\mathcal I} (g_i \phi) \rangle = \\
\sum_{i =1}^{\mathcal I} \langle \D^\lambda f, (g_i \phi) \rangle =
\sum_{i =1}^{\mathcal I} \langle (\D^\lambda f) \mid_{U_i}, (g_i \phi) \rangle = \\
\sum_{i =1}^{\mathcal I} \langle \D^\lambda (f \mid_{U_i}), (g_i \phi) \rangle =
\sum_{i =1}^{\mathcal I} \int_{U_i} \bm f_i (g_i \phi) dx = \\
\sum_{i =1}^{\mathcal I} \int_{U_i} \bm f (g_i \phi) dx =
\sum_{i =1}^{\mathcal I} \int_{D} \bm f (g_i \phi) dx = \\
\int_{D} \sum_{i =1}^{\mathcal I} \bm f (g_i \phi) dx =
\int_{D} \bm f (\sum_{i =1}^{\mathcal I} g_i) \phi dx = 
\int_{D} \bm f \phi dx.
\end{multline*}

А это значит, что $ \D^\lambda f \in L_q(D). $

Ввиду сказанного, применяя (2.2.22), (2.2.15), приходим к неравенству
\begin{multline*}
\| \D^\lambda f \|_{L_q(D)} = \| \bm f \|_{L_q(D)} \le
\sum_{i =1}^{\mathcal I} \| \bm f_i\|_{L_q(U_i)} = 
\sum_{i =1}^{\mathcal I} \| \D^\lambda (f \mid_{U_i}) \|_{L_q(U_i)} \le \\
\sum_{i =1}^{\mathcal I} c_4(d,\alpha,p,q,\lambda,D,U_i)
\| f \|_{(S_p^\alpha H)^\prime(D)} = c_5 \| f \|_{(S_p^\alpha H)^\prime(D)},
\end{multline*}
совпадающему с (2.2.20). $ \square $

Лемма 2.2.6

Пусть $ D $ -- область в $ \R^d, $ а $ \{U_i \subset D, i =1,\ldots,\mathcal I\} $ --
конечный набор открытых подмножеств $ D, $ для которых существует число
$ \bm \delta > 0 $ такое, что набор открытых множеств
$$
V_i = \{ x \in U_i: \rho_i(x) = \inf_{ y \in D \setminus U_i} \| x -y\| >
\bm \delta\}, i =1,\ldots,\mathcal I,
$$
образует покрытие области $ D. $
Тогда существует система функций $ g_i: \R^d \mapsto \R, i =1,\ldots,\mathcal I, $
обладающих следующими свойствами:

при $ i =1,\ldots,\mathcal I $ выполняются соотношения
\begin{equation*}
g_i(x) \ge 0, x \in \R^d;
\end{equation*}
\begin{equation*} \tag{2.2.23}
g_i \in C^\infty(\R^d), \\
\D^\lambda g_i \in L_\infty(\R^d), \lambda \in \Z_+^d;
\end{equation*}
\begin{equation*} \tag{2.2.24}
g_i(x) =0 \text{ при } x \in D \setminus U_i;
\end{equation*}
\begin{equation*} \tag{2.2.25}
\sum_{i =1}^{\mathcal I} g_i(x) =1, x \in D.
\end{equation*}

Доказательство.

В условиях леммы при $ i =1,\lдots,\mathcal I $ рассмотрим непрерывные
функции
$$
\rho_i(x) = \inf_{ y \in D \setminus U_i} \| x -y\|, x \in \R^d, \\
\rho_0(x) = \inf_{ y \in D} \| x -y\|, x \in \R^d,
$$
и определим множества
$$
\mathcal V_i = \{ x \in \R^d: \rho_i(x) > \bm \delta /3\}, i =0,1,\ldots,\mathcal I.
$$
Тогда для $ x \in \R^d \setminus \mathcal V_0 $ справедливо неравенство
$ \rho_0(x) \le \bm delta /3 < \bm \delta /2. $ Поэтому существует
$ y \in D, $ для которого верно неравенство $ \| x -y\| < \bm \delta /2. $
Учитывая, что
$$
D \subset (\cup_{i =1}^{\mathcal I} V_i) \subset (\cup_{i =1}^{\mathcal I} U_i)
\subset D,
$$
возьмём $ i \in \{1,\ldots, \mathcal I\} $ такое, что $ y \in V_i. $ Тогда для
любого $ z \in (D \setminus U_i) $ имеем
\begin{multline*}
\| x -z\| = \| x -y +y -z\| \ge \| y -z\| -\| y -x\| = \| y -z\| -\| x -y\| \ge\\
\rho_i(y) -\| x -y\| > \bm \delta -\bm \delta /2 = \bm \delta /2,
\end{multline*}
откуда
$$
\rho_i(x) = \inf_{z \in (D \setminus U_i)} \| x -z\| \ge \bm \delta /2 >
\bm \delta /3,
$$
т.е. $ x \in \mathcal V_i, $ а, значит, $ (\R^d \setminus \mathcal V_0) \subset
(\cup_{j =1}^{\mathcal I} \mathcal V_j), $ и, следовательно,
\begin{equation*} \tag{2.2.26}
\R^d = \cup_{i =0}^{\mathcal I} \mathcal V_i.
\end{equation*}

Далее, возьмём функцию $ w \in C^\infty(\R^d) $ такую,что
\begin{equation*} \tag{2.2.27}
w(x) > 0, \text{ при } \|x\| < 1, \\
w(x) = 0, \text{ при } \|x\| \ge 1,
\end{equation*}
и фиксируя $ \delta > 0, $ для которого выполняется условие
$$
\delta < \bm \delta /3,
$$
определим при $ i =0,1,\ldots, \mathcal I $ функцию $ \psi_i \in C^\infty(\R^d) $
равенством
$$
\psi_i(x) = \int_{\mathcal V_i} w_\delta(y -x) dy \ge 0, x \in \R^d,
$$
где
$$
w_\delta(x) = w(\delta^{-1} x), x \in \R^d.
$$
Обозначая через
$$
B_0^d = \{x \in \R^d: \|x\| < 1\},
$$
ввиду (2.2.27) для $ x, y \in \R^d $ получаем соотношение
\begin{equation*}
w_\delta(y -x) > 0, \text{ при } \delta^{-1} \|y -x\| < 1, \\
w_\delta(y -x) = 0, \text{ при } \delta^{-1} \|y -x\| \ge 1,
\end{equation*}
или
\begin{equation*}
w_\delta(y -x) > 0, \text{ при } \|y -x\| < \delta, \\
w_\delta(y -x) = 0, \text{ при } \|y -x\| \ge \delta,
\end{equation*}
т.е.
\begin{equation*} \tag{2.2.28}
w_\delta(y -x) > 0, \text{ при } y \in (x +\delta B_0^d), \\
w_\delta(y -x) = 0, \text{ при } y \notin (x +\delta B_0^d).
\end{equation*}
Принимая во внимание (2.2.28), при $ i =0,1,\ldots,\mathcal I $ для
$ x \in \R^d $ имеем
\begin{equation*} \tag{2.2.29}
\psi_i(x) = \int_{\mathcal V_i \cap (x +\delta B_0^d)} w_\delta(y -x) dy.
\end{equation*}

Установим некоторые полезные для нас свойства функций $ \psi_i, i =0,\ldots,\mathcal I. $
При $ i =0,1,\ldots,\mathcal I $ для $ x \in \mathcal V_i $ открытое
множество $ \mathcal V_i \cap (x +\delta B_0^d) \ne \emptyset, $ а, значит,
$ \mes(\mathcal V_i \cap (x +\delta B_0^d)) > 0, $ что в силу (2.2.28),
(2.2.29) даёт
\begin{equation*} \tag{2.2.30}
\psi_i(x) = \int_{ \mathcal V_i \cap (x +\delta B_0^d)} w_\delta(y -x) dy > 0.
\end{equation*}

Из (2.2.26) и (2.2.30) вытекает, что для каждого $ x \in \R^d $ существуеТ
$ i =0,1,\ldots,\mathcal I, $ для которого $ \psi_i(x) > 0, $ а, значит,
$ \sum_{i =0}^{\mathcal I} \psi_i(x) > 0, $ ибо $ \psi_j(x) \ge 0, j =0,\ldots,
\mathcal I. $
Отсюда, полагая
$$
\psi = \sum_{i =0}^{\mathcal I} \psi_i,
$$
и обознначая при $ i =0,\ldots,\mathcal I $ функцию
$$
g_i = \psi_i / \psi,
$$
видим, что при $ i =0,\ldots,\mathcal I $ имеют место соотношения
$ g_i(\cdot) \ge 0, g_i \in C^\infty(\R^d), $ т.е. соблюдается первое
включение в (2.2.23).
При этом для $ x \in \R^d $ выполняется равенство
\begin{equation*} \tag{2.2.31}
\sum_{i =0}^{\mathcal I} g_i(x) = \sum_{i =0}^{\mathcal I} (\psi_i(x)
/ \psi(x)) = (\sum_{i =0}^{\mathcal I} \psi_i(x)) / \psi(x) = 1.
\end{equation*}
Кроме того, заметим, что для $ x \in D $ пересечение $ \mathcal V_0 \cap
(x +\delta B_0^d) = \emptyset, $ поскольку, если это не так, то для
$ y \in (\mathcal V_0 \cap (x +\delta B_0^d)) $ получаем, что
$$
\delta > \| x -y\| \ge \rho_0(y) > \bm \delta /3 > \delta,
$$
а это неверно. Принимая во внимание сказанное, согласно (2.2.29) заключаем,
что для $ x \in D $ значение $ \psi_0(x) = 0, $ а, значит, и $ g_0(x) = 0, $
что в соединении с (2.2.31) приводит к (2.2.25).

Далее, при $ i =1,\ldots, \mathcal I $ для $ x \in (D \setminus U_i) $
пересечение $ \mathcal V_i \cap (x +\delta B_0^d) = \emptyset,$ поскольку,
если это не так, то для $ y \in (\mathcal V_i \cap (x +\delta B_0^d)) $
имеем
$$
\delta > \| x -y\| \ge \rho_i(y) > \bm \delta /3 > \delta,
$$
что неверно. Принимая во внимание сказанное, в силу (2.2.29) получаем,
что при $ i =1,\ldots,\mathcal I $ для $ x \in (D \setminus U_i) $ значение
$ \psi_i(x) = 0, $ а, следовательно, и $ g_i(x) = 0, $ что показывает
справедливость (2.2.24).

Остаётся проверить второе включение в (2.2.23). Для этого отметим некоторые факты.
Прежде всего, при $ i =0,\ldots,\mathcal I, \lambda \in \Z_+^d $ ввиду (2.1.4)
имеем
\begin{equation*} \tag{2.2.32}
\D^\lambda g_i = \D^\lambda (\psi_i / \psi) = \sum_{ \mu \in \Z_+^d(\lambda)}
C_\lambda^\mu \D^{\lambda -\mu} \psi_i \D^\mu (1 / \psi).
\end{equation*}

Кроме того, при $ i =0,\ldots,\mathcal I, \lambda \in \Z_+^d $ для $ x \in \R^d $
справедливо неравенство
\begin{multline*}
| \D^\lambda \psi_i(x) | = | \D^\lambda(\int_{\mathcal V_i} w_\delta(y -x) dy)| =
| \int_{\mathcal V_i} \frac{\D^{|\lambda|}} {\D x^\lambda} (w_\delta(y -x)) dy| = \\
| \int_{\mathcal V_i} (-\e)^\lambda (\D^\lambda w_\delta)(y -x) dy| =
| \int_{\mathcal V_i} \D^\lambda w_\delta(y -x) dy| \le \\
\int_{\mathcal V_i} | \D^\lambda w_\delta(y -x)| dy \le
\int_{\R^d} | \D^\lambda w_\delta(y -x)| dy =
\int_{\R^d} | \D^\lambda w_\delta(y)| dy < \infty,
\end{multline*}
т.е.
\begin{equation*} \tag{2.2.33}
\D^\lambda \psi_i \in L_\infty(\R^d), \\
\D^\lambda \psi = (\sum_{j =0}^{\mathcal I} \D^\lambda \psi_j) \in L_\infty(\R^d).
\end{equation*}

И ещё, для $ x \in \R^d, $ учитывая (2.2.26), выполняется неравенство
\begin{multline*}
\psi(x) = \sum_{i =0}^{\mathcal I} \psi_i(x) =
\sum_{i =0}^{\mathcal I} \int_{\mathcal V_i} w_\delta(y -x) dy = 
\sum_{i =0}^{\mathcal I} \int_{\R^d} \chi_{\mathcal V_i}(y) w_\delta(y -x) dy =\\
\int_{\R^d} (\sum_{i =0}^{\mathcal I} \chi_{\mathcal V_i}(y)) w_\delta(y -x) dy \ge 
\int_{\R^d} 1 \cdot w_\delta(y -x) dy = \int_{\R^d} w_\delta(y) dy > 0,
\end{multline*}
а, значит,
\begin{equation*} \tag{2.2.34}
1 / \psi \in L_\infty(\R^d).
\end{equation*}

Наконец, проверим, что при $ \lambda \in \Z_+^d $ имеет место равенство вида
\begin{multline*} \tag{2.2.35}
\D^\lambda (1 / \psi) = \sum_{\{(\mu^1,\ldots,\mu^{|\lambda|}): \mu^s \in
\Z_+^d, s =1,\ldots,|\lambda|; \sum_{s =1}^{|\lambda|} \mu^s = \lambda\}} 
a_{\mu^1,\ldots,\mu^{|\lambda|}}
(\prod_{s =1}^{|\lambda|} \D^{\mu^s} \psi) / \psi^{|\lambda| +1}, \\
a_{\mu^1,\ldots,\mu^{|\lambda|}} \in \R, \mu^s \in \Z_+^d, s =1,\ldots,|\lambda|:
\sum_{s =1}^{|\lambda|} \mu^s = \lambda.
\end{multline*}
Соотношение (2.2.35) проверим по индукции относительно $ l = |\lambda|. $
При $ l =0,1 $ равенство (2.2.35) очевидно. Пусть теперь $ l \in \N $ и
$ |\lambda| = l. $ Тогда при $ j =1,\ldots,d $ имеем
\begin{multline*} \tag{2.2.36}
\D^{\lambda +e_j}(1 / \psi) = \D^{e_j} (\D^\lambda (1 / \psi)) = \\
\D^{e_j}\biggl(\sum_{\{(\mu^1,\ldots,\mu^l): \mu^s \in \Z_+^d, s =1,\ldots,l; \sum_{s =1}^l \mu^s = \lambda\}}
a_{\mu^1,\ldots,\mu^l}
(\prod_{s =1}^l \D^{\mu^s} \psi) / \psi^{l +1}\biggr) = \\
\sum_{\{(\mu^1,\ldots,\mu^l): \mu^s \in \Z_+^d, s =1,\ldots,l; \sum_{s =1}^l \mu^s = \lambda\}}
a_{\mu^1,\ldots,\mu^l}
\D^{e_j} \biggl((\prod_{s =1}^l \D^{\mu^s} \psi) / \psi^{l +1}\biggr),
\end{multline*}
причём для $ (\mu^1,\ldots,\mu^l): \mu^s \in \Z_+^d, s =1,\ldots,l;
\sum_{s =1}^l \mu^s = \lambda, $ имеет место равенство
\begin{multline*} \tag{2.2.37}
\D^{e_j} ((\prod_{s =1}^l \D^{\mu^s} \psi) / \psi^{l +1}) = 
\sum_{k =1}^l (\D^{\mu^k +e_j} \psi (\prod_{s =1,\ldots,l: s \ne k} \D^{\mu^s} \psi) / \psi^{l +1}) + \\
(\prod_{s =1}^l \D^{\mu^s} \psi) \D^{e_j} (1 / \psi^{l +1}) = 
\sum_{k =1}^l (\D^{\mu^k +e_j} \psi (\prod_{s =1,\ldots,l: s \ne k} \D^{\mu^s} \psi) \psi / \psi^{l +2}) - \\
(\prod_{s =1}^l \D^{\mu^s} \psi) (l +1) \D^{e_j} \psi (1 / \psi^{l +2}).
\end{multline*}
Объединяя (2.2.36), (2.2.37), убеждаемся в справедливости (2.2.35) при
$ |\lambda| = l +1. $ Сопоставляя (2.2.32) -- (2.2.35), заключаем, что выполняется второе включение в (2.2.23). $ \square $

Теорема 2.2.7

При $ d \in \N $ пусть $ D $ -- область в $ \R^d, $ для которой существует
система открытых подмножеств $ \{ U_i \subset D, i =1,\ldots,\mathcal I\}, $
удовлетворяющих условиям теоремы 2.2.5 и леммы 2.2.6.
Тогда при любых $ l \in \Z_+^d, m \in \N^d $ существует семейство линейных
операторов
$$
\mathfrak E_\kappa^{d,l,m,D}: L_1^{\loc \square}(D) \mapsto L_1^{\loc}(\R^d), \kappa \in \Z_+^d,
$$
обладающих следующими свойствами:

1) при $ l \in \Z_+^d, m \in \N^d, \lambda \in \Z_+^d(m) $ и $ 1 \le p \le q
\le \infty, $ а также, если $ D $ -- ограниченная область, ещё и при $ 1 \le
q < p \le \infty, $ существует константа $ c_6(d,l,m,D,\lambda,p,q) > 0 $ такая,
что для $ f \in L_p(D) $ соблюдается неравенство
\begin{equation*} \tag{2.2.38}
\| \D^\lambda \mathfrak E_0^{d,l,m,D} f \|_{L_q(\R^d)} \le c_6 \| f \|_{L_p(D)};
\end{equation*}

2) при $ l \in \N^d, m \in \N^d, \lambda \in \Z_+^d(m), 1 \le p < \infty,
p \le q \le \infty, $ а также, если область $ D $ -- ограниченна, ещё и
при $ 1 \le q < p, $ существуют константы $ c_7(d,l,m,D,\lambda,p,q) > 0,
c_8(d,m,D) > 0 $ такие, что при $ \kappa \in \Z_+^d \setminus \{0\} $ для
$ f \in L_p(D) $ выполняется неравенство
\begin{equation*} \tag{2.2.39}
\| \D^\lambda \mathfrak E_{\kappa}^{d,l -\e,m,D} f \|_{L_q(\R^d)}
\le c_7 2^{(\kappa, \lambda +(p^{-1} -q^{-1})_+ \e)}
\Omega^{\prime l \chi_{\s(\kappa)}}(f, (c_8 2^{-\kappa})^{\s(\kappa)})_{L_p(D)};
\end{equation*}

3) при $ m \in \N^d, \alpha \in \R_+^d, 1 \le p < \infty, p \le q \le \infty, $
а также, если $ D $ -- ограниченная область, ещё и при $ 1 \le q < p,
\lambda \in \Z_+^d(m), $ удовлетворяющих условию (2.1.68), и $ l = l(\alpha) $
для любой функции $ f \in (S_p^\alpha H)^\prime(D) $ в $ L_q(D) $
соблюдается равенство
\begin{equation*} \tag{2.2.40}
\D^\lambda f = \sum_{\kappa \in \Z_+^d} (\D^\lambda (\mathfrak E_{\kappa}^{d, l -\e,m,D} f)) \mid_D.
\end{equation*}

Доказательство.

В условиях теоремы, опираясь на лемму 2.2.6, построим систему функций $ g_i,
i =1,\ldots,\mathcal I, $ обладающих описанными в лемме свойствами.
При $ l \in \Z_+^d, m \in \N^d, \kappa \in \Z_+^d $ определим линейный
оператор
$$
\mathfrak E_\kappa^{d,l,m,D}: L_1^{\loc \square}(D) \mapsto L_1^{\loc}(\R^d),
$$
полагая для $ f \in L_1^{\loc \square}(D) $ значение
\begin{equation*}
\mathfrak E_\kappa^{d,l,m,D} f = \sum_{i =1}^{\mathcal I} g_i
(\mathfrak E_\kappa^{d,l,m,D,U_i} f) \text{ (см. предложение 2.2.4 и (2.2.23))}.
\end{equation*}

Проверим, что определённое таким образом семейство операторов обладает
требуемыми свойствами.

В условиях п. 1) для $ f \in L_p(D), $ используя (2.1.4), (2.2.23), (2.2.12),
находим, что
\begin{multline*}
\| \D^\lambda \mathfrak E_0^{d,l,m,D} f \|_{L_q(\R^d)} =
\| \D^\lambda (\sum_{i =1}^{\mathcal I} g_i (\mathfrak E_0^{d,l,m,D,U_i} f))
\|_{L_q(\R^d)} = \\
\| \sum_{i =1}^{\mathcal I} \D^\lambda (g_i (\mathfrak E_0^{d,l,m,D,U_i} f))
\|_{L_q(\R^d)} = \\
\| \sum_{i =1}^{\mathcal I} \sum_{ \mu \in \Z_+^d(\lambda)} C_\lambda^\mu
\D^{\lambda -\mu} g_i \D^\mu(\mathfrak E_0^{d,l,m,D,U_i} f) \|_{L_q(\R^d)} \le \\
\sum_{i =1}^{\mathcal I} \sum_{ \mu \in \Z_+^d(\lambda)} C_\lambda^\mu
\| \D^{\lambda -\mu} g_i \D^\mu(\mathfrak E_0^{d,l,m,D,U_i} f) \|_{L_q(\R^d)} \le \\
\sum_{i =1}^{\mathcal I} \sum_{ \mu \in \Z_+^d(\lambda)} C_\lambda^\mu
\| \D^{\lambda -\mu} g_i\|_{L_\infty(\R^d)}
\| \D^\mu(\mathfrak E_0^{d,l,m,D,U_i} f) \|_{L_q(\R^d)} \le \\
\sum_{i =1}^{\mathcal I} \sum_{ \mu \in \Z_+^d(\lambda)} C_\lambda^\mu
\| \D^{\lambda -\mu} g_i\|_{L_\infty(\R^d)}
c_1(d,l,m,D,U_i,\mu,p,q) \| f \|_{L_p(D)} = \\
\biggl(\sum_{i =1}^{\mathcal I} \sum_{ \mu \in \Z_+^d(\lambda)} C_\lambda^\mu
\| \D^{\lambda -\mu} g_i\|_{L_\infty(\R^d)} c_1(d,l,m,D,U_i,\mu,p,q)\biggr)
\| f \|_{L_p(D)} = c_6 \| f \|_{L_p(D)},
\end{multline*}
что совпадает с (2.2.38).

Далее, в условиях п. 2) при $ \kappa \in \Z_+^d \setminus \{0\} $ для
$ f \in L_p(D), $ применяя (2.1.4), (2.2.23), (2.2.13), выводим

\begin{multline*}
\| \D^\lambda \mathfrak E_{\kappa}^{d,l -\e,m,D} f \|_{L_q(\R^d)} =
\| \D^\lambda (\sum_{i =1}^{\mathcal I} g_i (\mathfrak E_\kappa^{d,l -\e,m,D,U_i} f))
\|_{L_q(\R^d)} = \\
\| \sum_{i =1}^{\mathcal I} \D^\lambda (g_i (\mathfrak E_\kappa^{d,l -\e,m,D,U_i} f))
\|_{L_q(\R^d)} = \\
\| \sum_{i =1}^{\mathcal I} \sum_{ \mu \in \Z_+^d(\lambda)} C_\lambda^\mu
\D^{\lambda -\mu} g_i \D^\mu(\mathfrak E_\kappa^{d,l -\e,m,D,U_i} f) \|_{L_q(\R^d)} \le \\
\sum_{i =1}^{\mathcal I} \sum_{ \mu \in \Z_+^d(\lambda)} C_\lambda^\mu
\| \D^{\lambda -\mu} g_i \D^\mu(\mathfrak E_\kappa^{d,l -\e,m,D,U_i} f) \|_{L_q(\R^d)} \le \\
\sum_{i =1}^{\mathcal I} \sum_{ \mu \in \Z_+^d(\lambda)} C_\lambda^\mu
\| \D^{\lambda -\mu} g_i\|_{L_\infty(\R^d)}
\| \D^\mu(\mathfrak E_\kappa^{d,l -\e,m,D,U_i} f) \|_{L_q(\R^d)} \le \\
\sum_{i =1}^{\mathcal I} \sum_{ \mu \in \Z_+^d(\lambda)} C_\lambda^\mu
\| \D^{\lambda -\mu} g_i\|_{L_\infty(\R^d)}
c_2(d,l,m,D,U_i,\mu,p,q)\times\\ 
2^{(\kappa, \mu +(p^{-1} -q^{-1})_+ \e)}
\Omega^{\prime l \chi_{\s(\kappa)}}(f, (c_3(d,m,D,U_i) 2^{-\kappa})^{\s(\kappa)})_{L_p(D)} \le \\
\sum_{i =1}^{\mathcal I} \sum_{ \mu \in \Z_+^d(\lambda)} C_\lambda^\mu
\| \D^{\lambda -\mu} g_i\|_{L_\infty(\R^d)}
c_2(d,l,m,D,U_i,\mu,p,q)\times\\ 
2^{(\kappa, \lambda +(p^{-1} -q^{-1})_+ \e)}
\Omega^{\prime l \chi_{\s(\kappa)}}(f, (c_3(d,m,D,U_i) 2^{-\kappa})^{\s(\kappa)})_{L_p(D)} \le \\
\sum_{i =1}^{\mathcal I} \sum_{ \mu \in \Z_+^d(\lambda)} C_\lambda^\mu
\| \D^{\lambda -\mu} g_i\|_{L_\infty(\R^d)}
c_9(d,l,m,D,U_i,\mu,p,q)\times\\
 2^{(\kappa, \lambda +(p^{-1} -q^{-1})_+ \e)}
\Omega^{\prime l \chi_{\s(\kappa)}}(f, (c_8(d,m,D) 2^{-\kappa})^{\s(\kappa)})_{L_p(D)} = \\
(\sum_{i =1}^{\mathcal I} \sum_{ \mu \in \Z_+^d(\lambda)} C_\lambda^\mu
\| \D^{\lambda -\mu} g_i\|_{L_\infty(\R^d)}
c_9(d,l,m,D,U_i,\mu,p,q))\times\\ 
2^{(\kappa, \lambda +(p^{-1} -q^{-1})_+ \e)}
\Omega^{\prime l \chi_{\s(\kappa)}}(f, (c_8(d,m,D) 2^{-\kappa})^{\s(\kappa)})_{L_p(D)} = \\
c_7 2^{(\kappa, \lambda +(p^{-1} -q^{-1})_+ \e)}
\Omega^{\prime l \chi_{\s(\kappa)}}(f, (c_8 2^{-\kappa})^{\s(\kappa)})_{L_p(D)},
\end{multline*}
где $ c_8(d,m,D) = \max_{i =1,\ldots,\mathcal I} c_3(d,m,D,U_i), $ т.е. имеет
место (2.2.39).

Наконец, пусть $ m \in \N^d, \alpha \in \R_+^d, 1 \le p < \infty, p \le q \le \infty, $
а также, если $ D $ -- ограниченная область, пусть ещё и $ 1 \le q < p,
\lambda \in \Z_+^d(m) $ удовлетворяют условию (2.1.68), а $ l = l(\alpha). $
Сначала установим соблюдение в $ L_q(D) $ равенства (2.2.40) для любой функции
$ f \in (S_p^\alpha H)^\prime(D) $ при $ q = p, \lambda =0. $
Для этого, принимая во внимание (2.2.25), (2.2.24), (2.2.23), а также
(2.2.14), для $ f \in (S_p^\alpha H)^\prime(D) $ при $ k \in \Z_+^d $ имеем
\begin{multline*}
\| f -\sum_{\kappa \in \Z_+^d(k)} (\mathfrak E_{\kappa}^{d,l -\e,m,D} f) \mid_D\|_{L_p(D)} =\\
\| f (\sum_{i =1}^{\mathcal I} g_i) \mid_D -\sum_{\kappa \in \Z_+^d(k)}
(\sum_{i =1}^{\mathcal I} g_i
(\mathfrak E_\kappa^{d,l -\e,m,D,U_i} f)) \mid_D\|_{L_p(D)} = \\
\| \sum_{i =1}^{\mathcal I} ((g_i \mid_D) f) -\sum_{\kappa \in \Z_+^d(k)}
\sum_{i =1}^{\mathcal I} (g_i \mid_D)
(\mathfrak E_\kappa^{d,l -\e,m,D,U_i} f) \mid_D\|_{L_p(D)} = \\
\| \sum_{i =1}^{\mathcal I} ((g_i \mid_D) f) -\sum_{i =1}^{\mathcal I}
\sum_{\kappa \in \Z_+^d(k)} (g_i \mid_D)
(\mathfrak E_\kappa^{d,l -\e,m,D,U_i} f) \mid_D\|_{L_p(D)} = \\
\| \sum_{i =1}^{\mathcal I} ((g_i \mid_D) f) -\sum_{i =1}^{\mathcal I}
(g_i \mid_D) (\sum_{\kappa \in \Z_+^d(k)}
(\mathfrak E_\kappa^{d,l -\e,m,D,U_i} f) \mid_D)\|_{L_p(D)} = \\
\| \sum_{i =1}^{\mathcal I} (g_i \mid_D) (f -\sum_{\kappa \in \Z_+^d(k)}
(\mathfrak E_\kappa^{d,l -\e,m,D,U_i} f) \mid_D)\|_{L_p(D)} \le \\
\sum_{i =1}^{\mathcal I} \| (g_i \mid_D) (f -\sum_{\kappa \in \Z_+^d(k)}
(\mathfrak E_\kappa^{d,l -\e,m,D,U_i} f) \mid_D)\|_{L_p(D)} = \\
\sum_{i =1}^{\mathcal I} \| ((g_i \mid_D) (f -\sum_{\kappa \in \Z_+^d(k)}
(\mathfrak E_\kappa^{d,l -\e,m,D,U_i} f) \mid_D)) \mid_{U_i}\|_{L_p(U_i)} = \\
\sum_{i =1}^{\mathcal I} \| (g_i \mid_{U_i}) ((f \mid_{U_i}) -\sum_{\kappa \in \Z_+^d(k)}
(\mathfrak E_\kappa^{d,l -\e,m,D,U_i} f) \mid_{U_i})\|_{L_p(U_i)} \le \\
\sum_{i =1}^{\mathcal I} \| g_i \mid_{U_i} \|_{L_\infty(U_i)}
\| (f \mid_{U_i}) -\sum_{\kappa \in \Z_+^d(k)}
(\mathfrak E_\kappa^{d,l -\e,m,D,U_i} f) \mid_{U_i}\|_{L_p(U_i)} \to 0\\
\text{ при }\mn(k) \to \infty,
\end{multline*}
т.е. в $ L_p(D) $ верно равенство
\begin{equation*} \tag{2.2.41}
f = \sum_{\kappa \in \Z_+^d} (\mathfrak E_{\kappa}^{d, l -\e,m,D} f) \mid_D,
\end{equation*}
тем самым, соблюдается (2.2.40) при $ q = p, \lambda =0. $

Далее, в условиях п. 3) теоремы для $ f \in (S_p^\alpha H)^\prime(D) $ в силу
(2.2.39) справедлива оценка
\begin{multline*} \tag{2.2.42}
\| \D^\lambda ((\mathfrak E_{\kappa}^{d,l -\e,m,D} f) \mid_D) \|_{L_q(D)} =
\| (\D^\lambda (\mathfrak E_{\kappa}^{d,l -\e,m,D} f)) \mid_D \|_{L_q(D)} \le\\
\| \D^\lambda \mathfrak E_{\kappa}^{d,l -\e,m,D} f \|_{L_q(\R^d)}
\le c_7 2^{(\kappa, \lambda +(p^{-1} -q^{-1})_+ \e)}
\Omega^{\prime l \chi_{\s(\kappa)}}(f, (c_8 2^{-\kappa})^{\s(\kappa)})_{L_p(D)} \le \\
c_{10} 2^{-(\kappa, \alpha -\lambda -(p^{-1} -q^{-1})_+ \e)}
\| f \|_{(S_p^\alpha H)^\prime(D)},
\end{multline*}
а из (2.2.42) с учётом (2.1.68) следует, что ряд
\begin{equation*}
\sum_{ \kappa \in \Z_+^d} \| \D^\lambda ((\mathfrak E_{\kappa}^{d,l -\e,m,D} f) \mid_D) \|_{L_q(D)}
\end{equation*}
сходится, а, значит, в $ L_q(D) $ сходится ряд
\begin{equation*}
\sum_{ \kappa \in \Z_+^d} \D^\lambda ((\mathfrak E_{\kappa}^{d,l -\e,m,D} f) \mid_D).
\end{equation*}
Отсюда и из (2.2.41) вытекает (см. доказательство предложения 2.1.5), что
в $ L_q(D) $ имеет место равенство (2.2.40). $ \square $

В связи с применением леммы 2.2.6 приведём более простые условия на область $ D $ в случае её ограниченности.

Лемма 2.2.8

Пусть $ K $ -- компактное множество в $ \R^d, $ а $ \{\mathcal U_i \subset \R^d,
i \in \mathfrak I\} $ -- открытое покрытие $ K. $ Тогда для каждого
$ i \in \mathfrak I $ существуют ограниченное открытое множество
$ \mathcal V_i \subset \R^d $ и число $ \bm \delta_i > 0 $ такие, что для
$ x \in \mathcal V_i $ справедливо неравенство
\begin{equation*} \tag{2.2.43}
\inf_{y \in \R^d \setminus \mathcal U_i} \| x -y\| > \bm \delta_i, \\
\overline{\mathcal V}_i \subset \mathcal U_i;
\end{equation*}
и $ \{ \mathcal V_i, \mathcal U_j, j \in \mathfrak I \setminus \{i\}\} $
образует открытое покрытие $ K. $

Доказательство.

В условиях леммы при $ i \in \mathfrak I $ определим замкнутое (а, значит,
компактное) подмножество $ K_i \subset K, $ полагая
$$
K_i = K \setminus (\cup_{j \in \mathfrak I \setminus \{i\}} \mathcal U_j)
\subset (\cup_{j \in \mathfrak I} \mathcal U_j) \setminus
(\cup_{j \in \mathfrak I \setminus \{i\}} \mathcal U_j) \subset \mathcal U_i.
$$
Рассмотрим непрерывную в $ \R^d $ функцию
$$
\bm \rho_i(x) = \inf_{y \in \R^d \setminus \mathcal U_i} \| x -y\| =
\min_{y \in \R^d \setminus \mathcal U_i} \| x -y\|.
$$
В силу компактности $ K_i, $ непрерывности $ \bm \rho_i(\cdot) $
и того факта, что
$$
K_i \cap (\R^d \setminus \mathcal U_i) \subset \mathcal U_i \cap (\R^d \setminus
\mathcal U_i) = \emptyset,
$$
заключаем, что существует $ \bm \delta_i > 0 $ такое, что для $ x \in K_i $
выполняется неравенство $ \bm \rho_i(x) > \bm \delta_i. $
Фиксируя $ X^0 \in \R^d, \Delta \in \R_+^d, $ для которых $ K \subset
(X^0 +\Delta I^d), $ возьмём множество
$$
\mathcal V_I = \{x \in (X^0 +\Delta I^d): \bm \rho_i(x) > \bm \delta_i\}.
$$
Тогда ввиду непрерывности функции $ \rho_i(\cdot) $ ограниченное множество
$ \mathcal V_i $ -- открыто и соблюдается соотношение
$$
\overline{\mathcal V}_i \subset \{x \in \R^d: \bm \rho_i(x) \ge \bm \delta_i\}
\subset \mathcal U_i,
$$
ибо $ \bm \rho_i(x) = 0 $ для $ x \in \R^d \setminus \mathcal U_i,$
т.е. выполняется (2.2.43).

Замечая, что $ K_i \subset \mathcal V_i, $ для $ x \in K $ имеем, либо $ x \in
(\cup_{j \in \mathfrak I \setminus \{i\}} \mathcal U_j), $ либо $ x \in
K \setminus (\cup_{j \in \mathfrak I \setminus \{i\}} \mathcal U_j) = K_i
\subset \mathcal V_i, $ т.е.
$$
K \subset (\cup_{j \in \mathfrak I \setminus \{i\}} \mathcal U_j) \cup
\mathcal V_i. \square
$$

Лемма 2.2.9

При $ d \in \N $ пусть $ D $ -- ограниченная область в $ \R^d, $ а
система открытых множеств $ \{ \mathcal U_i \subset \R^d, i =1,\ldots,\mathcal I\} $
образует покрытие её замыкания, т.е.
$ \overline D \subset \cup_{i =1}^{\mathcal I} \mathcal U_i, $ и
$ U_I = \mathcal U_i \cap D, i =1,\ldots,\mathcal I. $ Тогда существует
$ \bm \delta > 0 $ такое, что набор открытых множеств
$$
V_i = \{ x \in U_i: \rho_i(x) = \inf_{ y \in D \setminus U_i} \| x -y\| >
\bm \delta\}, i =1,\ldots,\mathcal I,
$$
образует покрытие области $ D. $

Доказательство.

Сначала построим набор ограниченных открытых множеств $ \{ \mathcal V_i,
i =1,\ldots,\mathcal I\} $ и набор чисел $ \{ \bm \delta_i > 0, i =1,\ldots, \mathcal I\}, $
для которых при $ i =1,\ldots, \mathcal I $ соблюдается (2.2.43) и
$ \{\mathcal V_i, i =1,\ldots, \mathcal I\} $ образует открытое покрытие
$ \overline D. $

Для этого достаточно последовательно при $ i =1,\ldots, \mathcal I $
применить лемму 2.2.8 к открытому покрытию
$$
\{\mathcal V_1,\ldots,\mathcal V_{i -1}, \mathcal U_i,\ldots,\mathcal U_{\mathcal I}\}
$$
компактного множества $ K = \overline D $ и его элементу $ \mathcal U_i, $ и в
результате применения которой получить ограниченное открытое множество
$ \mathcal V_i, $ число $ \bm \delta_i > 0, $ удовлетворяющие (2.2.43),
и открытое покрытие
$$
\{\mathcal V_1,\ldots,\mathcal V_i, \mathcal U_{i +1}, \ldots,\mathcal U_{\mathcal I}\}
$$
множества $ \overline D. $

Теперь фиксируем $ \bm \delta > 0 $ такое, что $ \bm \delta <
\min_{i =1,\ldots,\mathcal I} \bm \delta_i, $ и рассмотрим множества
$$
V_i = \{ x \in U_i: \rho_i(x) = \inf_{ y \in D \setminus U_i} \| x -y\| >
\bm \delta\}, i =1,\ldots,\mathcal I.
$$
Тогда, пользуясь тем, что
$$
D \subset \overline D \subset \cup_{i =1}^{\mathcal I} \mathcal V_i,
$$
для $ x \in D $ возьмём $ i \in \{1,\ldots,\mathcal I\}, $ для которого
$ x \in \mathcal V_i. $ При этом ввиду (2.2.43) имеем также, что $ x \in \mathcal U_i, $
т.е. $ x \in D \cap \mathcal U_i = U_i. $
Далее, заметим, что $ (D \setminus U_i) \subset (\R^d \setminus \mathcal U_i),
i =1,\ldots,\mathcal I. $
В самом деле, для $ y \in D \setminus U_i $ выполняются соотношения
$ y \in \R^d $ и $ y \notin \mathcal U_i, $
поскольку, если $ y \in \mathcal U_i, $  то $ y \in D \cap \mathcal U_i = U_i, $
что неверно.

Учитывая отмеченные обстоятельства, а также (2.2.43), получаем, что для
$ x \in D \cap \mathcal V_i $ имеет место неравенство
\begin{equation*}
\inf_{y \in D \setminus U_i} \| x -y\| \ge
\inf_{y \in \R^d \setminus \mathcal U_i} \| x -y\| > \bm \delta_i > \bm \delta,
\end{equation*}
и, значит, $ x \in V_i. $ Таким образом,
$$
D \subset \cup_{i =1}^{\mathcal I} V_i. \square
$$
\bigskip

2.3. В этом пункте доказывается основной результат работы теорема 2.3.1.

Теорема 2.3.1

Пусть $ d \in \N, D \subset \R^d $ -- область, удовлетворяющая
условиям теоремы 2.2.7, и $ \alpha \in \R_+^d, 1 \le p < \infty, 1 \le \theta
\le \infty. $ Тогда существует непрерывное линейное отображение
$ \mathcal E^{d,\alpha,p,\theta,D}: (S_{p,\theta}^\alpha B)^\prime(D)
\mapsto L_p(\R^d), $ обладающее следующими свойствами:

1) для $ f \in (S_{p,\theta}^\alpha B)^\prime(D) $ выполняется равенство
\begin{equation*} \tag{2.3.1}
f = (\mathcal E^{d,\alpha,p,\theta,D} f) \mid_{D},
\end{equation*}

2) существует константа $ c_1(d,\alpha,p,\theta,D) >0 $ такая, что
при $ l = l(\alpha), \lambda \in \Z_+^d: \lambda < \alpha, $ для
$ f \in (S_{p,\theta}^\alpha B)^\prime(D) $ и любого множества
$ J \subset \{1,\ldots,d\} $ верно неравенство
\begin{multline*} \tag{2.3.2}
\biggl(\int_{(\R_+^d)^J} (t^J)^{-\e^J -\theta (\alpha^J -\lambda^J)}
(\Omega^{(l -\lambda) \chi_J}(\D^\lambda (\mathcal E^{d,\alpha,p,\theta,D} f),
t^J)_{L_p(\R^d)})^\theta dt^J\biggr)^{1/\theta} \le \\
c_1 \| f\|_{(S_{p,\theta}^\alpha B)^\prime(D)}, \ \theta \ne \infty; \\
\sup_{t^J \in (\R_+^d)^J} (t^J)^{-(\alpha^J -\lambda^J)}
\Omega^{(l -\lambda) \chi_J}(\D^\lambda (\mathcal E^{d,\alpha,p,\theta,D} f),
t^J)_{L_p(\R^d)} \le \\
c_1 \| f\|_{(S_{p,\theta}^\alpha B)^\prime(D)}, \ \theta = \infty.
\end{multline*}

Доказательство.

Проведём доказательство при $ \theta \ne \infty. $
Сначала отметим, что ориентируясь на выкладку перед (1.1.9), так же, как
(2.2.7) из [14], выводится используемое ниже соотношение
\begin{multline*} \tag{2.3.3}
\sum_{\kappa \in \Z_+^d \setminus \{0\}} (2^{(\kappa, \alpha)}
\Omega^{\prime l \chi_{\s(\kappa)}}(f,
c (2^{-\kappa})^{\s(\kappa)})_{L_p(D)})^\theta \le \\
c_2(d,\alpha,p,\theta,c) \sum_{ J \subset \Nu_{1,d}^1: J \ne \emptyset}
\int_{(\R_+^d)^J} (t^J)^{-\e^J -\theta \alpha^J} (\Omega^{\prime l \chi_J}(f,
t^J)_{L_p(D)})^\theta dt^J, \\
\alpha \in \R_+^d, 1 \le p < \infty, 1 \le \theta < \infty, l = l(\alpha),\\
f \in (S_{p, \theta}^\alpha B)^\prime(D), c \in \R_+, D 
\text{ -- область в } \R^d.
\end{multline*}

Теперь в условиях теоремы для построения оператора
$ \mathcal E^{d,\alpha,p,\theta,D} $ с требуемыми свойствами фиксируем $ m \in
\N^d: l = l(\alpha) \in \Z_+^d(m), $ и опираясь на теорему 2.2.7, заметим, что
для $ f \in (S_{p, \theta}^\alpha B)^\prime(D) $ при $ \lambda \in \Z_+^d:
\lambda < \alpha, $ ряд
\begin{equation*}
\sum_{\kappa \in \Z_+^d} \| \D^\lambda (\mathfrak E_{\kappa}^{d,l -\e,m,D} f)\|_{L_p(\R^d)}
\end{equation*}
сходится, ибо в силу (2.2.42), (1.1.10) 
при $ \kappa \in \Z_+^d \setminus \{0\} $ имеет место неравенство
\begin{equation*} \tag{2.3.4}
\| \D^\lambda \mathfrak E_{\kappa}^{d,l -\e,m,D} f \|_{L_p(\R^d)}
\le c_3 2^{-(\kappa, \alpha -\lambda)} \| f \|_{(S_p^\alpha H)^\prime(D)} \le \\
c_4 2^{-(\kappa, \alpha -\lambda)} \| f \|_{(S_{p,\theta}^\alpha B)^\prime(D)}.
\end{equation*}
Учитывая замечание после леммы 1.3.1, видим, что для $ f \in
(S_{p,\theta}^\alpha B)^\prime(D) $ при $ \lambda \in \Z_+^d: \lambda < \alpha, $
ряд $ \sum_{ \kappa \in \Z_+^d} \D^\lambda (\mathfrak E_{\kappa}^{d,l -\e,m,D}(f)) $
сходится в $ L_p(\R^d). $
Определим отображение $ \mathcal E^{d,\alpha,p,\theta,D}: (S_{p,\theta}^\alpha B)^\prime(D)
\mapsto L_p(\R^d), $ полагая для $ f \in (S_{p,\theta}^\alpha B)^\prime(D) $
значение
\begin{equation*} \tag{2.3.5}
\mathcal E^{d,\alpha,p,\theta,D} f = \sum_{ \kappa \in \Z_+^d}
\mathfrak E_{\kappa}^{d,l -\e,m,D} f.
\end{equation*}

Исходя из (2.3.5), с учётом сказанного выше заметим, что для
$ f \in (S_{p,\theta}^\alpha B)^\prime(D) $
при $ \lambda \in \Z_+^d: \lambda < \alpha, $ в $ L_p(\R^d) $ имеет
место равенство
\begin{equation*} \tag{2.3.6}
\D^\lambda (\mathcal E^{d,\alpha,p,\theta,D} f) = \sum_{ \kappa \in \Z_+^d}
\D^\lambda (\mathfrak E_{\kappa}^{d,l -\e,m,D} f).
\end{equation*}
Ясно, что отображение $ \mathcal E^{d,\alpha,p,\theta,D} $-- линейно, а
вследствие (2.3.5), (2.2.40) (см. (1.1.9)) в $ L_p(D) $ справедливо
равенство
\begin{equation*}
(\mathcal E^{d,\alpha,p,\theta,D} f) \mid_{D} = \sum_{ \kappa \in \Z_+^d}
(\mathfrak E_{\kappa}^{d,l -\e,m,D} f) \mid_{D} = f, \
\end{equation*}
т.е. соблюдается (2.3.1).

При проверке справедливости п. 2 теоремы, прежде всего, отметим, что
для $ f \in (S_{p,\theta}^\alpha B)^\prime(D) $
при $ \lambda \in \Z_+^d: \lambda < \alpha, $ ввиду (2.3.6), (2.2.38),
(2.3.4) выполняется неравенство
\begin{multline*} \tag{2.3.7}
\| \D^\lambda (\mathcal E^{d,\alpha,p,\theta,D} f) \|_{L_p(\R^d)} =
\| \sum_{ \kappa \in \Z_+^d }
\D^\lambda (\mathfrak E_{\kappa}^{d,l -\e,m,D} f) \|_{L_p(\R^d)} \le \\
\sum_{ \kappa \in \Z_+^d}
\| \D^\lambda (\mathfrak E_{\kappa}^{d,l -\e,m,D} f) \|_{L_p(\R^d)} = \\
\| \D^\lambda (\mathfrak E_0^{d,l -\e,m,D} f) \|_{L_p(\R^d)} +
\sum_{ \kappa \in \Z_+^d \setminus \{0\}}
\| \D^\lambda (\mathfrak E_{\kappa}^{d,l -\e,m,D} f) \|_{L_p(\R^d)} \le \\
c_5 \| f \|_{L_p(D)} +
\sum_{ \kappa \in \Z_+^d \setminus \{0\}}
c_4 2^{-(\kappa, \alpha -\lambda)} \| f \|_{(S_{p,\theta}^\alpha B)^\prime(D)} \le \\
c_6 \| f \|_{(S_{p, \theta}^\alpha B)^\prime(D)}.
\end{multline*}

Далее, при $ \lambda \in \Z_+^d: \lambda < \alpha, $ для функции $ F =
\mathcal E^{d,\alpha,p,\theta,D} f, \ f \in (S_{p,\theta}^\alpha B)^\prime(D), $
и любого непустого множества $ J \subset \{1,\ldots,d\} $ оценим \\
$ (\int_{(\R_+^d)^J} (t^J)^{-\e^J -\theta (\alpha^J -\lambda^J)}
(\Omega^{(l -\lambda) \chi_J}(\D^\lambda F,
t^J)_{L_p(\R^d)})^\theta dt^J)^{1/\theta}. $ Поскольку
для любого непустого множества $ J \subset \{1,\ldots,d\} $
справедливо представление $ (\R_+^d)^J = \cup_{ \J \subset J}
(I^d)^\J \times ([1, \infty)^d)^{J \setminus \J}, $ причём
множества в правой части последнего равенства попарно не пересекаются, то
\begin{multline*} \tag{2.3.8}
\int_{(\R_+^d)^J} (t^J)^{-\e^J -\theta (\alpha^J -\lambda^J)}
(\Omega^{(l -\lambda) \chi_J}(\D^\lambda F, t^J)_{L_p(\R^d)})^\theta dt^J = \\
\sum_{ \J \subset J} \int_{ (I^d)^{\J} \times ([1, \infty)^d)^{J
\setminus \J}} (t^J)^{-\e^J -\theta (\alpha^J -\lambda^J)}
(\Omega^{(l -\lambda) \chi_J}(\D^\lambda F, t^J)_{L_p(\R^d)})^\theta dt^J.
\end{multline*}

Учитывая, что для $  J \subset \{1,\ldots,d\}: J \ne \emptyset, \J
\subset J, $ при $ t^J \in (\R_+^d)^J $ выполняется неравенство
\begin{multline*}
\Omega^{(l -\lambda) \chi_J} (\D^\lambda F, t^J)_{L_p(\R^d)} =
\supvrai_{ \{ h \in \R^d: h^J \in t^J (B^d)^J \}} \| \Delta_h^{(l -\lambda) \chi_J}
\D^\lambda F \|_{L_p(\R^d)} = \\
\supvrai_{ \{ h \in \R^d: h^J \in t^J (B^d)^J \}}
\| \Delta_h^{(l -\lambda) \chi_{J \setminus \J}} (\Delta_h^{(l -\lambda) \chi_{\J}}
\D^\lambda F) \|_{L_p(\R^d)} \le \\
\supvrai_{ \{ h \in \R^d: h^J \in t^J (B^d)^J \}} \!c_7
\| \Delta_h^{(l -\lambda) \chi_{\J}} \D^\lambda F \|_{L_p(\R^d)}\! =\\ 
c_7\!\supvrai_{ \{ h \in \R^d: h^J \in t^J (B^d)^J \}} \!
\| \Delta_h^{(l -\lambda) \chi_{\J}}
\D^\lambda F \|_{L_p(\R^d)} = \\
 c_7 \supvrai_{ \{ h \in \R^d: h^{\J} \in t^{\J} (B^d)^{\J} \}}
\| \Delta_h^{(l -\lambda) \chi_{\J}} \D^\lambda F \|_{L_p(\R^d)} = c_7
\Omega^{(l -\lambda) \chi_{\J}} (\D^\lambda F, t^{\J})_{L_p(\R^d)},
\end{multline*}
находим, что
\begin{multline*}
\int_{ (I^d)^{\J} \times ([1, \infty)^d)^{J \setminus \J}}
(t^J)^{-\e^J -\theta (\alpha^J -\lambda^J)}
(\Omega^{(l -\lambda) \chi_J}(\D^\lambda F, t^J)_{L_p(\R^d)})^\theta dt^J \le\\
c_7^\theta \int_{ (I^d)^{\J} \times ([1, \infty)^d)^{J \setminus \J}}
(t^J)^{-\e^J -\theta (\alpha^J -\lambda^J)}
(\Omega^{(l -\lambda) \chi_{\J}}(\D^\lambda F,
t^{\J})_{L_p(\R^d)})^\theta dt^J = \\
c_7^\theta \int_{ (I^d)^{\J} \times ([1, \infty)^d)^{J \setminus \J}}
(t^{\J}\!)^{\!-\! \e^{\J} \!-\! \theta (\alpha^{\J} -\lambda^{\J})} \!
(t^{J \setminus \J}\!)^{\!-\! \e^{J \setminus \J} \!-\! \theta
(\alpha^{J \setminus \J} -\lambda^{J \setminus \J})}\\
\times (\Omega^{(l -\lambda) \chi_{\J}}(\D^\lambda F,
t^{\J}\!)_{L_p(\R^d)})^\theta dt^{\J} dt^{J \setminus \J}\! = \\
c_7^\theta ( \int_{ ([1, \infty)^d)^{J \setminus \J}} \!
(t^{J \setminus \J}\!)^{\!-\! \e^{J \setminus \J} \!-\! \theta
(\alpha^{J \setminus \J} -\lambda^{J \setminus \J})} \!dt^{J \setminus \J}\!) \!\int_{ (I^d)^{\J} }
\!(t^{\J}\!)^{\!-\! \e^{\J} \!-\! \theta (\alpha^{\J} -\lambda^{\J})}\\
\times (\Omega^{(l -\lambda) \chi_{\J}}(\D^\lambda F, t^{\J}\!)_{L_p(\R^d)})^\theta
\!dt^{\J} \!\le\\
c_8 \int_{ (I^d)^{\J} } (t^{\J})^{-\e^{\J} -\theta (\alpha^{\J} -\lambda^{\J})}
(\Omega^{(l -\lambda) \chi_{\J}}(\D^\lambda F,
t^{\J})_{L_p(\R^d)})^\theta dt^{\J}.
\end{multline*}

Подставляя последнюю оценку в (2.3.8), получаем, что
\begin{multline*} \tag{2.3.9}
\int_{(\R_+^d)^J} (t^J)^{-\e^J -\theta (\alpha^J -\lambda^J)}
(\Omega^{(l -\lambda) \chi_J}(\D^\lambda F, t^J)_{L_p(\R^d)})^\theta dt^J \le\\
c_8 \sum_{ \J \subset J} \int_{ (I^d)^{\J} } (t^{\J})^{-\e^{\J}
-\theta (\alpha^{\J} -\lambda^{\J})}
(\Omega^{(l -\lambda) \chi_{\J}}(\D^\lambda F,
t^{\J})_{L_p(\R^d)})^\theta dt^{\J} = \\
c_8 ( \| \D^\lambda F \|_{L_p(\R^d)}^\theta +
\sum_{ \J \subset J: \J \ne \emptyset} \int_{ (I^d)^{\J} }
(t^{\J})^{-\e^{\J} -\theta (\alpha^{\J} -\lambda^{\J})} \\
\times (\Omega^{(l -\lambda) \chi_{\J}}(\D^\lambda F,
t^{\J})_{L_p(\R^d)})^\theta dt^{\J}).
\end{multline*}

Таким образом, приходим к необходимости оценки $ \int_{ (I^d)^{\J} }
(t^{\J})^{-\e^{\J} -\theta (\alpha^{\J} -\lambda^{\J})}$
$(\Omega^{(l -\lambda) \chi_{\J}}(\D^\lambda F,
t^{\J})_{L_p(\R^d)})^\theta dt^{\J} $ для $ \J \subset \{1, \ldots,d\}:
\J \ne \emptyset. $  Проведём эту оценку.

Пользуясь тем, что для любого непустого множества $ J \subset \{1,\ldots, d\} $
имеет место представление $ (I^d)^J = (\cup_{k^J \in (\N^d)^J}
(2^{-k^J} +2^{-k^J} (I^d)^J)) \cup M^J, \ $
причём множества в правой части последнего равенства попарно не
пересекаются и $ \mes M^J =0, $  получаем, что для любого непустого
множества $ J \subset \{1,\ldots,d\}$ имеет место неравенство
\begin{multline*} \tag{2.3.10}
\int_{ (I^d)^J } (t^J)^{-\e^J -\theta (\alpha^J -\lambda^J)}
(\Omega^{(l -\lambda) \chi_J}(\D^\lambda F, t^J)_{L_p(\R^d)})^\theta dt^J  = \\
\sum_{ k^J \in (\N^d)^J} \!\int_{(2^{-k^J} \!+2^{-k^J} (I^d)^J)} \!
(t^J)^{-\e^J -\theta (\alpha^J -\lambda^J)}
(\Omega^{(l -\lambda) \chi_J}(\D^\lambda F, t^J)_{L_p(\R^d)})^\theta dt^J \le \\
\sum_{ k^J \in (\N^d)^J} \!\int_{(2^{-k^J} \!+2^{-k^J} (I^d)^J)} \!
(2^{-k^J}\!)^{-\e^J -\theta (\alpha^J -\lambda^J)}
(\Omega^{(l -\lambda) \chi_J}(\D^\lambda F,
2 \!\cdot\! 2^{-k^J}\!)_{L_p(\R^d)})^\theta dt^J\! = \\
\sum_{ k^J \in (\N^d)^J} 2^{\theta (k^J, \alpha^J -\lambda^J)}
(\Omega^{(l -\lambda) \chi_J}(\D^\lambda F, 2 \cdot 2^{-k^J})_{L_p(\R^d)})^\theta.
\end{multline*}

Оценивая $ \Omega^{(l -\lambda) \chi_J} (\D^\lambda F, 2 \cdot 2^{-k^J})_{L_p(\R^d)} $
при $ k^J \in (\N^d)^J, J \subset \{1,\ldots,d\}: J \ne \emptyset, $
для $ h \in \R^d: h^J \in 2 \cdot 2^{-k^J} (B^d)^J $ с учётом (2.3.6),
(1.1.8) и (2.2.39), (2.2.38) имеем
\begin{multline*}
\| \Delta_h^{(l -\lambda) \chi_J} \D^\lambda F \|_{L_p(\R^d)} =
\| \Delta_h^{(l -\lambda) \chi_J} (\sum_{ \kappa \in \Z_+^d}
\D^\lambda (\mathfrak E_{\kappa}^{d,l -\e,m,D} f)) \|_{L_p(\R^d)} = \\
\| \sum_{ \kappa \in \Z_+^d} \Delta_h^{(l -\lambda) \chi_J}
(\D^\lambda (\mathfrak E_{\kappa}^{d,l -\e,m,D} f)) \|_{L_p(\R^d)} \le \\
\sum_{ \kappa \in \Z_+^d} \| \Delta_h^{(l -\lambda) \chi_J}
(\D^\lambda (\mathfrak E_{\kappa}^{d,l -\e,m,D} f)) \|_{L_p(\R^d)} = \\
\sum_{ \J \subset J} \sum_{\substack{ \kappa \in \Z_+^d:
\kappa^{\J} \le k^{\J},\\ \kappa^{J \setminus \J} > k^{J \setminus \J}}}
\| \Delta_h^{(l -\lambda) \chi_J}
(\D^\lambda (\mathfrak E_{\kappa}^{d,l -\e,m,D} f)) \|_{L_p(\R^d)} \le \\
\sum_{ \J \subset J} \sum_{\substack{ \kappa \in \Z_+^d:
\kappa^{\J} \le k^{\J},\\ \kappa^{J \setminus \J} > k^{J \setminus \J}}}
c_7 \| \Delta_h^{(l -\lambda) \chi_{\J}}
(\D^\lambda (\mathfrak E_{\kappa}^{d,l -\e,m,D} f)) \|_{L_p(\R^d)} \le \\
\sum_{ \J \subset J} \sum_{\substack{ \kappa \in \Z_+^d:
\kappa^{\J} \le k^{\J},\\ \kappa^{J \setminus \J} > k^{J \setminus \J}} }
c_7 (\prod_{ j \in \J} | h_j |^{ l_j -\lambda_j}) \| \D^{(l -\lambda) \chi_{\J}}
(\D^\lambda (\mathfrak E_{\kappa}^{d,l -\e,m,D} f)) \|_{L_p(\R^d)} \le \\
\sum_{ \J \subset J} \sum_{\substack{ \kappa \in \Z_+^d:
\kappa^{\J} \le k^{\J},\\ \kappa^{J \setminus \J} > k^{J \setminus \J}}}
c_7 (\prod_{ j \in \J} (2^{-k_j +1})^{ l_j -\lambda_j})
\| \D^{(l -\lambda) \chi_{\J} +\lambda}
(\mathfrak E_{\kappa}^{d,l -\e,m,D} f) \|_{L_p(\R^d)} \le \\
\sum_{ \J \subset J} \sum_{\substack{ \kappa \in \Z_+^d:
\kappa^{\J} \le k^{\J},\\ \kappa^{J \setminus \J} > k^{J \setminus \J}} }
c_9 (\prod_{ j \in \J} 2^{-k_j (l_j -\lambda_j)})
2^{ ( \kappa, (l -\lambda) \chi_{\J} +\lambda)}
\Omega^{\prime l \chi_{\s(\kappa)}}(f, c_{10} (2^{-\kappa})^{\s(\kappa)})_{L_p(D)} = \\
c_9 \sum_{ \J \subset J} 2^{ -( k^{\J}, l^{\J} -\lambda^{\J})}
\sum_{\substack{ \kappa \in \Z_+^d: \kappa^{\J} \le k^{\J},\\
\kappa^{J \setminus \J} > k^{J \setminus \J}} }
2^{ ( \kappa^{\J}, l^{\J} -\lambda^{\J})} 2^{(\kappa, \lambda)}
\Omega^{\prime l \chi_{\s(\kappa)}}(f, c_{10} (2^{-\kappa})^{\s(\kappa)})_{L_p(D)}.
\end{multline*}

Отсюда получаем, что для $ k^J \in (\N^d)^J, J \subset
\{1,\ldots,d\}: J \ne \emptyset, $ выполняется неравенство
\begin{multline*}
\Omega^{(l -\lambda) \chi_J} (\D^\lambda F, 2 \cdot 2^{-k^J})_{L_p(\R^d)} \le \\
c_9 \sum_{ \J \subset J} 2^{ -( k^{\J}, l^{\J} -\lambda^{\J})}
\sum_{\substack{ \kappa \in \Z_+^d: \kappa^{\J} \le k^{\J},\\
\kappa^{J \setminus \J} > k^{J \setminus \J}}}
2^{ ( \kappa^{\J}, l^{\J} -\lambda^{\J})} 2^{(\kappa, \lambda)}
\Omega^{\prime l \chi_{\s(\kappa)}}(f, c_{10} (2^{-\kappa})^{\s(\kappa)})_{L_p(D)}.
\end{multline*}

Подставляя эту оценку в (2.3.10) и применяя неравенство Гёльдера,
для $  J \subset \{1,\ldots,d\}: J \ne \emptyset, $ выводим
\begin{multline*} \tag{2.3.11}
\int_{ (I^d)^J } (t^J)^{-\e^J -\theta (\alpha^J -\lambda^J)}
(\Omega^{(l -\lambda) \chi_J}(\D^\lambda F, t^J)_{L_p(\R^d)})^\theta dt^J \le
\sum_{ k^J \in (\N^d)^J} 2^{\theta (k^J, \alpha^J -\lambda^J)} \\
\times \biggl( c_9 \sum_{ \J \subset J} 2^{ -( k^{\J}, l^{\J} -\lambda^{\J})}
\sum_{\substack{ \kappa \in \Z_+^d: \kappa^{\J} \le k^{\J}, \\
\kappa^{J \setminus \J} > k^{J \setminus \J}} }
2^{ ( \kappa^{\J}, l^{\J} -\lambda^{\J})} 2^{(\kappa, \lambda)}
\Omega^{\prime l \chi_{\s(\kappa)}}(f, c_{10} (2^{-\kappa})^{\s(\kappa)})_{L_p(D)} \biggr)^\theta \le \\
c_{11} \sum_{ k^J \in (\N^d)^J} 2^{\theta (k^J, \alpha^J -\lambda^J)}
\sum_{ \J \subset J} \biggl(2^{ -( k^{\J}, l^{\J} -\lambda^{\J})}
\sum_{ \kappa \in \Z_+^d: \kappa^{\J} \le k^{\J},
\kappa^{J \setminus \J} > k^{J \setminus \J}}
2^{ ( \kappa^{\J}, l^{\J} -\lambda^{\J})} \\
\times 2^{(\kappa, \lambda)}
\Omega^{\prime l \chi_{\s(\kappa)}}(f, c_{10} (2^{-\kappa})^{\s(\kappa)})_{L_p(D)} \biggr)^\theta = \\
c_{11} \sum_{ k^J \in (\N^d)^J} \sum_{ \J \subset J}
2^{\theta (k^J, \alpha^J -\lambda^J) -\theta ( k^{\J}, l^{\J} -\lambda^{\J})} \\
\times \biggl(\sum_{ \kappa \in \Z_+^d: \kappa^{\J} \le k^{\J},
\kappa^{J \setminus \J} > k^{J \setminus \J}}
2^{ ( \kappa^{\J}, l^{\J} -\lambda^{\J})} 2^{(\kappa, \lambda)}
\Omega^{\prime l \chi_{\s(\kappa)}}(f, c_{10} (2^{-\kappa})^{\s(\kappa)})_{L_p(D)} \biggr)^\theta = \\
c_{11} \sum_{ \J \subset J} \sum_{ k^J \in (\N^d)^J}
2^{\theta (k^J, \alpha^J -\lambda^J) -\theta ( k^{\J}, l^{\J} -\lambda^{\J})} \\
\times \biggl(\sum_{ \kappa \in \Z_+^d: \kappa^{\J} \le k^{\J},
\kappa^{J \setminus \J} > k^{J \setminus \J}}
2^{ ( \kappa^{\J}, l^{\J} -\lambda^{\J})} 2^{(\kappa, \lambda)}
\Omega^{\prime l \chi_{\s(\kappa)}}(f, c_{10} (2^{-\kappa})^{\s(\kappa)})_{L_p(D)} \biggr)^\theta.
\end{multline*}

Для оценки правой части (2.3.11) при $ \lambda \in \Z_+^d: \lambda < \alpha, $
фиксируем $ \epsilon \in \R_+^d, $ для которого соблюдаются условия
$ \alpha -\lambda -\epsilon >0, l -\alpha -\epsilon >0, $ и положим
$ \overline J = \{1, \ldots, d\} \setminus J. $ Используя неравенство
Гёльдера, для $ J \subset \{1,\ldots,d\}: J \ne \emptyset, \J \subset J,
k^J \in (\N^d)^J $ приходим к неравенству
\begin{multline*} \tag{2.3.12}
\biggl(\sum_{ \kappa \in \Z_+^d: \kappa^{\J} \le k^{\J},
\kappa^{J \setminus \J} > k^{J \setminus \J}}
2^{ ( \kappa^{\J}, l^{\J} -\lambda^{\J})} 2^{(\kappa, \lambda)}
\Omega^{\prime l \chi_{\s(\kappa)}}(f,
c_{10} (2^{-\kappa})^{\s(\kappa)})_{L_p(D)}\biggr)^\theta = \\
\biggl(\sum_{ \kappa \in \Z_+^d: \kappa^{\J} \le k^{\J},
\kappa^{J \setminus \J} > k^{J \setminus \J}}
2^{-\!(\kappa, \alpha)} 2^{ (\kappa^{\J} \!, l^{\J} -\lambda^{\J})}
2^{(\kappa, \lambda)} 2^{(\kappa^{J \setminus \J} \!, \epsilon^{J \setminus \J})} \\
\times 2^{-\!(\kappa^{\J} \!, \epsilon^{\J})} 2^{-\!(\kappa^{J \setminus \J} \!,
\epsilon^{J \setminus \J})} 2^{ (\kappa^{\J} \!, \epsilon^{\J})}
2^{(\kappa, \alpha)} \Omega^{\prime l \chi_{\s(\kappa)}}(f,
c_{10} (2^{-\kappa})^{\s(\kappa)})_{L_p(D)} \biggr)^\theta \le \\
\biggl(\sum_{ \kappa \in \Z_+^d: \kappa^{\J} \le k^{\J},
\kappa^{J \setminus \J} > k^{J \setminus \J}}\! 2^{-\theta^\prime (\kappa,
\alpha -\lambda)} 2^{\theta^\prime (\kappa^{J \setminus \J},
\epsilon^{J \setminus \J})}
2^{ \theta^\prime (\kappa^{\J}, l^{\J} -\lambda^{\J} -\epsilon^{\J})}\biggr)^{\theta / \theta^\prime} \\
\times \biggl(\sum_{ \kappa \in \Z_+^d: \kappa^{\J} \le k^{\J},
\kappa^{J \setminus \J} > k^{J \setminus \J}} (2^{-(\kappa^{J \setminus \J},
\epsilon^{J \setminus \J})} 2^{ (\kappa^{\J}, \epsilon^{\J})}
2^{(\kappa, \alpha)} \Omega^{\prime l \chi_{\s(\kappa)}}(f,
c_{10} (2^{-\kappa})^{\s(\kappa)})_{L_p(D)})^\theta \biggr).
\end{multline*}

Оценивая правую часть (2.3.12), для $ J \subset \{1,\ldots,d\}: J \ne \emptyset,
\J \subset J, k^J \in (\N^d)^J $ имеем
\begin{multline*} \tag{2.3.13}
\biggl(\sum_{ \kappa \in \Z_+^d: \kappa^{\J} \le k^{\J},
\kappa^{J \setminus \J} > k^{J \setminus \J}} 2^{-\theta^\prime (\kappa, \alpha -\lambda)}
2^{\theta^\prime (\kappa^{J \setminus \J}, \epsilon^{J \setminus \J})}
2^{ \theta^\prime (\kappa^{\J}, l^{\J} -\lambda^{\J} -\epsilon^{\J})}\biggr)^{\theta / \theta^\prime} = \\
\biggl(\sum_{ \kappa^{\overline J} \in (\Z_+^d)^{\overline J},
\kappa^{\J} \in (\Z_+^d)^{\J}: \kappa^{\J} \le k^{\J},
\kappa^{J \setminus \J} \in (\Z_+^d)^{J \setminus \J}:
\kappa^{J \setminus \J} > k^{J \setminus \J}}
2^{-\theta^\prime (\kappa^{\overline J}, \alpha^{\overline J} -\lambda^{\overline J})}  \\
\times 2^{-\theta^\prime (\kappa^{J \setminus \J}, \alpha^{J \setminus \J} -
\lambda^{J \setminus \J} -\epsilon^{J \setminus \J})}
2^{ \theta^\prime (\kappa^{\J}, l^{\J} -\alpha^{\J} -\epsilon^{\J})} \biggr)^{\theta /
\theta^\prime} =
\biggl(\biggl(\sum_{ \kappa^{\overline J} \in
(\Z_+^d)^{\overline J}} 2^{-\theta^\prime (\kappa^{\overline J},
\alpha^{\overline J} -\lambda^{\overline J})} \biggr)  \\
\times \biggl(\sum_{\kappa^{J \setminus \J} \in (\Z_+^d)^{J \setminus \J}:
\kappa^{J \setminus \J} > k^{J \setminus \J}}
2^{-\theta^\prime (\kappa^{J \setminus \J}, \alpha^{J \setminus \J}
-\lambda^{J \setminus \J} -\epsilon^{J \setminus \J})} \biggr)  \\
\times \biggl(\sum_{\kappa^{\J} \in (\Z_+^d)^{\J}: \kappa^{\J} \le k^{\J}}
2^{ \theta^\prime (\kappa^{\J}, l^{\J} -\alpha^{\J} -\epsilon^{\J})} \biggr)\biggr)^{\theta / \theta^\prime} =
\biggl((\prod_{j \in \overline J}
(\sum_{\kappa_j =0}^\infty
2^{-\theta^\prime \kappa_j (\alpha_j -\lambda_j)}))  \\
\times (\prod_{j \in J \setminus \J} (\sum_{\kappa_j = k_j +1}^\infty
2^{-\theta^\prime \kappa_j (\alpha_j -\lambda_j -\epsilon_j)})) \times
(\prod_{j \in \J} (\sum_{\kappa_j =0}^{ k_j}
2^{ \theta^\prime \kappa_j (l_j -\alpha_j -\epsilon_j)}))\biggr)^{\theta / \theta^\prime} \le \\
\biggl( c_{12} (\prod_{j \in J \setminus \J}
2^{-\theta^\prime k_j (\alpha_j -\lambda_j -\epsilon_j)})
(\prod_{j \in \J} 2^{ \theta^\prime k_j
(l_j -\alpha_j -\epsilon_j)} )\biggr)^{\theta / \theta^\prime} = \\
c_{13} \biggl(2^{-\theta^\prime (k^{J \setminus \J}, \alpha^{J \setminus \J} -
\lambda^{J \setminus \J} -\epsilon^{J \setminus \J})}
2^{ \theta^\prime (k^{\J}, l^{\J} -\alpha^{\J} -\epsilon^{\J})}\biggr)^{\theta / \theta^\prime} = \\
c_{13} 2^{-\theta (k^{J \setminus \J}, \alpha^{J \setminus \J}
-\lambda^{J \setminus \J} -\epsilon^{J \setminus \J})}
2^{ \theta (k^{\J}, l^{\J} -\alpha^{\J} -\epsilon^{\J})}.
\end{multline*}

Соединяя (2.3.11), (2.3.12), (2.3.13), получаем неравенство
\begin{multline*} \tag{2.3.14}
\int_{ (I^d)^J } (t^J)^{-\e^J -\theta (\alpha^J -\lambda^J)}
(\Omega^{(l -\lambda) \chi_J}(\D^\lambda F, t^J)_{L_p(\R^d)})^\theta dt^J \le \\
c_{11} \sum_{ \J \subset J} \sum_{ k^J \in (\N^d)^J}
2^{\theta (k^J, \alpha^J -\lambda^J) -\theta ( k^{\J}, l^{\J} -\lambda^{\J})} \\
\times c_{13} 2^{-\theta (k^{J \setminus \J}, \alpha^{J \setminus \J}
-\lambda^{J \setminus \J} -\epsilon^{J \setminus \J})}
2^{ \theta (k^{\J}, l^{\J} -\alpha^{\J} -\epsilon^{\J})} \\
\biggl(\sum_{ \kappa \in \Z_+^d: \kappa^{\J} \le k^{\J},
\kappa^{J \setminus \J} > k^{J \setminus \J}} (2^{-(\kappa^{J \setminus \J},
\epsilon^{J \setminus \J})} 2^{ (\kappa^{\J}, \epsilon^{\J})}
2^{(\kappa, \alpha)}
\Omega^{\prime l \chi_{\s(\kappa)}}(f, c_{10} (2^{-\kappa})^{\s(\kappa)})_{L_p(D)} )^\theta\biggr) = \\
c_{14} \sum_{ \J \subset J} \sum_{ k^J \in (\N^d)^J}
(2^{\theta (k^{J \setminus \J}, \epsilon^{J \setminus \J})}
2^{ -\theta (k^{\J}, \epsilon^{\J})}) \\
\sum_{ \kappa \in \Z_+^d: \kappa^{\J} \le k^{\J},
\kappa^{J \setminus \J} > k^{J \setminus \J}} \biggl(2^{-(\kappa^{J \setminus \J},
\epsilon^{J \setminus \J})} 2^{ (\kappa^{\J}, \epsilon^{\J})}
2^{(\kappa, \alpha)}
\Omega^{\prime l \chi_{\s(\kappa)}}(f, c_{10} (2^{-\kappa})^{\s(\kappa)})_{L_p(D)} \biggr)^\theta = \\
c_{14} \sum_{ \J \subset J} \sum_{ k^J \in (\N^d)^J}
\sum_{\kappa \in \Z_+^d: \kappa^{\J} \le k^{\J},
\kappa^{J \setminus \J} > k^{J \setminus \J}}
(2^{\theta (k^{J \setminus \J}, \epsilon^{J \setminus \J})}
2^{ -\theta (k^{\J}, \epsilon^{\J})})  \\
\times (2^{-(\kappa^{J \setminus \J}, \epsilon^{J \setminus \J})}
2^{ (\kappa^{\J}, \epsilon^{\J})} 2^{(\kappa, \alpha)}
\Omega^{\prime l \chi_{\s(\kappa)}}(f, c_{10} (2^{-\kappa})^{\s(\kappa)})_{L_p(D)} )^\theta = \\
c_{14} \sum_{ \J \subset J} \sum_{ k^J \in (\N^d)^J, \kappa \in \Z_+^d:
\kappa^{\J} \le k^{\J}, \kappa^{J \setminus \J} > k^{J \setminus \J}}
(2^{\theta (k^{J \setminus \J}, \epsilon^{J \setminus \J})}
2^{ -\theta (k^{\J}, \epsilon^{\J})})  \\
\times (2^{-(\kappa^{J \setminus \J}, \epsilon^{J \setminus \J})}
2^{ (\kappa^{\J}, \epsilon^{\J})} 2^{(\kappa, \alpha)}
\Omega^{\prime l \chi_{\s(\kappa)}}(f, c_{10} (2^{-\kappa})^{\s(\kappa)})_{L_p(D)} )^\theta = \\
c_{14} \sum_{ \J \subset J} \sum_{\kappa \in \Z_+^d} \sum_{ k^J \in (\N^d)^J:
\kappa^{\J} \le k^{\J}, \kappa^{J \setminus \J} > k^{J \setminus \J}}
(2^{\theta (k^{J \setminus \J}, \epsilon^{J \setminus \J})}
2^{ -\theta (k^{\J}, \epsilon^{\J})})  \\
\times (2^{-(\kappa^{J \setminus \J}, \epsilon^{J \setminus \J})}
2^{ (\kappa^{\J}, \epsilon^{\J})} 2^{(\kappa, \alpha)}
\Omega^{\prime l \chi_{\s(\kappa)}}(f, c_{10} (2^{-\kappa})^{\s(\kappa)})_{L_p(D)} )^\theta = \\
c_{14} \sum_{ \J \subset J} \sum_{\kappa \in \Z_+^d}
(2^{-(\kappa^{J \setminus \J}, \epsilon^{J \setminus \J})}
2^{(\kappa^{\J}, \epsilon^{\J})} 2^{(\kappa, \alpha)}
\Omega^{\prime l \chi_{\s(\kappa)}}(f, c_{10} (2^{-\kappa})^{\s(\kappa)})_{L_p(D)} )^\theta \\
\biggl(\sum_{ k^J \in (\N^d)^J: \kappa^{\J} \le k^{\J},
\kappa^{J \setminus \J} > k^{J \setminus \J}}
(2^{\theta (k^{J \setminus \J}, \epsilon^{J \setminus \J})}
2^{ -\theta (k^{\J}, \epsilon^{\J})})\biggr), \\
J \subset \{1,\ldots,d\}: J \ne \emptyset.
\end{multline*}

Оценивая правую часть (2.3.14), для $  J \subset \{1,\ldots,d\}: J
\ne \emptyset, \J \subset J, \kappa \in \Z_+^d $ имеем
\begin{multline*}
\sum_{ k^J \in (\N^d)^J: \kappa^{\J} \le k^{\J},
\kappa^{J \setminus \J} > k^{J \setminus \J}}
(2^{\theta (k^{J \setminus \J}, \epsilon^{J \setminus \J})}
2^{ -\theta (k^{\J}, \epsilon^{\J})}) = \\
\sum_{ k^{\J} \in (\N^d)^{\J}: \kappa^{\J} \le k^{\J},
k^{J \setminus \J} \in (\N^d)^{J \setminus \J}:
\kappa^{J \setminus \J} > k^{J \setminus \J}}
2^{\theta (k^{J \setminus \J}, \epsilon^{J \setminus \J})}
2^{ -\theta (k^{\J}, \epsilon^{\J})} = \\
\biggl(\sum_{ k^{\J} \in (\N^d)^{\J}: \kappa^{\J} \le k^{\J}}
2^{ -\theta (k^{\J}, \epsilon^{\J})}\biggr)
\biggl(\sum_{k^{J \setminus \J} \in (\N^d)^{J \setminus \J}:
\kappa^{J \setminus \J} > k^{J \setminus \J}}
2^{\theta (k^{J \setminus \J}, \epsilon^{J \setminus \J})}\biggr) = \\
\biggl( \prod_{j \in \J} (\sum_{ k_j = \kappa_j}^\infty
2^{ -\theta k_j \epsilon_j})\biggr) \biggl(\prod_{j \in J \setminus \J}
(\sum_{k_j =1}^{\kappa_j -1} 2^{\theta k_j \epsilon_j})\biggr) \le \\
c_{15} (\prod_{j \in \J} 2^{ -\theta \kappa_j \epsilon_j})
(\prod_{j \in J \setminus \J}
2^{\theta \kappa_j \epsilon_j}) =
c_{15} 2^{ -\theta (\kappa^{\J}, \epsilon^{\J})}
2^{\theta(\kappa^{J \setminus \J}, \epsilon^{J \setminus \J})}.
\end{multline*}

Подставляя эту оценку в (2.3.14) и учитывая (2.3.3), получаем, что для
$ J \subset \{1,\ldots,d\}: J \ne \emptyset, $ справедливо неравенство
\begin{multline*} \tag{2.3.15}
\int_{ (I^d)^J } (t^J)^{-\e^J -\theta (\alpha^J -\lambda^J)}
(\Omega^{(l -\lambda) \chi_J}(\D^\lambda F, t^J)_{L_p(\R^d)})^\theta dt^J \le \\
c_{14} \!\sum_{ \J \subset J} \sum_{\kappa \in \Z_+^d}
\!(2^{-(\kappa^{J \setminus \J}, \epsilon^{J \setminus \J})}
2^{(\kappa^{\J}, \epsilon^{\J})} 2^{(\kappa, \alpha)}
\Omega^{\prime l \chi_{\s(\kappa)}}(f,
c_{10} (2^{-\kappa})^{\s(\kappa)})_{L_p(D)})^\theta \\
\times c_{15} 2^{ -\theta (\kappa^{\J}, \epsilon^{\J})}
2^{\theta (\kappa^{J \setminus \J}, \epsilon^{J \setminus \J})} = \\
c_{16} \!\sum_{ \J \subset J} \sum_{\kappa \in \Z_+^d}
\!(2^{(\kappa, \alpha)} \Omega^{\prime l \chi_{\s(\kappa)}}(f,
c_{10} (2^{-\kappa})^{\s(\kappa)})_{L_p(D)} )^\theta \le \\
c_{17} \sum_{\kappa \in \Z_+^d} (2^{(\kappa, \alpha)}
\Omega^{\prime l \chi_{\s(\kappa)}}(f, c_{10} (2^{-\kappa})^{\s(\kappa)})_{L_p(D)})^\theta = \\
c_{17} (\| f \|_{L_p(D)}^\theta +\sum_{\kappa \in \Z_+^d \setminus \{0\}}
(2^{(\kappa, \alpha)} \Omega^{\prime l \chi_{\s(\kappa)}}(f,
c_{10} (2^{-\kappa})^{\s(\kappa)})_{L_p(D)})^\theta) \le \\
c_{18} \biggl(\| f\|_{L_p(D)}^\theta +\sum_{ \J \subset \Nu_{1,d}^1: \J \ne \emptyset}
\int_{(\R_+^d)^{\J}} (t^{\J})^{-\e^{\J} -\theta \alpha^{\J}} (\Omega^{\prime l \chi_{\J}}(f,
t^{\J})_{L_p(D)})^\theta dt^{\J}\biggr) \le \\
c_{19} \| f \|_{(S_{p, \theta}^\alpha B)^\prime(D)}^\theta.
\end{multline*}

Соединяя (2.3.7), (2.3.9), (2.3.15), приходим к (2.3.2), чем завершаем
доказательство теоремы при $ \theta \ne \infty. $ При $ \theta = \infty $
доказательство теоремы проводится по той же схеме с заменой в соответствующих
местах операции суммирования на операцию взятия супремума или максимума.$ \square $

Следствие

В условиях теоремы 2.3.1 при $ \l \in \Z_+^d: \l < \alpha, $ имеет место
включение
\begin{equation*} \tag{2.3.16}
(S_{p, \theta}^\alpha B)^\prime(D) \subset (S_{p, \theta}^\alpha B)^{\l}(D), \
\end{equation*}
и для любой функции $ f \in (S_{p, \theta}^\alpha B)^\prime(D) $
выполняется неравенство
\begin{equation*} \tag{2.3.17}
\| f\|_{(S_{p, \theta}^\alpha B)^{\l}(D)} \le c_1
\| f\|_{(S_{p, \theta}^\alpha B)^\prime(D)}.
\end{equation*}

Для получения (2.3.16), (2.3.17) в условиях теоремы 2.3.1 достаточно
применить (2.3.1), оценку
\begin{equation*}
\| (\mathcal E^{d,\alpha,p,\theta,D} f) \mid_{D} \|_{(S_{p, \theta}^\alpha B)^{\l}(D)} \le
\| \mathcal E^{d,\alpha,p,\theta,D} f\|_{(S_{p, \theta}^\alpha B)^{\l}(\R^d)}
\end{equation*}
и неравенство (2.3.2) при
$ \lambda = \l \chi_J, J \subset \{1,\ldots,d\}. $
Из (1.1.11), (1.1.12) и (2.3.16), (2.3.17) вытекает, что в условиях
теоремы 2.3.1 $ (S_{p, \theta}^\alpha B)^\prime(D) =
(S_{p, \theta}^\alpha B)^{\l}(D), $ и нормы на этих пространствах
$ \| \cdot \|_{(S_{p, \theta}^\alpha B)^\prime(D)},
\| \cdot \|_{(S_{p, \theta}^\alpha B)^{\l}(D)} $ эквивалентны
при $ \l \in \Z_+^d: \l < \alpha. $

Отметим также, что из теоремы 2.3.1 и (1.1.7) следует, что при
$ \lambda \in \Z_+^d: \lambda < \alpha, $ в условиях теоремы 2.3.1 имеет место
включение $ \D^\lambda \mid_{(S_{p, \theta}^\alpha B)^\prime(D)} \in
\mathcal B((S_{p, \theta}^\alpha B)^\prime(D),
(S_{p, \theta}^{\alpha -\lambda} B)^\prime(D)). $

\end{document}